%% file: main.tex
\newif\ifdraft
\draftfalse

\ifdraft
\documentclass[aos]{imsart-no-journal-name}
\else
\documentclass[aos]{imsart-no-journal-name}
\fi

\RequirePackage{amsthm,amsmath,amsfonts,amssymb}
\RequirePackage[numbers]{natbib}
\RequirePackage[colorlinks,citecolor=blue,urlcolor=blue]{hyperref}
\usepackage{lslmath}
\RequirePackage{graphicx}
\usepackage{psfig,epsf}
\usepackage{caption,subcaption}
\usepackage[ruled,vlined]{algorithm2e}
\usepackage{import}

\renewcommand{\new}{\emph}

\ifdraft
\newlength{\offt}
\newlength{\offb}
\newlength{\offh}
\fi

\startlocaldefs

\SetKwInput{KwInput}{Input}
\SetKwInput{KwOutput}{Output}
\SetKwInput{KwStart}{Initialization step}
\SetKwInput{KwMinStep}{Minimization step}
\SetKwInput{KwEnd}{Final step}
\SetKwInput{KwInitializationStep}{Initialize}
\SetKwInput{KwUpdateStep}{Update}
\SetKwInput{KwDynUpdateStep}{Dynamic update}
\SetKwInput{KwProcessStep}{Process}
\SetKwInput{KwReturn}{Return}
\SetKwInput{KwConsensus}{Consensus}

\newcommand{\V}{\mathbb{V}} 

\newcommand{\insymbol}{{\rm in}}
\newcommand{\outsymbol}{{\rm out}}

\newcommand{\fin}{f}
\newcommand{\fout}{g}
\newcommand{\pin}{p}
\newcommand{\pout}{q}

\newcommand{\cin}{c_\insymbol}
\newcommand{\cout}{c_\outsymbol}
\newcommand{\vin}{v^{\rm in}}
\newcommand{\vout}{v^{\rm out}}

\newcommand{\alphamin}{\alpha_{\rm min}}
\newcommand{\alphamax}{\alpha_{\rm max}}

\newcommand{\xon}{{x}_{\rm on}}

\newcommand{\hsigma}{\hat{\sigma}}
\newcommand{\hr}{\hat{r}}
\newcommand{\hR}{\hat{R}}

\newcommand{\hP}{\widehat{P}}
\newcommand{\hQ}{\widehat{Q}}

\newcommand{\hZ}{\hat{Z}}
\newcommand{\hmu}{\hat \mu}
\newcommand{\hnu}{\hat \nu}

\newcommand{\tf}{\tilde f}
\newcommand{\tg}{\tilde g}
\newcommand{\tI}{\tilde I}

\newcommand{\tX}{\tilde{X}}

\newcommand{\Nmin}{N_{\rm min}}
\newcommand{\Nmax}{N_{\rm max}}
\newcommand{\Nmini}{N^{\rm min}_{-i}}

\newcommand{\dren}{D} 
\newcommand{\drenh}{\dren_{1/2}}
\newcommand{\dkl}{d_{\rm KL}}

\newcommand{\vkl}{v_{\rm KL}}
\newcommand{\Hel}{\mathrm{Hel}}
\newcommand{\Ham}{\operatorname{Ham}}
\newcommand{\Hami}{\operatorname{Ham}_{-i}}
\newcommand{\dham}{\Ham}
\newcommand{\dhams}{\Ham^*}
\newcommand{\dhamsi}{d_{i}^*}
\newcommand{\ace}{\dhams}

\newcommand{\Mir}{\operatorname{Mir}}
\newcommand{\Rand}{\operatorname{Rand}}

\DeclareMathOperator*{\argmax}{arg\,max}

\newcommand{\Zgeo}{W}

\newcommand{\nik}{\frac{NI}{K}}
\newcommand{\fracnk}{\frac{N}{K}}

\newcommand{\Opt}{\operatorname{Opt}}

\newcommand{\Crit}{\operatorname{Crit}}

\newcommand{\nCrit}{L^+}

\newcommand{\LLR}{\Lambda}

\newcommand{\card}{\operatorname{card}}

\newcommand{\renyiMarkov}{\tilde I}

\endlocaldefs

\begin{document}

\begin{frontmatter}

\title{Community recovery in non-binary and temporal stochastic block models}

\runtitle{Community recovery in non-binary and temporal stochastic block models}

\begin{aug}
\author[A]{\fnms{Konstantin} \snm{Avrachenkov}\ead[label=e1]{k.avrachenkov@inria.fr}},
\author[A]{\fnms{Maximilien} \snm{Dreveton}\ead[label=e2]{maximilien.dreveton@gmail.com}}\\
\and
\author[B]{\fnms{Lasse} \snm{Leskel\"a}\ead[label=e3]{lasse.leskela@aalto.fi}}
\address[A]{Inria Sophia Antipolis, \printead{e1,e2}}
\address[B]{Aalto University, \printead{e3}}
\end{aug}

\begin{abstract}
This article studies the estimation of latent community memberships from pairwise interactions in a network of $N$ nodes, where the observed interactions can be of arbitrary type, including binary, categorical, and vector-valued, and not excluding even more general objects such as time series or spatial point patterns. As a generative model for such data, we introduce a stochastic block model with a general measurable interaction space $\cS$, for which we derive information-theoretic bounds for the minimum achievable error rate. These bounds yield sharp criteria for the existence of consistent and strongly consistent estimators in terms of data sparsity, statistical similarity between intra- and inter-block interaction distributions, and the shape and size of the interaction space. The general framework makes it possible to study temporal and multiplex networks with $\cS = \{0,1\}^T$, in settings where both $N \to \infty$ and $T \to \infty$, and the temporal interaction patterns are correlated over time.  For temporal Markov interactions, we derive sharp consistency thresholds. We also present fast online estimation algorithms which fully utilise the non-binary nature of the observed data. Numerical experiments on synthetic and real data show that these algorithms rapidly produce accurate estimates even for very sparse data arrays.
\end{abstract}

\begin{keyword}[class=MSC]
\kwd[Primary ]{62H30}
\kwd{60J10}
\kwd[; secondary ]{90B15} \kwd{91D30}
\end{keyword}
\begin{keyword}
\kwd{longitudinal network data}
\kwd{random graphs}
\kwd{temporal networks}
\kwd{multiplex networks}
\kwd{multilayer networks}
\kwd{stochastic block model}
\kwd{planted partition model}
\kwd{community recovery}
\end{keyword}
\end{frontmatter}

\tableofcontents
\import{text}{text_main}

\begin{funding}
This work has been done within the project of Inria - Nokia Bell Labs ``Distributed Learning and Control for Network Analysis'' and was partially supported by COSTNET Cost Action CA15109. \end{funding}

\ifdraft
\bibliographystyle{alpha}
\else
\bibliographystyle{imsart-number}
\fi
\bibliography{main}

\clearpage

\begin{appendix}
\import{text}{text_appendix}
\end{appendix}

\end{document}

%% file: text/text_main.tex
\section{Introduction}

Data sets in many application domains consist of non-binary pairwise interactions. Examples include human interactions in sociology and epidemiology \cite{Lewis_Gonzalez_Kaufman_2012, Mastrandrea_Fournet_Barrat_2015,Zhao_Wang_Li_Wang_Wang_Gao_2014}, brain activity measurements in neuroscience \cite{Bassett_Wymbs_Porter_Mucha_Carlson_Grafton_2011}, and financial interactions in economics \cite{Mazzarisi_Barucca_Lillo_Tantari_2020}. 
Pair interactions are usually characterised by types (attributes, labels, features) of interacting objects (nodes, agents, individuals), and a set of objects with a common type is called a community (block, group, cluster). An important unsupervised learning problem is to infer the community memberships from the observed pair interactions, a task commonly known as \new{community recovery} or \new{clustering} \cite{Fortunato_2010}.
	
Temporal interactions are an important particular case of non-binary interactions. The longitudinal nature of such data calls for replacing classical graph-based models by temporal and multiplex network models  \cite{Hartle_Papadopoulos_Krioukov_2021,Holme_Saramaki_2012, Kivela_etal_2014}. Although many powerful clustering methods exist for static networks (spectral methods~\cite{Lei_Rinaldo_2015}, semidefinite programming~\cite{Hajek_Wu_Xu_2016}, modularity maximisation \cite{Bickel_Chen_2009}, belief propagation \cite{Mezard_Montanari_2009}, Bayesian methods \cite{Peixoto_2019}, likelihood-based methods \cite{Wang_Bickel_2017}),
their extension to dynamic networks is not necessarily straightforward. In particular, simple approaches employing a static clustering method to a temporally aggregated network may lead to a severe loss of information~\cite{Avrachenkov_Dreveton_Leskela_2021}, and they are in general ill-suited to online updating.


The stochastic block model (SBM), first explicitly defined in \cite{Holland_Laskey_Leinhardt_1983}, has become a standard framework for analysing network data with binary interactions.
The present article extends the definition of the stochastic block model to its most general form
in which the observed interactions can be of arbitrary type, including binary, categorical, and vector-valued, and not excluding even more general objects such as time series or spatial point patterns. The observed data are represented by an $N$-by-$N$ symmetric array with entries in a general measurable space $\cS$. The binary case with $\cS=\{0,1\}$ corresponds to the most studied setting of random graphs. Temporal and multiplex networks can be represented by choosing $\cS = \{0,1\}^T$ where $T$ equals the number of snapshots or layers. Other important choices for the interaction space include $\cS = \{0, 1 \dots, L\}$ (link-labelled SBMs) and $\cS = \R$ (weighted SBMs).

\subsection{Related work}

Existing works on community recovery in binary networks provide a strong information-theoretic foundation
\cite{Gao_Ma_Zhang_Zhou_2017,Mossel_Neeman_Sly_2016,Zhang_Zhou_2016}. In particular, for $\cS = \{0,1\}$ it is known that communities can be consistently recovered if the difference between intra- and inter-block link probabilities is large enough. Similar conclusions have been extended to models with categorical ($\cS=\{0,\dots,L\}$) interactions
\cite{Heimlicher_Lelarge_Massoulie_2012,Jog_Loh_2015,Lelarge_Massoulie_Xu_2015,Xu_Jog_Loh_2020,Yun_Proutiere_2016} and real-valued ($\cS=\R$) interactions \cite{Xu_Jog_Loh_2020}. 
In principle, temporal and multiplex network data with $\cS = \{0,1\}^T$ could be modelled as categorical interactions with $2^T$ categories, but such approaches suffer from the following limitations.
First, existing theoretical results are mainly limited to models with a bounded or slowly growing number of categories. For example, the results in \cite{Xu_Jog_Loh_2020} will directly apply only for $L = o(N)$.
Second, algorithms designed for categorical interactions typically have complexity linear in $L$, and are hence inefficient even for a modest number of snapshots.

The present article is motivated by the inference of community structures from temporal network data; see \cite{Hartle_Papadopoulos_Krioukov_2021} for a comprehensive review of dynamic network models.  Earlier works on models, algorithms, and data experiments on temporal networks include
\cite{%
Ghasemian_Zhang_Clauset_Moore_Peel_2016,%
Longepierre_Matias_2019,%
Matias_Miele_2017,%
Pensky_2019,%
Xu_Hero_2014,%
Yang_Chi_Zhu_Gong_Jin_2011%
}, where interactions are assumed temporally uncorrelated given the community memberships.
Some of the aforementioned works also allow for time-varying community memberships. Because time-varying community memberships are known to involve model identifiability issues \cite{Matias_Miele_2017}, this feature is left out of the scope of the present article. Information-theoretic studies on multiplex networks with independent layers include \cite{Han_Xu_Airoldi_2015} presenting a strongly consistent estimator for models with $N = O(1)$ and $T \gg 1$,  
\cite{Paul_Chen_2016} establishing minimax error rates for models with $N,T \gg 1$ and balanced community sizes, \cite{Alaluusua_Leskela_2022} establishing posterior consistency in a Bayesian framework, and \cite{Bhattacharyya_Chatterjee_2018-05-27,Bhattacharyya_Chatterjee_2020-04-06,Lei_Lin_2022,Paul_Chen_2020,Pensky_Zhang_2019} presenting consistent estimators based on spectral clustering.
%
%
Dynamic networks with temporally correlated interactions, or persistent edges, have so far attracted much less attention.
Articles \cite{Barucca_Lillo_Mazzarisi_Tantari_2018,Mazzarisi_Barucca_Lillo_Tantari_2020} present numerical algorithms for estimating community memberships in temporally correlated SBMs in which the interaction patterns between nodes are positively correlated discrete-time Markov chains.
%
%
Recently,
\cite{Rastelli_Fop_2020,Suveges_Olhede_2022} presented EM algorithms for temporal SBMs where interactions are continuous-time Markov processes. 

A detailed technical discussion of our contributions with respect to the most closely related earlier works
is postponed to Section~\ref{section:comparison_other_work}. 

\subsection{Main contributions}

The main contributions of the present article can be summarised as follows:
\begin{enumerate}
\item We extend the SBM analysis to a general framework which allows the size and shape of the space of interactions $\cS$ to vary with scale, making it possible to analyse vector-valued and functional interactions with dimension growing with scale, and 
multiplex and temporal networks where the number of layers or snapshots goes to infinity.

\item We derive a lower bound on the minimum achievable error rate of community recovery in a SBM with general interactions, including binary, categorical, weighted, and temporal patterns. This result extends in a natural but non-trivial way earlier results for binary and real-valued SBMs, 
by allowing the space of interactions $\cS$ and the interaction distributions to be arbitrary. This is one of the first explicit quantitative lower bounds in this context.

\item 
We show that the maximum likelihood estimator recovers the true communities up to the information-theoretic lower bound. Combined with the lower bound, this yields sharp thresholds for community recovery in terms of the \Renyi divergence between the interaction distributions.
We also propose a polynomial-time algorithm which attains the desired lower bound under mild additional regularity assumptions.

\item We analyse temporal SBMs where interactions between nodes are correlated over time, and both the number of nodes $N$ and the number of time slots $T$ may tend to infinity. For sparse networks with Markov interactions, we derive information-theoretic consistency thresholds. The thresholds are presented in terms of an asymptotic formula for the \Renyi divergence between two sparse Markov chains, which could be of independent interest.

\item We provide online algorithms for temporal networks in situations where the interaction parameters are known or unknown, with complexity linear in the number of layers $T$. 
%
In particular, a numerical study demonstrates that in a typical situation, we recover the correct communities starting from a blind random guess, even in very sparse regimes.
\end{enumerate}

\subsection{Outline}
The rest of the article is structured as follows. 
Section~\ref{subsection:model} describes model details and notations.
Section~\ref{section:recovery_threshold} summarises the main theoretical results for general non-binary network models,
and Section~\ref{section:temporal_networks} specialises to temporally correlated networks.
%
Section~\ref{section:numerical_results}
describes numerical experiments on synthetic and real data sets.
Section~\ref{section:comparison_other_work} provides a technical discussion on our main contributions with respect to the state of the art. Finally, Section~\ref{section:conclusion} 
describes avenues for future research. The proofs of the main theorems are presented in the appendices.

\section{Model description and notations}
\label{subsection:model}

\subsection{General stochastic block model}

The objective of study is a population of $N \ge 1$ mutually interacting nodes partitioned into $K \ge 2$ disjoint sets called blocks. The partition is represented by a node labelling $\sigma: [N] \to [K]$, so that $\sigma(i)$ indicates the block which contains node $i$. In line with the classical definition of a stochastic block model \cite{Holland_Laskey_Leinhardt_1983}, we assume that interactions between node pairs can be of arbitrary type, and the set of possible interaction types is a measurable space $\cS$. This general setup allows to model usual random graphs ($\cS = \{0,1\}$), edge-labelled random graphs ($\cS = \{0,\dots,L\}$, $\cS = \R$), multilayer and temporal networks ($\cS = \{0,1\}^T$, $\cS = \{0,1\}^\infty$), and many other settings such as nodes interacting over a continuous time interval. In full generality, such a \new{stochastic block model (SBM)} is parameterised by a node labelling $\sigma: [N] \to [K]$ and an interaction kernel $(f_{k\ell})$ which is a collection of probability density functions with respect to a common sigma-finite reference measure $\mu$ on $\cS$, such that $f_{k\ell} = f_{\ell k}$ for all $k, \ell = 1,\dots,K$. These parameters specify a probability measure on a space of observations
\[
 \cX
 \weq \Big\{x: [N]\times [N] \to \cS: \ x_{ij} = x_{ji}, \ x_{ii}=0 \ \text{for all $i,j$} \Big\}
\] 
with probability density function
\begin{equation}
 \label{eq:PairwiseInteractionModel}
 P_{\sigma}(x)
 \ = \prod_{1 \le i < j \le N} f_{\sigma(i) \sigma(j)} \left( x_{ij} \right)
\end{equation}
with respect to the $N(N-1)/2$-fold product of the reference measure $\mu$. Our main focus is on homogeneous models in which the interaction kernel can be represented as
\begin{equation}
 \label{eq:Homogeneous}
 f_{k\ell}
 \weq
 \begin{cases}
  f & \text{ if } k = \ell, \\
  g & \text{ otherwise,}
 \end{cases}
\end{equation} 
for some probability densities $f$ and $g$ on $\cS$, called the intra-block and inter-block interaction distribution, respectively. A \new{homogeneous SBM} is hence a probability density $P_\sigma$ on $\cX$ specified by \eqref{eq:PairwiseInteractionModel}--\eqref{eq:Homogeneous} and parameterised by a 5-tuple $(N, K, \sigma, f, g)$. For an observation $X$ distributed according to such $P_\sigma$, the entries $X_{ij}$, $1 \le i < j \le N$, are mutually independent, and $X_{ij}$ is distributed according to $f$ when $\sigma(i)=\sigma(j)$, and according to $g$ otherwise.

The node labelling $\sigma$ representing the block membership structure is considered an unknown parameter to be estimated. When studying the average error rate of estimators, it is natural to regard the node labelling as a random variable distributed according to the uniform distribution $\pi(\sigma) = K^{-N}$ on parameter space
$
 \cZ
 = \big\{ \sigma: [N]\to [K] \big\}.
$
In this case the joint distribution of the node labelling and the observed data is characterised by a probability density
\begin{equation}
\label{eq:JointProbability}
 \pr(\sigma,x) 
 \weq \pi_\sigma P_\sigma(x)
\end{equation}
on $\cZ \times \cX$ 
with respect to $\card_\cZ \times \mu$, where $\card_\cZ$ is the counting measure on $\cZ$.

\subsection{Classification error}
	
The community recovery problem is the task of developing an algorithm $\phi: \cX \to \cZ$ which maps an observed data array $X = (X_{ij})$ into an estimated node labelling $\hsigma = \phi(X)$.
Stated like this, the recovery problem is ill-posed because the map $\sigma \mapsto P_\sigma$ defined by \eqref{eq:PairwiseInteractionModel} is in general non-injective. Therefore, we adopt the common approach in which the goal is to recover the unlabelled block structure, that is, the partition $[\sigma] = \{ \sigma^{-1}(k): k \in [K]\}$, and the estimation error is considered small when $[\hsigma]$ is close to $[\sigma]$. Accordingly, we define for node labellings $\sigma_1, \sigma_2: [N] \to [K]$ an error quantity by
\[
 \dhams (\sigma_1, \sigma_2)
 \weq \min_{\rho \in \Sym(K)} \dham( \rho \circ \sigma_1, \sigma_2)
\]
where $\Sym(K)$ denotes the group of permutations on $[K]$ and $\dham$ refers to the Hamming distance. The above error takes values in $\{0,\dots,N\}$ and depends on its inputs only via the partitions $[\sigma_1]$ and $[\sigma_2]$. The normalised error quantity $N^{-1} \dhams (\sigma_1, \sigma_2)$ is known as the classification error
\cite{Fortunato_2010,Meila_2007}.

When analysing the average performance of an estimator, we can view $\hsigma$ as a $\cZ$-valued random variable defined on the observation space $\cX$. Then $E_\sigma \dhams(\hsigma, \sigma)$ equals the expected clustering error given a true parameter $\sigma$, and 
\[
 \E \dhams(\hsigma)
 \weq \sum_{\sigma \in \cZ} \pi_\sigma E_\sigma \dhams(\hsigma, \sigma)
\]
is the average clustering error with respect to the uniform distribution $\pi_\sigma = K^{-N}$ on the parameter space.

\subsection{Consistent estimators}

A large-scale network is represented as a sequence of models $P_\sigma^{(\eta)}$ indexed by a scale parameter $\eta=1,2,\dots$ In this setting the model dimensions $N^{(\eta)}, K^{(\eta)}$, the node labelling $\sigma^{(\eta)}$, the interaction densities $f^{(\eta)}, g^{(\eta)}$, as well as the spaces $\cS^{(\eta)}, \cX^{(\eta)}, \cZ^{(\eta)}$ all depend on the scale parameter $\eta$. In this setup, an estimator is viewed as a map $\phi^{(\eta)}: \cX^{(n)} \to \cZ^{(\eta)}$. For nonnegative sequences $a=a^{(\eta)}$ and $b=b^{(\eta)}$ we denote $a = o(b)$ when $\limsup_{\eta \to \infty} a^{(\eta)}/b^{(\eta)} = 0$, and $a = O(b)$ when $\limsup_{\eta \to \infty} a^{(\eta)}/b^{(\eta)} < \infty$. We write $a = \omega(b)$ when $b = o(a)$, $a = \Omega(b)$ when $b = O(a)$, and $a=\Theta(b)$ when $a=O(b)$ and $b=O(a)$. We also denote $a \ll b$ for $a = o(b)$, $a \lesim b$ for $a = O(b)$, $a \asymp b$ for $a = \Theta(b)$, and $a \sim b$ for $a = (1+o(1)b$.
To avoid overburdening the notation, the scale parameter is mostly omitted from the notation in what follows.

For a large-scale model with $N \gg 1$ nodes, an estimator $\hsigma = \hsigma^{(\eta)}$ is called \new{consistent} if $\E \dhams(\hsigma) = o(N)$, and \new{strongly consistent} if
 $\E \dhams(\hsigma) = o(1)$.
A strongly consistent estimator is also said to achieve \new{exact recovery}, and a consistent estimator is said to achieve \new{almost exact recovery}~\cite{Abbe_2018_JMLR}.

\subsection{Information-theoretic divergences and distances}

Let us recall basic information divergences and distances associated with probability distributions $f$ and $g$ on a general measurable space $\cS$ \cite{Ghosal_VanDerVaart_2017,vanErven_Harremoes_2014}. 
The \new{\Renyi divergence} of positive order $\alpha \not = 1$ is defined as 
\begin{align*}
 \dren_\alpha ( f \| g ) = (\alpha-1)^{-1} \log \int \left( \frac{df}{d \mu} \right)^{\alpha}  \left( \frac{dg}{d \mu} \right)^{1-\alpha} d \mu,
\end{align*}
where $\mu$ is an arbitrary measure which dominates $f$ and $g$. We use the conventions $\log 0 = -\infty$, $0/0 = 0$ and $x/0 = \infty$ for $x>0$. In particular, if $f \not \perp g$ and $\alpha < 1$, then $\dren_\alpha(f\|g) < \infty$. In the symmetric case with $\alpha=\frac12$ we write
\[
 \drenh(f,g) \weq -2 \log \int \sqrt{ \frac{df}{d \mu} } \sqrt{\frac{dg}{d \mu} } d \mu,
\]
and note that this quantity is related to the \new{Hellinger distance} defined by
\[
 \Hel^2(f,g)
 \weq \frac12 \int \left( \sqrt{ \frac{df}{d \mu} }
   - \sqrt{ \frac{dg}{d \mu} } \right)^2 d \mu,
\]
via the formula $\drenh(f,g) = - 2 \log \left( 1 - \Hel^2(f,g) \right)$. In what follows, we assume that a sigma-finite reference measure $\mu$ on $\cS$ is fixed once and for all, and we write $\frac{df}{d\mu}, \frac{dg}{d\mu}$ simply as $f,g$, and we omit $d\mu$ from the integral signs, so that $\dren_\alpha(f\|g) = (\alpha-1)^{-1} \log \int f^\alpha g^{1-\alpha}$. When $\cS$ is countable, $\mu$ is always chosen as the counting measure, in which case write $\dren_\alpha(f\|g) = (\alpha-1)^{-1} \log \sum_{x\in\cS} f^\alpha(x) g^{1-\alpha}(x)$, and so on.
We also denote symmetrised \Renyi divergences by $\dren^s_\alpha(f,g) = \frac12 ( \dren_\alpha(f \| g) + \dren_\alpha(g \| f))$.

\section{Results for general SBMs}
\label{section:recovery_threshold}

Section~\ref{subsection:recovery_conditions_homogeneous_sbm} describes information-theoretic thresholds for consistent community recovery.
Section~\ref{sec:SparseNetworks} specialises to sparse networks.
Section~\ref{subsubsection:polynomial_time_algo} describes a polynomial-time algorithm and discusses its accuracy.

\subsection{General information thresholds}
\label{subsection:recovery_conditions_homogeneous_sbm}

The following theorem characterises fundamental information-theoretic limits for the recovery of block memberships from data generated by a homogeneous $\cS$-valued SBM. It does not make any scaling assumptions on the model dimensions $N$ and $K$, or on the space of interaction types $\cS$, and its proof indicates that maximum likelihood estimators achieve the upper bound.

\begin{theorem}
\label{thm:asymptotic_exponential_bound_recovery}
For a homogeneous SBM with $N$ nodes, $K$ blocks, and interaction distributions $\fin, \fout$ on a general measurable space $\cS$ having \Renyi divergence $I = \drenh(f,g)$, the minimum average classification error among all estimators $\hsigma: \cX \to \cZ$ is bounded from below by
\[
 \min_{\hsigma}
 \E \left( \frac{\dhams(\hsigma)}{N} \right)
 \wge \frac{1}{84} K^{-3} e^{ - \frac{N}{K} I - \sqrt{ 8 N I_{21} } } - \frac16 e^{ - \frac{N}{8K} }
\]
and from above by
\[
 \min_{\hsigma} \E \left( \frac{\dhams(\hsigma)}{N}   \right)
 \wle 8 e (K-1) e^{-(1 - \zeta - \kappa) \frac{N}{K} I}
 + K^N e^{-\frac14 (\frac{\zeta}{K-1} - \epsilon) (N/K)^2 I }
 + 2 K e^{-\frac13 \epsilon^2 \frac{N}{K}},
\]
for all $0 \le \epsilon \le \zeta \le \frac{1}{21}$, where $\kappa = 56 \max\{ K^2 e^{- \frac{NI}{8K}}, \, K N^{-1} \}$ and another auxiliary parameter is defined by $I_{21} = \left( \frac12 - K^{-1} \right) K^{-1} I + \frac12 K^{-1} J$
with $J = \left( \int \sqrt{fg} \right)^{-1} \int \sqrt{fg} \log^2 \frac{f}{g}$.
\end{theorem}

\begin{proof}
The lower bound is established in Proposition~\ref{the:LowerBoundHom} in Appendix~\ref{appendix:lower_bound_proof}, while the upper bound is analysed in Appendix~\ref{appendix:mle_consistency} and follows from Proposition~\ref{the:MLEUpperAverage}.
\end{proof}

The next key result characterises information-theoretic recovery conditions in large-scale networks, for which we emphasise that the model dimensions $N = N^{(\eta)}$ and $K= K^{(\eta)}$, the interaction distributions $f = f^{(\eta)}$ and $g = g^{(\eta)}$, and also the interaction type space $\cS = \cS^{(\eta)}$, are allowed to depend on a scale parameter $\eta$ which omitted from notation for clarity.

\begin{theorem}
\label{cor:recovery_conditions}
For a homogeneous SBM with $N \gg 1$ nodes, $K \asymp 1$ blocks, and interaction distributions $\fin$ and $\fout$ having \Renyi divergence $I = \drenh(f,g)$:
\begin{enumerate}[(i)]
\item a consistent estimator exists if
$I \gg N^{-1}$, and does not exist if $I \lesim N^{-1}$;
\item a strongly consistent estimator exists if
$I \ge (1+\Omega(1)) \frac{K \log N}{N}$, and does not exist if
$I \le (1-\Omega(1)) \frac{K \log N}{N}$.
\end{enumerate}
Furthermore, if $N^{-1} \ll I \ll 1$ and $K \asymp 1$, the optimal achievable misclassification rate equals
\begin{equation}
 \label{eq:OptimalRate}
 \min_{\hsigma} \E \left( \frac{\dhams(\hsigma)}{N} \right)
 \wasymp e^{-(1-o(1))NI/K}.
\end{equation}
\end{theorem}
\begin{proof}
The nonexistence statements are a direct consequence of the lower bound in Theorem~\ref{thm:asymptotic_exponential_bound_recovery} combined with Lemma~\ref{lemma:bounding_J_over_I} to guarantee that $J \lesim I$.
The existence results follow by analysing the upper bound of Theorem~\ref{thm:asymptotic_exponential_bound_recovery}, which is done in Proposition~\ref{the:MLEConsistent} in Appendix~\ref{appendix:mle_consistency}.
Formula \eqref{eq:OptimalRate} follows from the bounds of 
Theorem~\ref{thm:asymptotic_exponential_bound_recovery} by choosing $\zeta = (NI)^{-1/2}$ and $\epsilon = \frac12 \frac{\zeta}{K-1}$, and recalling $J \lesim I$ by Lemma~\ref{lemma:bounding_J_over_I}.
\end{proof}

The following examples illustrate how Theorem~\ref{cor:recovery_conditions} can be applied to various types of SBMs in sparse and dense regimes.

\begin{example}[Binary interactions]
\label{exa:BinarySBM}
A graph in which two nodes in the same community (resp.\ different communities) are linked with probability $p$ (resp.\ $q$) forms an instance of a homogeneous SBM where the $\frac12$-order \Renyi divergence between Bernoulli interaction distributions equals
$I = -2 \log( (1-p)^{1/2} (1-q)^{1/2} + p^{1/2} q^{1/2} )$. In a sparse regime where
$p = p_0 \frac{\log N}{N}$
and
$q = q_0 \frac{\log N}{N}$
for scale-independent constants $p_0,q_0>0$,
this is approximated by
$I \sim ( \sqrt{p_0} - \sqrt{q_0})^2 \frac{\log N}{N}$.
Theorem~\ref{cor:recovery_conditions} tells that
a strongly consistent estimator exists if $( \sqrt{p_0} - \sqrt{q_0})^2 > K$ and does not if $( \sqrt{p_0} - \sqrt{q_0})^2 < K$. This is the well-known threshold for strong consistency in sparse binary SBMs \cite{Abbe_Bandeira_Hall_2016,Mossel_Neeman_Sly_2016}. Alternatively, in a dense regime where 
$p = p_0 + \epsilon$ and $q = p_0$ for
$\epsilon = o(1)$ and some scale-independent constant $0 < p_0 < 1$, we find that
$I \sim \frac{\epsilon^2}{4p_0(1-p_0)}$. Especially, when
$\epsilon = c (\frac{\log N}{N})^{1/2}$ for some scale-independent constant $c>0$, then we find that
a strongly consistent estimator exists if
$(p_0-\frac12)^2 > \frac{1-c^2/K}{4}$
and does not if
$(p_0-\frac12)^2 < \frac{1-c^2/K}{4}$.
\end{example}

\begin{example}[Poisson interactions]
\label{exa:PoissonSBM}
Consider an integer-valued SBM where the interaction between two nodes in the same block (resp.\ different blocks) is a Poisson-distributed random integer with mean $\lambda$ (resp.\ $\mu$).  The $\frac12$-order \Renyi divergence between such Poisson distributions equals
$I = ( \sqrt{\lambda} - \sqrt{\mu})^2$.
In a sparse regime where
$\lambda = \lambda_0 \frac{\log N}{N}$
and
$\mu = \mu_0 \frac{\log N}{N}$
for scale-independent constants
$\lambda_0, \mu_0 > 0$,
Theorem~\ref{cor:recovery_conditions} tells that
a strongly consistent estimator exists if
$( \sqrt{\lambda_0} - \sqrt{\mu_0})^2 > K$ and does not if
$(\sqrt{\lambda_0} - \sqrt{\mu_0})^2 < K$.  In a dense regime where
$\lambda = \lambda_0 + \epsilon$ and $\mu = \lambda_0$ for $\epsilon = o(1)$ and some scale-independent constant $\lambda_0 > 0$, we see that
$I \sim \frac{\epsilon^2}{4\lambda_0}$. Especially, if $\epsilon = c (\frac{\log N}{N})^{1/2}$ for some scale-independent constant $c>0$, then 
a strongly consistent estimator exists when
$\frac{c^2}{4\lambda_0} > K$ 
and does not when
$\frac{c^2}{4\lambda_0} < K$.
\end{example}

\begin{example}[Normal interactions]
\label{exa:NormalSBM}
Consider a real-valued SBM where the interaction between two nodes in the same block (resp.\ different blocks) follows a normal distribution with mean zero and standard deviation $\sigma$ (resp.\ $\tau$).  The $\frac12$-order \Renyi divergence between such normal distributions equals
$I
= \log(1 + \frac{(\sigma-\tau)^2}{2\sigma\tau})
$.
Theorem~\ref{cor:recovery_conditions} combined with Taylor's approximation $\log(1+t) = t + O(t^2)$ 
tells that
a consistent estimator exists if $\frac{(\sigma-\tau)^2}{2\sigma\tau} \gg N^{-1}$
and does not if
$\frac{(\sigma-\tau)^2}{2\sigma\tau} \lesim N^{-1}$; and that
a strongly consistent estimator exists if
$\frac{(\sigma-\tau)^2}{2\sigma\tau} \ge (1+\Omega(1)) K \frac{\log N}{N}$
and does not if
$\frac{(\sigma-\tau)^2}{2\sigma\tau} \le (1-\Omega(1)) K \frac{\log N}{N}$.
\end{example}

\begin{example}[Multiplex networks]
\label{exa:MultiplexSBM}
An SBM with product-form intra- and inter-block interaction distributions $f = \prod_{t=1}^T f_t$ and $g = \prod_{t=1}^T g_t$ on $\cS = \cS_1 \times \cdots \times \cS_T$ corresponds to observing $T$ mutually independent network layers over a common node set, where data on the $t$-th layer are distributed according to an SBM with interaction distributions $f_t$ and $g_t$ on $\cS_t$. By observing that $I = \sum_{t=1}^T I_t$ for
$I = \dren_{1/2}(f,g)$ and
$I_t = \dren_{1/2}(f_t,g_t)$, 
Theorem~\ref{cor:recovery_conditions} tells that
strong consistency is possible when 
$
\sum_{t=1}^T I_t
\ge ( 1 + \Omega(1) ) \frac{K \log N}{N},
$
and impossible when
$
\sum_{t=1}^T I_t
\le ( 1 - \Omega(1) ) \frac{K \log N}{N}.
$
Similarly, consistency is possible when $\sum_{t=1}^T I_t \gg N^{-1}$ and impossible when $\sum_{t=1}^T I_t \lesim N^{-1}$.
In the binary case where $\cS_t = \{0,1\}$ for all $t$, corresponding thresholds have been derived in \cite{Alaluusua_Leskela_2022,Paul_Chen_2016}. The present example is an important extension allowing to analyse heterogeneous multiplex networks in which some layers may only be partially observed (cf.\ Example~\ref{exa:BinaryCensoredSBM}) and some may carry real-valued edge labels (cf.\ Example~\ref{exa:NormalSBMAlgo}).
\end{example}

\subsection{Sparse networks}
\label{sec:SparseNetworks}
Sparse networks can be modelled using intra-block and inter-block interaction distributions of form
\begin{align}
\label{eq:zero_inflated_distributions}
 f \weq (1-p_0 \rho )\delta_0 + p_0 \rho \tf
 \qquad \text{and} \qquad 
 g \weq (1-q_0 \rho ) \delta_0 + q_0 \rho \tg,
\end{align}
where $\delta_0$ is the Dirac measure at
an element $0 \in \cS$ representing no-interaction, probability measures $\tf,\tg$ on $\cS \setminus \{0\}$ are conditional distributions of interaction types given that there is an interaction, $p_0,q_0 > 0$ are scale-independent constants, and $\rho \ll 1$ describes the overall network density. The following result describes how the regimes for consistent and strongly consistent community recovery are characterised by a fundamental information quantity
\begin{equation}
 \label{eq:InformationSparse}
 \tI
 \weq \left( \sqrt{p_0} - \sqrt{q_0} \right)^2
 + 2\sqrt{p_0 q_0} \, \Hel^2(\tf, \tg).
\end{equation}
In the above formula, the quantity
$\left( \sqrt{p_0} - \sqrt{q_0} \right)^2$ corresponds to information gained from observing whether or not there is an interaction, and the Hellinger distance $\Hel(\tf,\tg)$ characterises the additional information gained by observing the types of interactions between node pairs.

\begin{theorem}
\label{the:ConsistencySparse}
For a homogeneous SBM with $N \gg 1$ nodes, $K \asymp 1$ blocks, and interaction distributions of form \eqref{eq:zero_inflated_distributions} 
where $\rho \ll 1$, and $p_0,q_0$ are scale-independent constants:
\begin{enumerate}[(i)]
\item \label{ite:Consistent}
a consistent estimator exists if
$\rho \tilde I \gg N^{-1}$, and does not exist if $\rho \tilde I \lesim N^{-1}$;
\item \label{ite:StronglyConsistent}
a strongly consistent estimator exists if
$\rho \tilde I \ge (1+\Omega(1)) \frac{K \log N}{N}$, and does not exist if
$\rho \tilde I \le (1-\Omega(1)) \frac{K \log N}{N}$;
\end{enumerate}
\end{theorem}

\begin{proof}
Taylor's approximations show that the \Renyi divergence of order $\alpha \ne 1$ for probability distributions of form \eqref{eq:zero_inflated_distributions} is approximated by
\begin{equation}
 \label{eq:renyi_divergence_zero_inflated_distributions_gen}
 \dren_\alpha (f \| g)
 \weq \frac{p^\alpha q^{1-\alpha}}{\alpha-1}
 e^{ (\alpha-1) \dren_\alpha( \tf \| \tg ) }
 - \frac{\alpha p + (1-\alpha)q }{ \alpha-1 } + O(\rho^2),
\end{equation}
where $p = p_0 \rho$ and $q = q_0 \rho$.
In particular, the formula $1-\Hel(\tf,\tg)^2 = e^{-\frac12 D_{1/2}(\tf,\tg)}$ implies that the
\Renyi divergence of order half is given by
\begin{align}
 \label{eq:renyi_divergence_zero_inflated_distributions}
 D_{1/2}(f,g)
 \weq (\sqrt{p}-\sqrt{q})^2
 + 2 \sqrt{pq} \, \Hel^2(\tf,\tg)
 + O(\rho^2),
\end{align}
so that $D_{1/2}(f,g) = (1+o(1)) \rho \tilde I$.  Statements (i) and (ii) hence follow from Theorem~\ref{cor:recovery_conditions}. 
\end{proof}

The following three examples illustrate the applicability of Theorem \ref{the:ConsistencySparse} for finite and real-valued interaction spaces.

\begin{example}[Sparse categorical interactions]
\label{exa:SparseCategoricalSBM}
Consider a categorical stochastic block model with intra- and inter-block interactions distributed according to \eqref{eq:zero_inflated_distributions}
in which $\tf$ and $\tg$ are probability distributions on $\{1,\dots,L\}$. 
The critical information quantity defined in~\eqref{eq:InformationSparse} can then be written as
\[
 \tI
 \weq \left( \sqrt{p_0} - \sqrt{q_0} \right)^2 + \sqrt{p_0 q_0} \sum_{\ell=1}^L \left( \sqrt{\tf(\ell)} - \sqrt{\tg(\ell)} \right)^2.
\]
This model was studied in \cite{Jog_Loh_2015} in a parameter regime with
$\rho=\frac{\log N/K}{N/K} \asymp K \frac{\log N}{N}$, where it was assumed that
neither $L$ nor the probabilities $\tf(\ell)$, $\tg(\ell)>0$ depend on the scale parameter. In this case $\tI$ is a scale-independent constant, and by applying Theorem~\ref{the:ConsistencySparse}:\eqref{ite:StronglyConsistent} we recover the main results of \cite{Jog_Loh_2015} stating that a strongly consistent estimator exists if $\tI > 1$ and does not exist if $\tI < 1$. Unlike \cite{Jog_Loh_2015}, Theorem~\ref{the:ConsistencySparse} does not require any regularity conditions on $\tf$ and $\tg$.

\end{example}

\begin{example}[Censored binary SBM]
\label{exa:BinaryCensoredSBM}
Assume that between any pair of nodes in the same community (resp.\ different communities), there is an edge with probability $a$ (resp.\ $b$) and the edge status of the node pair is observed with probability $p = p_0\rho$ (resp.\ $q = q_0\rho$) regardless of whether an edge is present or not.
We assume that $\rho \ll 1$, and that $a, b, p_0, q_0$ are scale-independent constants. The observed data can be modelled as an instance of \eqref{eq:zero_inflated_distributions} with interaction type space $\cS = \{0,10,11\}$ where 0 = censored, 10 = observed\&absent, and 11 = observed\&present, and
\begin{alignat*}{2}
 \tf(10) &= 1-a, & \qquad\qquad \tg(10) &= 1-b, \\
 \tf(11) &= a,   & \qquad\qquad \tg(11) &= b.
\end{alignat*}
The fundamental information quantity in \eqref{eq:InformationSparse} equals
$
 \tilde I
 = \left( \sqrt{p_0} - \sqrt{q_0} \right)^2
 + 2\sqrt{p_0 q_0} \, \Hel^2(\tf, \tg),
$
where
$
 \Hel(\tf, \tg)^2
 = \frac12 (\sqrt{1-a} - \sqrt{1-b})^2
  + \frac12 (\sqrt{a} - \sqrt{b})^2.
$
Assume that $t=p_0=q_0$ equals a common observation rate and $\rho = \frac{\log N}{N}$.  Theorem~\ref{the:ConsistencySparse} then tells that exact recovery is possible if
$t > t_{\rm crit}$ and impossible
if $t < t_{\rm crit}$.
where $t_{\rm crit} = \frac{K}{2 \Hel^2(\tf,\tg)}$.
For $K=2$, this coincides with the exact recovery threshold 
recently presented in Dhara et al.\ \cite{Dhara_Souvik_Mossel_Sandon_2022}, and extends their criterion into models with $K > 2$ and $p_0 \ne q_0$.
\end{example}

\begin{example}[Censored real-valued SBM]
\label{exa:NormalCensoredSBM}
Assume that associated to each pair of nodes in the same community (resp.\ different communities), there is a random variable following a normal distribution with mean zero and standard deviation $\sigma$ (resp.\ $\tau$), and this variable is observed with probability $p = p_0\rho$ (resp.\ $q = q_0\rho$), where $\rho \ll 1$ and $\sigma,\tau,p_0,q_0 > 0$ are scale-independent constants. The observed data can be modelled as an instance of \eqref{eq:zero_inflated_distributions} with interaction type space $\cS = \R$ in which $\tf = \Nor(0,\sigma^2)$ and $\tg = \Nor(0,\tau^2)$, and the value 0 represents no-observation.
The fundamental information quantity in \eqref{eq:InformationSparse} then equals
\[
 \tI
 \weq \left( \sqrt{p_0} - \sqrt{q_0} \right)^2
 + 2\sqrt{p_0 q_0} \, \left( 1 - \sqrt{\frac{2 \sigma \tau}{\sigma^2+\tau^2} } \right).
\]
If $t=p_0=q_0$ equals a common observation rate and $\rho = \frac{\log N}{N}$, then Theorem~\ref{the:ConsistencySparse} then tells that exact recovery is possible if $t > t_{\rm crit}$
and impossible if $t < t_{\rm crit}$ were $t_{\rm crit} = \frac{K/2}{1 - \sqrt{\frac{2 \sigma \tau}{\sigma^2+\tau^2} }}$.
\end{example}

\subsection{Polynomial-time algorithm}
\label{subsubsection:polynomial_time_algo}

To cluster non-binary SBMs in a polynomial time in $N$, we propose Algorithm~\ref{algo:likelihood_based_algo} which employs spectral clustering as a subroutine to produce a moderately accurate initial clustering, and then performs a refinement step through node-wise likelihood maximisation.
Similarly to \cite{Gao_Ma_Zhang_Zhou_2017, Xu_Jog_Loh_2020}, for technical reasons related to the proofs, the initialisation step of Algorithm~\ref{algo:likelihood_based_algo} involves $N$ separate spectral clustering steps. A consensus step is therefore needed at the end, to correctly permute the individual predictions. Numerical experiments indicate that in practice it often suffices to do one spectral clustering on a binary matrix, and remove this consensus step. We will discuss practical aspects in more detail in Section~\ref{section:numerical_results}. 

\begin{algorithm}[!ht]
\small
\DontPrintSemicolon
\KwIn{$\cS$-valued interaction array $X_{ij}$; interaction distributions $f, g$; set $\cA \subset \cS$.}
\KwOut{Estimated node labelling $\hsigma$.}
\BlankLine
\BlankLine

{\it Step 1: Coarse clustering using binary interaction data}
\BlankLine
Compute a binary matrix $\tilde X$ by setting $\tX_{ij} = 1( X_{ij} \not \in \cA )$.
\\
\For{$i=1,\dots,N$}
{
Let $\tilde X_{-i}$ be the submatrix of $\tilde X$ with row $i$ and column $i$ removed.\\
Compute a node labelling $\tilde\sigma_i$ on $[N] \setminus \{i\}$ by applying a standard graph clustering algorithm with adjacency matrix $\tilde X_{-i}$.
}
\BlankLine
\BlankLine

{\it Step 2: Refined clustering using full interaction data}
\BlankLine
\For{$i=1,\dots,N$}
{
Compute $h_i(k) = \sum_{j: \tilde\sigma_{i}(j) = k } \log \frac{f(X_{ij})}{g(X_{ij})} $
for all $k \in [K]$. \\
Set $\hsigma_i(i) = \argmax_{k \in [K]} h_i(k)$ with arbitrary tie breaks.\\
Set $\hsigma_i(j) = \tilde\sigma_{i}(j)$ for $j \ne i$. 
}
\BlankLine
\BlankLine

{\it Step 3: Consensus}
\BlankLine
Select $\hsigma_1$ as a baseline node labelling and set $\hsigma(1) = \hsigma_1(1)$. \\
\For{$i=2,\dots,N$}
{
Set $\hsigma(i) = \argmax_\ell \abs{ \hsigma_i^{-1}( \hsigma_i(i) ) \cap \hsigma_1^{-1}(\ell) }$
with arbitrary tie breaks.
}
\caption{Clustering using general $\cS$-valued interaction data}
\label{algo:likelihood_based_algo}
\end{algorithm}

The following theorem characterises the accuracy of Algorithm~\ref{algo:likelihood_based_algo} for large-scale models, and implies
that under mild technical conditions this algorithm achieves the optimal error rate in Theorem~\ref{cor:recovery_conditions}. The proof of Theorem~\ref{thm:general_algo_consistency} is given in Appendix~\ref{appendix:consistency_likelihood_based_algo}.

\begin{theorem}
\label{thm:general_algo_consistency}
Consider a homogeneous SBM with $N \gg 1$ nodes, $K \asymp 1$ blocks, and interaction distributions $f$ and $g$ having \Renyi divergence $I = \drenh(f,g)$. 
If $(f(\cA)^{1/2} - g(\cA)^{1/2})^2
\gg N^{-1} \frac{\dren_{1+r}^s(f,g)}{\dren_{r}^s(f,g)}$
for some $0 < r \le \frac12$,
then the classification error of Algorithm~\ref{algo:likelihood_based_algo}
applied with $\cA \subset \cS$ is bounded by
\begin{align*}
 \E \left( \frac{\dhams(\hsigma)}{N} \right)
 \wle K e^{-(1 - o(1)) 2r \frac{N}{K}I} + o(1).
\end{align*}
\end{theorem}

The following three examples illustrate how Theorem~\ref{thm:general_algo_consistency} can be applied to analyse the performance of Algorithm~\ref{algo:likelihood_based_algo} in sparse and dense settings.

\begin{example}[Normal interactions]
\label{exa:NormalSBMAlgo}
Consider the real-valued SBM in Example~\ref{exa:NormalSBM} with intra- and inter-block interaction distributions $f = \Nor( 0, \sigma^2 )$ and $g = \Nor( 0, \tau^2)$.
Assume that $\sigma>0$ is scale-independent and $\tau = \sigma(1+\epsilon)$ for some $\epsilon = o(1)$. Then the $\frac12$-order \Renyi divergence equals $I \sim \frac12 \epsilon^2$, and the symmetrised $\frac32$-order \Renyi divergence is finite and approximated by $\dren_{3/2}^s(f,g) \sim \frac32 \epsilon^2$.
For an interval $\cA = [-x,x]$ we find that
$f(\cA) - g(\cA)
= 2( \Phi(\frac{x}{\sigma}) - \Phi(\frac{x}{\tau}))
= 2 \frac{x}{\sigma} \Phi'(\frac{x}{\sigma})  \epsilon + O(\epsilon^2)$,
where $\Phi$ is the standard normal cdf. 
Thus, $(\sqrt{f(\cA)} - \sqrt{g(\cA)})^2 \asymp \epsilon^2$ for any $x = O(1)$.  Theorem~\ref{thm:general_algo_consistency} hence tells that Algorithm~\ref{algo:likelihood_based_algo} applied with $\cA = [-1,1]$ is consistent if $\epsilon^2 \gg N^{-1}$ and strongly consistent if $\frac12 \epsilon^2 \ge ( 1 + \Omega(1)) \frac{K \log N}{N}$. 
In light of Example~\ref{exa:NormalSBM}, we see that Algorithm~\ref{algo:likelihood_based_algo} recovers communities up to the information-theoretic boundaries in Theorem~\ref{cor:recovery_conditions}.
\end{example}

\begin{example}[Geometric interactions]
\label{exa:GeometricSBMAlgo}
Suppose $f = \Geo(a)$ and $g = \Geo( a + \epsilon )$ with $a$ scale-independent and $\epsilon \ll 1$.
Then the $\frac12$-order \Renyi divergence is $I \sim \frac{\epsilon^2 }{4 a (1-a)^2}$. Theorem~\ref{cor:recovery_conditions} tells that a consistent estimator exists if $\epsilon \gg N^{-1/2}$ and a strongly consistent estimator exists if $\frac{ \epsilon^2 }{ a(1-a)^2 } \ge ( 1 + \Omega(1)) \frac{K \log N}{N}$.
We also note that $( \sqrt{f(0)} - \sqrt{g(0)})^2 \asymp \epsilon^2$, and that the symmetrised $\frac32$-order \Renyi divergence is finite and satisfies $\dren^s_{3/2}(f,g) \asymp \epsilon^2$.
Theorem~\ref{thm:general_algo_consistency} hence tells that Algorithm~\ref{algo:likelihood_based_algo} applied with $\cA = \{0\}$ recovers the communities up to the information-theoretic boundaries in Theorem~\ref{cor:recovery_conditions}.
\end{example}

\begin{figure}[!ht]
\centering
\includegraphics[
height=30mm
]{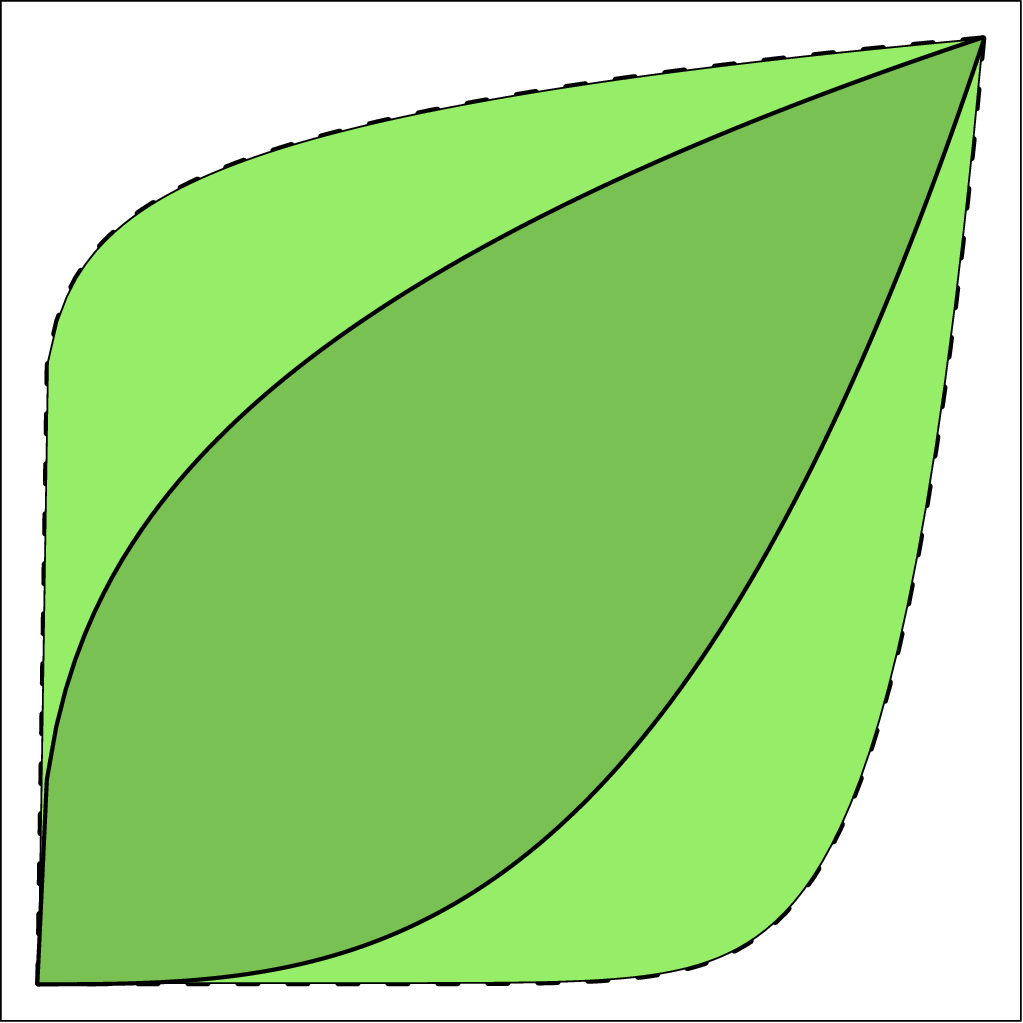}
\caption{$(a,b)$-pairs
such that
$b^{ \frac{1+r}{r} } < a < b^{ \frac{r}{1+r}}$ $r=0.1$ (green) and $r = 0.5$ (dark green).
\label{fig:leaf}}
\end{figure}

\begin{example}[Zero-inflated geometric distributions]
\label{exa:ZIGeometricSBMAlgo}
Consider an integer-valued stochastic block model with intra- and inter-block interactions distributed according to
\begin{align*}
 f(x) =
 \begin{cases}
  1-\rho p_0, &\quad x=0,\\
  \rho p_0 (1-a)a^x, &\quad x \ge 1,
 \end{cases}
 \qquad\text{and}\qquad
 g(x) = 
 \begin{cases}
  1-\rho q_0, &\quad x=0,\\
  \rho q_0 (1-b)b^x, &\quad x \ge 1,
 \end{cases}
\end{align*}
for some $\rho \ll 1$ and some scale-independent constants $p_0,q_0 > 0$ and $0 < a,b < 1$.
This is an instance of model \eqref{eq:zero_inflated_distributions} in which $\tf$ and $\tg$ are geometric distributions with parameters $a$ and $b$, and the critical information quantity in \eqref{eq:InformationSparse} equals
\[
 \tI
 \weq \left( \sqrt{p_0} - \sqrt{q_0} \right)^2
 + 2 \sqrt{p_0 q_0} \left( 1 -  \frac{(1-a)^{1/2} (1-b)^{1/2}}{1 - a^{1/2} b^{1/2}} \right).
\]
In this case, a higher order symmetrised \Renyi divergence $\dren_{1+r}^s (\tf, \tg)$ is finite if and only if
$b^{ \frac{1+r}{r} } < a < b^{ \frac{r}{1+r} }$, see Figure~\ref{fig:leaf}. This condition holds for a small enough $r>0$.
Let $\cA = \{0\}$ if $p_0 \not= q_0$ and $\cA = \{0,1\}$ otherwise. Theorem~\ref{the:ConsistencySparse} then tells that
Algorithm~\ref{algo:likelihood_based_algo} is consistent in the full information-theoretically feasible parameter range with $\rho \tI \gg 1$, and strongly consistent when $\rho \tI \ge (1+\Omega(1)) \frac{K \log N}{N}$ and $b^3 < a < b^{1/3}$. 
For strong consistency, the latter somewhat counterintuitive extra condition is needed to guarantee that the log-likelihood ratios used in Algorithm~\ref{algo:likelihood_based_algo}
are sufficiently well concentrated around their expected values.
\end{example}

\section{Results for temporal SBMs}
\label{section:temporal_networks}

This section is devoted to clustering nodes using temporally correlated network data. Section~\ref{subsection:temporal_networks_markov_interactions} provides consistency results for models where interaction patterns between node pairs are Markov chains over time. Asymptotic results for sparse interactions are based on information-theoretic divergences between binary Markov chains, which can be of independent interest.
Section~\ref{section:online_likelihood_algos} describes two online clustering algorithms for temporal networks, one assuming known interaction parameters, and the other adaptively learning the interaction parameters from data.

\subsection{Information thresholds for Markov SBMs}
\label{subsection:temporal_networks_markov_interactions}

As an instance of a network where interactions are correlated over time, we investigate an SBM with interaction space $\cS = \{0,1\}^T$ in which intra- and inter-block distributions are given by
\begin{equation}
\label{eq:Markov}
 f
 \weq \mu_{x_1 } P_{x_1,x_2} \cdots P_{x_{T-1},x_T},
 \qquad \text{and} \qquad 
 g
 \weq \nu_{x_1} Q_{x_1,x_2} \cdots Q_{x_{T-1},x_T},
\end{equation}
where $\mu, \nu$ are initial probability distributions, and $P,Q$ are stochastic matrices on~$\{0,1\}$.
This is an instance of the general SBM model in which the symmetric \Renyi divergence between interaction distributions, a key quantity in Theorems \ref{thm:asymptotic_exponential_bound_recovery}--\ref{cor:recovery_conditions}, equals
\begin{equation}
 \label{eq:RenyiMarkov}
 D_{1/2}(f,g)
 \weq -2\log \left( \sum_{x \in \{0,1\}^T} ( \mu_{x_1} \nu_{x_1})^{1/2}
 \prod_{t=2}^T (P_{x_{t-1} x_t} Q_{x_{t-1} x_t})^{1/2}
 \right).
\end{equation}

For fixed instances of transition parameters, $D_{1/2}(f,g)$ may be numerically computed by formula \eqref{eq:RenyiMarkov}. To gain analytical insight, we will derive simplified expressions corresponding to sparse chains where
$\mu_1, \nu_1, P_{01}, Q_{01} \lesim \rho$ for some $\rho \ll 1$.  Under this assumption, the expected number of 1's in any particular interaction pattern is $O(\rho T)$. Therefore, for $\rho T \ll 1$, the probability of observing an interaction between any particular node pair is small.
The following result presents
a key approximation formula
with proof
provided in Appendix~\ref{appendix:information_theoretic_divergences_sparse_markov_chains}.

\begin{proposition}
\label{thm:renyi_divergence_order_less_1_sparse_markov_chains}
Consider binary Markov chains with initial distributions $\mu,\nu$ and
transition probability matrices $P,Q$. Assume that
$\mu_1, \nu_1, P_{01}, Q_{01} \le \rho$ for some $\rho$ such that $\rho T \le 0.01$. Then the \Renyi divergence~\eqref{eq:RenyiMarkov} is approximated by $\abs{D_{1/2}(f,g) - I} \le 92 (\rho T)^2$, where
\begin{equation}    
 \label{eq:RenyiSparseMC}
 \begin{aligned}
 I
 &\weq (\sqrt{\mu_1}-\sqrt{\nu_1})^2
 +
 \bigg(
 (\sqrt{P_{01}}-\sqrt{Q_{01}})^2
 + 2 H_{11}^2 \sqrt{P_{01}  Q_{01}}
 \bigg) (T-1) \\
 &\qquad + 2 \Big(\Gamma\sqrt{\mu_1\nu_1} - \sqrt{P_{01}Q_{01}} \Big) H_{11}^2 \sum_{t=0}^{T-2} (1-\Gamma)^t
 \end{aligned}
\end{equation}
is defined in terms of 
$
 H_{11}^2
 = 1 - \frac{\sqrt{1-P_{11}} \sqrt{ 1-Q_{11}}}
 {1 - \sqrt{P_{11} Q_{11}}}
$,
and
$\Gamma = 1-\sqrt{P_{11}Q_{11}}$.
\end{proposition}

The quantities in \eqref{eq:RenyiSparseMC}
can be understood as follows. With the help of Taylor's approximations we see that
\begin{align*}
 (\sqrt{\mu_1}-\sqrt{\nu_1})^2
 &\weq \drenh(\Ber(\mu_1) \| \Ber(\nu_1)) + O(\rho^2), \\
 (\sqrt{P_{01}}-\sqrt{Q_{01}})^2
 &\weq \drenh(\Ber(P_{01}) \| \Ber(Q_{01})) + O(\rho^2), \\
 H_{11}
 &\weq \Hel( \Geo(P_{11}),   \Geo(Q_{11}) ).
\end{align*}
We also note that $\Gamma = 1-\sqrt{P_{11}Q_{11}}$ may be interpreted as an effective spectral gap averaged
%
over the two Markov
chains%
%
\footnote{
The nontrivial eigenvalues of transition matrices $P$ and $Q$ can be written as $\Lambda_P = {P_{11} - P_{01}}$ and $\Lambda_Q = {Q_{11} - Q_{01}}$. These are nonnegative when $P_{01} \le P_{11}$ and $Q_{01} \le Q_{11}$. The absolute spectral gaps characterising the mixing rates of these chains \cite{Levin_Peres_Wilmer_2008} are then $\Gamma_P = 1-\Lambda_P$ and $\Gamma_Q = 1-\Lambda_Q$. When $P_{01} \ll P_{11}$ and $Q_{01} \ll Q_{11}$, we find that $\Gamma = 1 - \sqrt{P_{11}Q_{11}} = 1 - (\Lambda_P \Lambda_Q)^{1/2} + o(1)$.
}.
%

\subsubsection{Short time horizon}

Consider a Markov SBM in which $T = O(1)$ is a scale-independent constant, and
\begin{equation}
 \label{eq:SparseMCShortParam}
 \begin{aligned}
 \mu_1 &\weq u \rho + o(\rho), &\quad 
 P_{01} &\weq p_{01} \rho + o(\rho), &\quad
 H_{11} &\weq h_{11} + o(1), \\
 \nu_1 &\weq v \rho + o(\rho), &\quad 
 Q_{01} &\weq q_{01} \rho + o(\rho), &\quad
 \Gamma &\weq \gamma + o(1), \\
 \end{aligned}
\end{equation}
for some constants $u, v, p_{01}, q_{01}, h_{11}, \gamma$, and define a constant $\renyiMarkov$ by
\begin{equation}
 \label{eq:iShortMarkov}
 \begin{aligned}    
 \renyiMarkov
 &\weq \left( \sqrt{u}-\sqrt{v} \right)^2
 +
 \left(
 \left( \sqrt{p_{01}}-\sqrt{q_{01}} \right)^2
 + 2 h_{11}^2 \sqrt{p_{01} q_{01}}
 \right) (T-1) \\
 &\qquad + 2 h_{11}^2 \Big(\gamma\sqrt{uv} - \sqrt{p_{01}q_{01}} \Big) \sum_{t=0}^{T-2} (1-\gamma)^t.
 \end{aligned}
\end{equation}

\begin{theorem}
\label{the:ShortMarkovSBM}
Consider a Markov SBM with $N \gg 1$ nodes, $K = O(1)$ blocks, and $T = O(1)$ snapshots, and assume that \eqref{eq:SparseMCShortParam} holds for some constants $u, v, p_{01}, q_{01}, h_{11}, \gamma \ge 0$ such that $\renyiMarkov \ne 0$, and some $\rho \ll 1$. Then:
\begin{enumerate}[(i)]
\item A consistent estimator does not exist for
$\rho \lesim \frac{1}{N}$ and does exist for $\rho \gg \frac{1}{N}$.
\item A strongly consistent estimator does not exist for $\rho \ll \frac{\log N}{N}$ and does exist for 
$\rho \gg \frac{\log N}{N}$.
\item In a critical regime with $\rho = \frac{\log N}{N}$, a strongly consistent estimator does not exist for $\renyiMarkov < K$ and does exist for $\renyiMarkov > K$.
\end{enumerate}
If we further assume that
$u, v, p_{01}, q_{01} > 0$, 
$u + (T-1)p_{01} \ne v + (T-1)q_{01}$,
$P_{10} \asymp Q_{10}$,
and $P_{11} \asymp Q_{11}$,
then Algorithm~\ref{algo:likelihood_based_algo} is consistent when $\rho \gg \frac{1}{N}$; and strongly consistent when $\rho \gg \frac{\log N}{N}$, or when
$\rho = \frac{\log N}{N}$ and $\renyiMarkov > K$.
\end{theorem}

\begin{proof}
By Proposition~\ref{thm:renyi_divergence_order_less_1_sparse_markov_chains}, we find that
\[
 D_{1/2}(f,g)
 \weq (1+o(1)) \renyiMarkov \rho + O\left( \rho^2 \right).
\]
The assumption that $\renyiMarkov \ne 0$ now implies that
$D_{1/2}(f,g) = (1+o(1)) \renyiMarkov \rho$. The claims (i)--(iii) now follow Theorem~\ref{cor:recovery_conditions}.

Let us now impose the further extra assumptions of the theorem.
In this case may fix a constant $M \ge 1$ such that
$M^{-1} \le \frac{\mu_1}{\nu_1}, \frac{P_{01}}{Q_{01}},
\frac{P_{10}}{Q_{10}} \le M$. 
Moreover, the assumption $\gamma >0$ implies that $P_{11}$ and $Q_{11}$ cannot both go to one. Thus, we may choose a $\beta \in [0,1]$ such that $P_{11}^{3/2} Q_{11}^{\beta - 3/2} \not= 1 +o(1)$. Denote $\Lambda = P_{11}^{3/2} Q_{11}^{-1/2}$. Because $P_{11} \asymp Q_{11}$, we find that $\Lambda \lesim 1$.
Proposition~\ref{prop:renyi_divergence_bound_order_more_than_one_markov_chains} then implies that
$D_{3/2}(f\|g) \lesim \rho$. A similar argument shows that $D_{3/2}(g\|f) \lesim \rho$ as well. Therefore, $\frac{\dren_{3/2}^s(f,g)}{\drenh(f,g)} \lesim 1$.
Taylor's approximations further show that
the intra- and inter-block probabilities
$p = 1-(1-\mu_1)(1-P_{01})^{T-1}$ and
$q = 1-(1-\nu_1)(1-Q_{01})^{T-1}$ of observing a nonzero interaction pattern satisfy
$p = (u + (T-1) p_{01}) \rho + o(\rho)$ and
$q = (v + (T-1)q_{01}) \rho + o(\rho)$.
It follows that
$p,q \asymp \rho$ and $p-q \asymp \rho$.
When we assume that $\rho \gg \frac{1}{N}$, it follows that $p \vee q \gg N^{-1}$ and 
$\frac{(p-q)^2}{p \vee q} \asymp \rho$. 
We will apply Theorem~\ref{thm:general_algo_consistency} to conclude that 
Algorithm~\ref{algo:likelihood_based_algo} is consistent when $\rho \gg \frac{1}{N}$, and strongly consistent when $\rho \gg \frac{\log N}{N}$, or
when $\rho = \frac{\log N}{N}$ and $\renyiMarkov > K$.
\end{proof}

\begin{remark}
\label{second_remark:threshold_finite_T}
Theorem~\ref{the:ShortMarkovSBM} shows that the critical network density for strong consistency is $\rho = \frac{\log N}{N}$. In this regime, the existence of a strongly consistent estimator is determined by $\renyiMarkov$ defined in \eqref{eq:iShortMarkov}. The first term of $\renyiMarkov$ equals $(\sqrt{u}-\sqrt{v})^2$ and accounts for the first snapshot: for $T=1$ we recover the known threshold for strong consistency in the binary SBM~\cite{Abbe_Bandeira_Hall_2016, Mossel_Neeman_Sly_2016}. Each additional snapshot adds to $\renyiMarkov$ an extra term of size $\renyiMarkov_t$ bounded~by
\[
 (\sqrt{p_{01}}-\sqrt{q_{01}})^2
 + 2 c_1 h_{11}^2
 \wle \renyiMarkov_t
 \wle (\sqrt{p_{01}}-\sqrt{q_{01}})^2
 + 2 c_2 h_{11}^2
\]
with
$c_1 = \min\{ \sqrt{p_{01} q_{01}}, \gamma \sqrt{uv}\}$
and 
$c_2 = \max\{ \sqrt{p_{01} q_{01}}, \gamma \sqrt{uv}\}$.
The extra term is zero when $p_{01}=q_{01}$ and $h_{11}=0$. Notably, if the left side above is nonzero, then there exists a finite threshold $T^*$ such that strong consistency is possible for $T \ge T^*$. We illustrate this phase transition numerically in Section~\ref{subsection:experiments_temporalMarkovInteractions}.
\end{remark}

\begin{remark}
\label{exa:threshold_iid_outside_community}
In a special case of \eqref{eq:SparseMCShortParam} with $p_{01}=u$, $q_{01}=v$, $h_{11}=0$, and $\gamma=1$, the critical information quantity in \eqref{eq:iShortMarkov} equals $\renyiMarkov = T (\sqrt{u}-\sqrt{v})^2$. This coincides with multiplex networks composed of $T$ independent layers studied in Example~\ref{exa:MultiplexSBM}. This is also what we would obtain when studying transition matrices $P$ and $Q$ corresponding to independent Bernoulli sequences with means $\mu_1 = u\rho + o(\rho)$ and $\nu_1 = v\rho + o(\rho)$, because in this case 
$P_{11}=\mu_1$ and $Q_{11}=\nu_1$, leading to
$\Gamma = 1-\sqrt{P_{11}Q_{11}} = 1-O(\rho)$ and $H_{11} = \Hel(\Geo(\mu_1), \Geo(\nu_1)) = O(\rho)$.
\end{remark}


\subsubsection{Long time horizon}

Consider a Markov SBM with $T \gg 1$ snapshots in which
\begin{equation}
 \label{eq:SparseMCLongParam}
  P_{01} \weq p_{01} \rho + o(\rho), \qquad 
  Q_{01} \weq q_{01} \rho + o(\rho), \qquad
  H_{11} \weq h_{11} + o(1),
\end{equation}
for some constants $p_{01}, q_{01}, h_{11}$, and define
\begin{equation}
 \label{eq:iLongMarkov}
 \renyiMarkov
 \weq
 (\sqrt{p_{01}}-\sqrt{q_{01}})^2
 + 2 h_{11}^2 \sqrt{p_{01} q_{01}}.
\end{equation}
In the following result we assume that the effective spectral gap $\Gamma = 1-\sqrt{P_{11}Q_{11}}$ satisfies $\Gamma \gg T^{-1}$ which guarantees that both Markov chains mix fast enough, and we may ignore the role of initial states.

\begin{theorem}
\label{the:LongMarkovSBM}
Consider a Markov SBM with $N \gg 1$ nodes, $K = O(1)$ blocks, $T \gg 1$ snapshots, and assume that $\mu_1,\nu_1 \lesim \rho$ and \eqref{eq:SparseMCLongParam} holds for some constants $p_{01}, q_{01}, h_{11} \ge 0$ such that $\renyiMarkov \ne 0$. Assume also that $\rho \ll T^{-1} \ll 1-\sqrt{P_{11}Q_{11}}$. Then:
\begin{enumerate}[(i)]
\item a consistent estimator does not exist for
$\rho \lesim \frac{1}{NT}$ and does exist for $\rho \gg \frac{1}{NT}$;
\item a strongly consistent estimator does not exist for $\rho \ll \frac{\log N}{NT}$ and does exist for 
$\rho \gg \frac{\log N}{NT}$;
\item in a critical regime with $\rho = \frac{\log N}{NT}$, a strongly consistent estimator does not exist for 
$\renyiMarkov < K$ and does exist for $\renyiMarkov > K$.

\end{enumerate}
If we further assume that
$p_{01},q_{01} > 0$ and $p_{01} \ne q_{01}$, $\mu_1 \asymp \nu_1$,
$P_{10} \asymp Q_{10}$,
and that
\begin{equation}
 \label{eq:LeafPQ}
 (1+\Omega(1)) P_{11}^3
 \wle Q_{11}
 \wle (1-\Omega(1)) P_{11}^{1/3},
\end{equation}
then Algorithm~\ref{algo:likelihood_based_algo} applied with $\cA = \{0\}$ is consistent when $\rho \gg \frac{1}{NT}$; and strongly consistent when $\rho \gg \frac{\log N}{NT}$, or when $\rho = (1+o(1)) \tau \frac{\log N}{NT}$ for some constant $\tau$ and $\tau \renyiMarkov > K$.
\end{theorem}

\begin{proof}
By Proposition~\ref{thm:renyi_divergence_order_less_1_sparse_markov_chains}, we find that
\[
 D_{1/2}(f,g)
 \weq (1+o(1)) \renyiMarkov \rho T + 2 \Big(\Gamma\sqrt{\mu_1\nu_1} - \sqrt{P_{01}Q_{01}} \Big) H_{11}^2 \Gamma_T
 + O( (\rho T)^2),
\]
where $\Gamma_T = \sum_{t=0}^{T-2} (1-\Gamma)^t$.
Because $\Gamma_T \le \Gamma^{-1}$ and $H_{11} \le 1$, we see that the middle term on the right is bounded in absolute value by $2 \Gamma^{-1} \rho$. The assumption that $\rho T \ll 1 \ll \Gamma T$, combined with the assumption that $\renyiMarkov \ne 0$, now implies that
$D_{1/2}(f,g) = (1+o(1)) \renyiMarkov \rho T$. The claims (i)--(iii) now follow from Theorem~\ref{cor:recovery_conditions}.

Let us now impose the extra assumptions that
$p_{01},q_{01} > 0$ and $p_{01} \ne q_{01}$, $\mu_1 \asymp \nu_1$,
$P_{10} \asymp Q_{10}$, and \eqref{eq:LeafPQ}.
In this case may fix a constant $M \ge 1$ such that
$M^{-1} \le \frac{\mu_1}{\nu_1}, \frac{P_{01}}{Q_{01}},
\frac{P_{10}}{Q_{10}} \le M$.
Furthermore, \eqref{eq:LeafPQ} implies that $\Lambda \le 1-\Omega(1)$. Proposition~\ref{prop:renyi_divergence_bound_order_more_than_one_markov_chains} then implies that 
\[
 D_\alpha(f || g)
 \wle 8 C \rho T e^{5 C \rho T}
 \qquad \text{with $C = \frac{M^3}{1-\Lambda}$}.
\]
Because $C \lesim 1$ and $\rho T \ll 1$, we conclude that $D_\alpha(f \| g) \lesim \rho T$. A similar argument shows that
$D_{3/2}(g\|f) \lesim \rho T$ as well. Therefore, $\frac{\dren_{3/2}^s(f,g)}{\drenh(f,g)} \lesim 1$.
Taylor's approximations further show that
the intra- and inter-block probabilities
$p = 1-(1-\mu_1)(1-P_{01})^{T-1}$ and
$q = 1-(1-\nu_1)(1-Q_{01})^{T-1}$ of observing a nonzero interaction pattern satisfy
$p = p_{01} \rho T + o(\rho T)$ and
$q = q_{01} \rho T + o(\rho T)$.
It follows that
$p,q \asymp \rho T$ and $p-q \asymp \rho T$.
When we assume that $\rho \gg \frac{1}{NT}$, it follows that $p \vee q \gg N^{-1}$ and 
$\frac{(p-q)^2}{p \vee q} \asymp \rho T$. 
We will apply Theorem~\ref{thm:general_algo_consistency} to conclude that 
Algorithm~\ref{algo:likelihood_based_algo} 
applied with $\cA = \{0\}$ is consistent when $\rho \gg \frac{1}{NT}$, and strongly consistent when $\rho \gg \frac{\log N}{N T}$, or when $\rho = \frac{\log N}{NT}$ and $\renyiMarkov > K$.
\end{proof}

\begin{remark}
\label{rem:consistency_in-constantDegreeRegime}
Theorem~\ref{the:LongMarkovSBM} shows that consistent recovery may be possible even in cases where individual snapshots are very sparse, for example in regimes with $\rho \asymp \frac{1}{N}$ and $T \gg 1$. 
This is in stark contrast with standard binary SBMs, where in the constant-degree regime with $\rho \asymp \frac{1}{N}$, the best one can achieve is detection~\cite{Massoulie_2014,Mossel_Neeman_Sly_2015,Mossel_Neeman_Sly_2018}.
Similarly, when $\rho = \frac{1}{N}$ and $T = \tau \log N$, strong consistency is possible if $\tau > K \renyiMarkov^{-1}$.
\end{remark}

\begin{remark}
The conditions in Theorem~\ref{the:LongMarkovSBM} are similar to those derived for an integer-valued SBM with zero-inflated geometrically distributed interactions. Indeed, the critical quantity $\tilde{I}$ in \eqref{eq:iLongMarkov}
%
corresponds (up to second-order terms) to the \Renyi divergence between two zero-inflated geometric distributions (equation~\eqref{eq:renyi_divergence_zero_inflated_distributions}).
\end{remark}

\begin{example}[Markov SBM with persistence parameter]
\label{exa:SBMPersistence}
A temporal network model in \cite{Barucca_Lillo_Mazzarisi_Tantari_2018} is characterised by link density $\rho \ll 1$ and parameters $0 \le a, \xi, \eta \le 1$ corresponding to assortativity, link persistence, and community persistence. For $\eta=1$, the model corresponds to a Markov SBM with intra- and inter-block node pairs interacting according to stationary Markov chains having transition matrices
$P = \xi
\left[\begin{smallmatrix}
1 & \ 0 \\
0 & \ 1 \\
\end{smallmatrix}\right]
+ (1-\xi)
\left[\begin{smallmatrix}
1-\mu_1 & \ \mu_1 \\
1-\mu_1 & \ \mu_1 \\
\end{smallmatrix}\right]
$
and
$Q = \xi
\left[\begin{smallmatrix}
1 & \ 0 \\
0 & \ 1 \\
\end{smallmatrix}\right]
+ (1-\xi)
\left[\begin{smallmatrix}
1-\nu_1 & \ \nu_1 \\
1-\nu_1 & \ \nu_1 \\
\end{smallmatrix}\right]
$
and marginal link probabilities
$\mu_1 = (1-a+Ka)\rho$ and
$\nu_1 = (1-a)\rho$, respectively.
When $K=2$ and $0 < a,\xi \le 1$ are constants, conditions \eqref{eq:SparseMCLongParam} are valid with $p_{01}=(1-\xi)(1+a)$,
$q_{01}=(1-\xi)(1-a)$, and
$h_{11}=0$, 
and the critical information quantity in \eqref{eq:iLongMarkov} equals
\begin{equation}
\label{eq:iMarkovBarucca}
 \renyiMarkov
 \weq 2 (1-\xi) \left(1-\sqrt{(1-a)(1+a)}\right).
\end{equation}
By Theorem~\ref{the:LongMarkovSBM}, strong consistency in the critical regime with $\rho = \frac{\log N}{NT}$ is possible for $\renyiMarkov > 2$ and impossible for $\renyiMarkov < 2$. Formula \eqref{eq:iMarkovBarucca} quantifies how higher link persistence $\xi$ makes community recovery harder, whereas higher assortativity $a$ makes it easier. The model in \cite{Barucca_Lillo_Mazzarisi_Tantari_2018} assumes that intra-block and inter-block links have equal persistence $\xi$, leading to $h_{11}=0$. 
\end{example}

\subsection{Online algorithms}
\label{section:online_likelihood_algos}

\subsubsection{Known interaction parameters}

Given $X^{1:t} =(X^{1}, \dots, X^{t})$, we define a log-likelihood ratio matrix by
\begin{equation}
\label{eq:CumLogLikelihoodRatio}
M^{(t)}_{ij}
\weq \log \frac{ f }{ g } \left( X^{1:t}_{ij} \right).
\end{equation}
Then the log of the probability of observing a graph sequence $X^{1:t}$ given node labelling $\sigma$ equals $ \frac12 \sum_i \sum_{j \ne i} M^{(t)}_{ij} \delta_{\sigma_j, \sigma_i} + \frac12 \sum_i \sum_{j \ne i} g(X_{ij}^{1:t})$.
Therefore, given an assignment $\hsigma^{(t-1)}$ computed from the observation of the $t-1$ first snapshots, one can compute a new assignment~$\hsigma^{(t)}$ such that node~$i$ is assigned to any block~$k$ which maximizes
\begin{eqnarray}
L_{i,k}^{(t)} = \sum_{j \ne i} M_{ij}^{(t)} \delta_{\hsigma_j^{(t-1)} k }.
\end{eqnarray}
	
This formula brings computational benefits only if the computation of $M^{(t)}$ can be easily done from $M^{(t-1)}$. This is in particular the case of the Markov evolution. 
Indeed, if $\mu$ and $\nu$ are the initial probability distributions, and $P, Q$ are the transition matrices, then the cumulative log-likelihood matrices defined in equation~\eqref{eq:CumLogLikelihoodRatio} can be computed recursively by $M^{(t)} = M^{(t-1)} + \Delta^{(t)}$ with
$
M^{(1)}_{ij} = \log \frac{ \mu }{ \nu } \left( X^{1}_{ij} \right)
$
and 
$
\Delta^{(t)}_{ij} = \log \frac{ P }{ Q }\left( X^{t-1}_{ij}, X^{t}_{ij} \right)
$.
We summarise this in Algorithm~\ref{alg:MarkovClustering}. 

\begin{algorithm}[!ht]
    \DontPrintSemicolon
	\KwInput{Observed interaction array $(X_{ij}^t)$; dynamic block interaction parameters $\mu, \nu, P, Q$; number of communities $K$; static graph clustering algorithm \texttt{algo}.
	}
		
	\KwOutput{Node labelling $\hsigma = \left( \hsigma_1,\dots, \hsigma_N \right) \in [N]^K$.
	}
	~\\
	\KwInitializationStep{Compute $\hsigma \leftarrow \texttt{algo}(X^1)$, and
	$M_{ij} \leftarrow \log \frac{\mu \left( X^{1}_{ij} \right) }{\nu \left( X^{1}_{ij} \right)}$ for $i,j=1,\dots,N$. 
	}
	~\\
	\For{$t=2$, \dots, $T$}
	{   Compute $\Delta_{ij} \leftarrow \log \frac{P\left( X^{t-1}_{ij}, X^{t}_{ij} \right)}{Q \left(X^{t-1}_{ij}, X^{t}_{ij} \right)}$
	for $i,j=1,\dots,N$.
	\\
	Update $M \leftarrow M + \Delta$. \\
	\For{$i=1,\dots,N$}
	{Set $L_{ik} \leftarrow \sum_{j \ne i} M_{ij} \, \delta_{\hsigma_jk}$ for $k=1,\dots,K$. \\
	Set $\hsigma_i \leftarrow \argmax_{1 \le k \le K} L_{ik}$.
	}
    }
    ~\\
	\KwReturn{ $\hsigma$ }
\caption{\label{alg:MarkovClustering} Online clustering with known interaction parameters.}
\end{algorithm}
	
The time complexity of Algorithm~\ref{alg:MarkovClustering} is $O(K N^2 T)$ plus the time complexity of the initial clustering step. The space complexity is $O(N^2)$. Algorithm~\ref{alg:MarkovClustering} can be optimised in the following ways: 
\begin{enumerate}[(i)]
\item Since at each time step, $\Delta$ can take only one of four values, these four different values of $\Delta$ can be precomputed and stored to avoid computing $N^2T$ logarithms.
\item The $N$-by-$K$ matrix $(L_{ik})$ can be computed as a matrix product $L = M^0 \Sigma$, where $M^0$ is the matrix obtained by zeroing out the diagonal of $M$, and $\Sigma$ is the one-hot representation of $\hsigma$ such that $\Sigma_{ik} = 1$ if $\hsigma_i = k$ and zero otherwise.
\item For sparse networks, the time and space complexity can be reduced by a factor of $d/N$ where $d$ is the average node degree in a single snapshot, by neglecting the $0 \to 0$ transitions and only storing nonzero entries (similarly to what is often done for belief propagation in the static SBM \cite{Moore_2017}).
\end{enumerate}

\subsubsection{Unknown interaction parameters}
\label{subsection:unknown_parameters}
	
Algorithm~\ref{alg:MarkovClustering} requires a priori knowledge of the block interaction parameters.  This is often not the case in practice, and one has to learn the parameters during the process of recovering communities \cite{Billingsley_1961,Pensky_2019}. In this section, we adapt Algorithm~\ref{alg:MarkovClustering} to estimate the parameters on the fly.
	
Let~${n}_{ab}(i,j)$ be the observed number of $a \to b$ transitions in the interaction pattern between nodes $i$ and $j$, and let $n_a(i,j) = \sum_{b} n_{ab}(i,j)$.  
Let $P(i,j)$ be the 2-by-2 transition probability matrix for the interaction pattern between node pair $\{i,j\}$. 
By the law of large numbers (for stationary and ergodic random processes), the empirical transition probabilities
\[
 \hP_{ab}(i,j)
 \weq \frac{{n}_{ab}(i,j)}{ n_{a}(i,j)}
\]
are with high probability close to $P(i,j)$ for $T \gg 1$.

An estimator of the intra-block transition matrix $P$ is obtained by averaging those probabilities over the pairs of nodes predicted to belong to the same community. 
More precisely, after~$t$ observed snapshots ($t \ge 2$), given a predicted community assignment~$\hsigma^{(t)}$, we define for $a,b \in \{0,1\}$,
\begin{equation}
\label{eq:estimator_Pin}
\hP^{(t)}_{ab}
\weq \frac{1}{ |\{ (i,j) \ : \ \hsigma^{(t)}_i = \hsigma^{(t)}_j \} |  } \sum_{(i,j) \ : \ \hsigma^{(t)}_i = \hsigma^{(t)}_j} \dfrac{n_{ab}^{(t)}(i,j) }{n_{a}^{(t) }(i,j)  },
\end{equation}
where
\begin{align*}
 n_{ab}^{(t)}(i,j)
 \weq \sum_{t'=1}^{t-1} 1\big( X_{ij}^{t'} = a \big)  1\big( X_{ij}^{t'+1} = b \big)
\end{align*}
is the number of $a \to b$ transitions in the interaction pattern between nodes $i$ and $j$ (with $a, b \in \{0,1\}$) seen during the $t$ first snapshots, and 
$
n_{a}^{(t)} (i,j) = \sum_{b = 0}^1 n_{ab}^{(t)}(i,j).
$
Similarly, 
\begin{align}
 \label{eq:estimator_Pout}
 \hQ^{(t)} _{ab} & 
 \weq \frac{1}{ |\{ (i,j) \ : \ \hsigma^{(t)}_i \not= \hsigma^{(t)}_j \} |  } \sum_{(i,j) \ : \ \hsigma^{(t)}_i \not= \hsigma^{(t)}_j} \dfrac{n_{ab}^{(t)}(i,j) }{n_{a}^{(t)} (i,j) },
\end{align}
is an estimator of $Q_{ab}$. 
Moreover, the quantities $n_{a,b}^{(t)}(i,j)$ can be updated recursively according to
\begin{align}
 \label{eq:estimator_number_transitions}
 n_{ab}^{(t+1)}(i,j)
 = n_{ab}^{(t)}(i,j) + 1\big( X_{ij}^{t} = a \big) \, 1\big( X_{ij}^{t+1} = b \big).
\end{align}

This leads to Algorithm~\ref{alg:MarkovClustering_parameters_unknown} for clustering a Markov SBM when only the number of communities $K$ is known. Note that to save computation time, we can choose not to update the parameters at each time step.

\begin{algorithm}[!ht]
\DontPrintSemicolon
\KwInput{Observed graph sequence $X^{1:T} = \left( X^{1}, \dots, X^{T} \right)$; number of communities $K$; static graph clustering algorithm \texttt{algo}.}
	
\KwOutput{Node labelling $\hsigma = \left( \hsigma_1,\dots, \hsigma_n \right)$.}
~\\
	
\KwInitializationStep{
\begin{itemize}
 \item Compute $\hsigma \leftarrow \texttt{algo} \left( X^1 \right) $;
 \item Set $n_{ab}(i,j) \leftarrow 0$ for $i,j \in [N]$ and $a,b \in \{0,1\}$.
\end{itemize}
}
\KwUpdateStep{}
 \For{$t=2,\cdots,T$}
 {
   For every node pair $(ij)$, update $n_{ab}(i,j)$ using  \eqref{eq:estimator_number_transitions}; \\
  Compute $\hP, \hQ$ using \eqref{eq:estimator_Pin} and \eqref{eq:estimator_Pout}; \\
  Compute $M$ such that $M_{ij} = \sum_{a,b} n_{ab}(i,j) \log \frac{ \hP_{ab} }{ \hQ_{ab} }$. \\
  \For{$i=1, \dots, N$}
  {
   Set $L_{i,k} \leftarrow \sum_{j \ne i} M_{ij} 1 \left( \hsigma_j = k \right)$ for all $k=1,\dots,K$ \\
   Set $\hsigma_i \leftarrow \argmax_{1 \le k \le K} L_{i,k}$
  }
  }
\caption{\label{alg:MarkovClustering_parameters_unknown} Online clustering with unknown interaction parameters.}
\end{algorithm}

\section{Numerical experiments}
\label{section:numerical_results}
This section presents numerical experiments of the different algorithms presented in this paper\footnote{Source code for the algorithms is available at \\ \url{https://github.com/mdreveton/clusteringNonBinaryAndTemporalSBM}.}.

\subsection{Static networks with numerical interactions}

Let us study the performance of Algorithm~\ref{algo:likelihood_based_algo} on synthetic data sampled from real-valued and nonnegative integer-valued SBMs. As input to the algorithm, the set $\cA$ is chosen as a continuous interval 
$[-x,x]$ (real-valued interaction space) or a set 
$\{0,\dots,x\}$ (nonnegative integer-valued SBM), so that the Hellinger distance $\cA \mapsto \Hel \left( \Ber(f(\cA)), \Ber(g(\cA)) \right)$ is maximised.

\begin{figure}[!ht]
\centering
\captionsetup[subfigure]{justification=centering}
\begin{subfigure}[b]{0.45\textwidth}
\includegraphics[width=\textwidth]{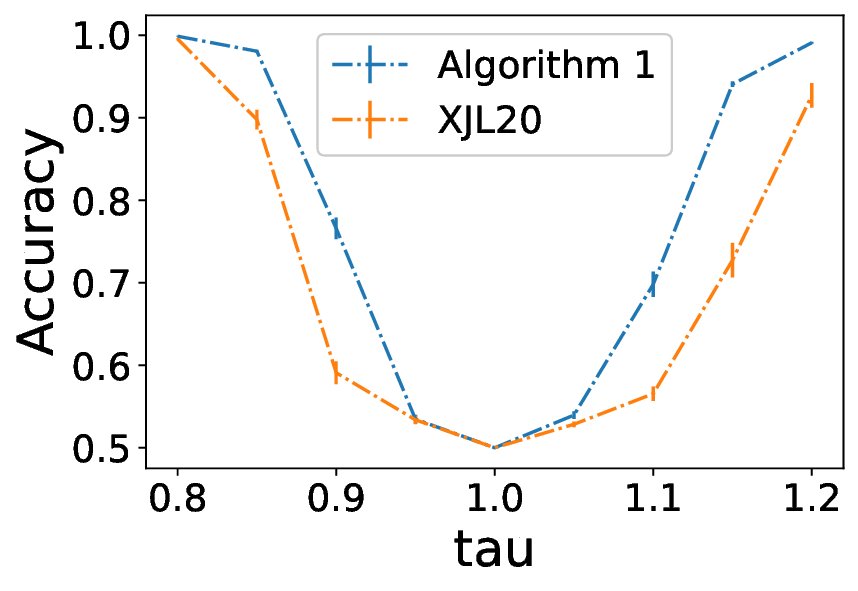}
\caption{$f = \Nor(0,1)$, $g = \Nor(0,\tau)$.}
\end{subfigure}
\hspace{3mm}
\begin{subfigure}[b]{0.45\textwidth}
\includegraphics[width=\textwidth]{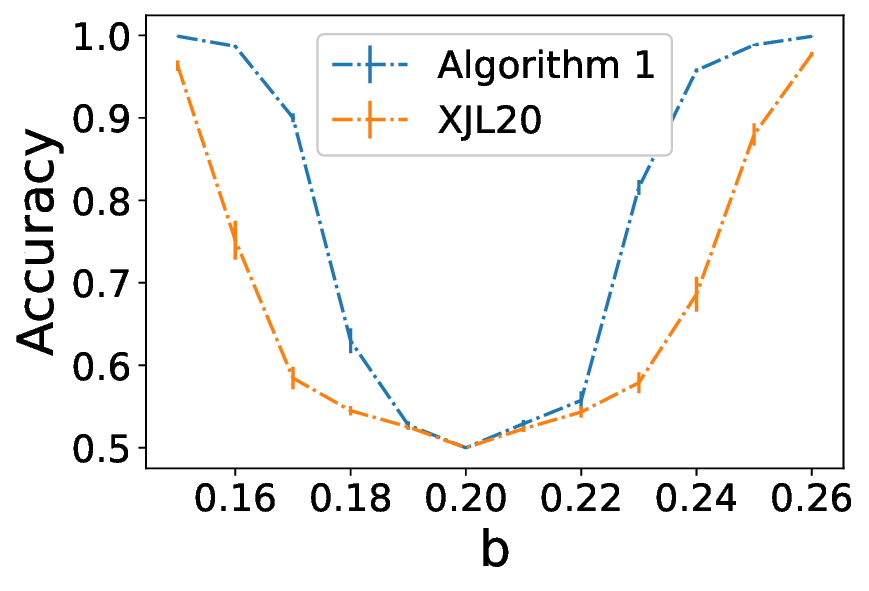}
\caption{$f = \Geo(0.2)$, $g = \Geo(b)$.}
\end{subfigure}
\caption{
Accuracy of Algorithm~\ref{algo:likelihood_based_algo} and the algorithm in~\cite{Xu_Jog_Loh_2020} on data sampled from numerical SBMs of $N = 400$ nodes and $K=2$ blocks; with (a) normal and (b) geometric interaction distributions. Results are averaged over $25$ samples. Error bars display empirical standard deviations.
}
\label{fig:experiments_dense_networks}
\end{figure}

\begin{figure}
\centering
\captionsetup[subfigure]{justification=centering}
\begin{subfigure}[b]{0.45\textwidth}
\includegraphics[width=\textwidth]{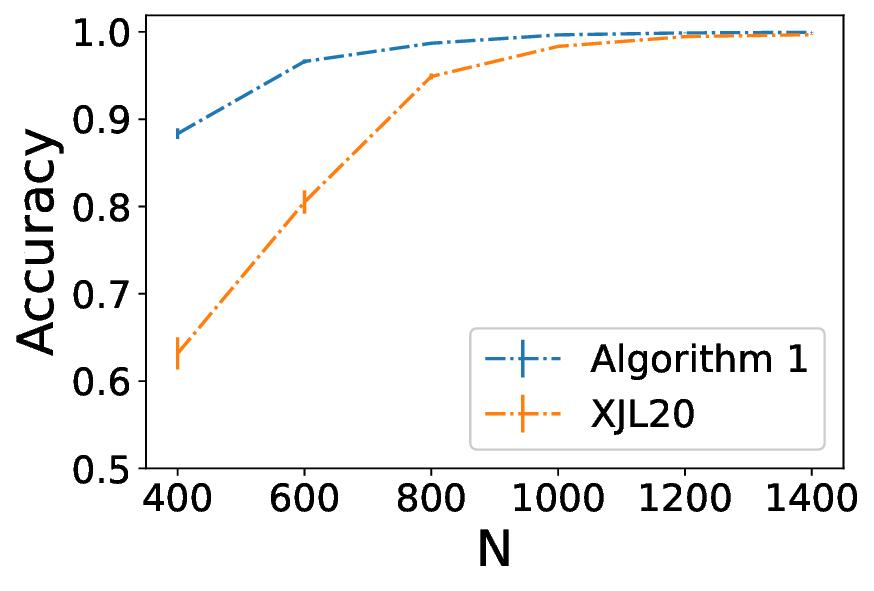}
\caption{$\tf = \Nor(0,1)$, \\ \hspace{7.0mm} $\tg = \Nor(0,1.2)$.}
\end{subfigure}
\hspace{1mm}
\begin{subfigure}[b]{0.45\textwidth}
\includegraphics[width=\textwidth]{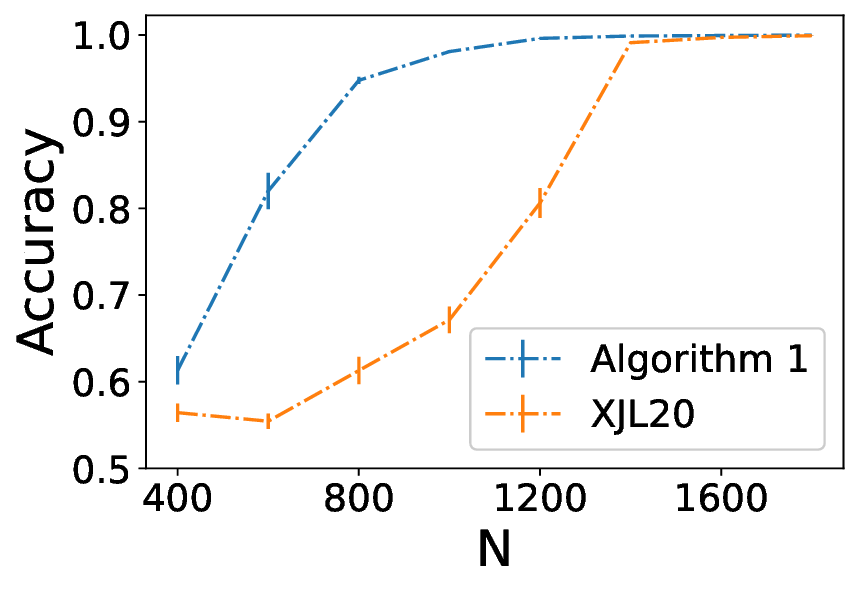}
\caption{$\tf = \Nor(0, 1.3^2+1 )$, \\ $\tg = \frac12 \Nor(-1.3,1) + \frac12 \Nor(1.3,1)$.}
\label{fig:theoreticalAlgo_variation_with_N_sameAsXJL20}
\end{subfigure}
\caption{
Accuracy of Algorithm~\ref{algo:likelihood_based_algo} and the algorithm in~\cite{Xu_Jog_Loh_2020}
as a function of the number of nodes $N$ for $K=2$ blocks. Data are sampled from zero-inflated interaction distributions
$f = \frac12 \delta_0 + \frac12 \tf$
and 
$g = \frac12 \delta_0 + \frac12 \tg$ with 
(a) normal, (b) mixed normal distributions $\tf,\tg$.
Results are averaged over $25$ samples.
}

\label{fig:theoreticalAlgo_variation_with_N}
\end{figure}

We compare the performance of Algorithm~\ref{algo:likelihood_based_algo} with the algorithm in~\cite{Xu_Jog_Loh_2020} which to best of our knowledge is the only other clustering algorithm that works both with discrete and continuous edge labels. 
For a fair comparison, we implemented a version of the algorithm in~\cite{Xu_Jog_Loh_2020} in which the interaction distributions are given as input.
Figure~\ref{fig:experiments_dense_networks} compares the accuracy\footnote{
We define accuracy as the proportion of correctly labeled nodes $1-N^{-1} \dhams (\sigma_1, \sigma_2)$.}
of the algorithms on networks with normal (Example~\ref{exa:NormalSBMAlgo}) and geometric (Example~\ref{exa:GeometricSBMAlgo}) interaction distributions. Figure~\ref{fig:theoreticalAlgo_variation_with_N} compares the algorithms for zero-inflated normal and mixed normal interaction distributions (cf.\ Example~\ref{exa:NormalCensoredSBM}) with parameters in  Figure~\ref{fig:theoreticalAlgo_variation_with_N_sameAsXJL20}
matching the simulation experiments in \cite[Section 7]{Xu_Jog_Loh_2020}. 
Overall, Algorithm~\ref{algo:likelihood_based_algo} achieves improved accuracy for all studied parameter combinations, with most remarkable improvements obtained in cases involving non-normal interaction distributions.
%

\subsection{Temporal networks}
\label{subsection:experiments_temporalMarkovInteractions}


We study community recovery from temporal network data sampled from a stationary Markov SBM described in Section~\ref{section:temporal_networks}.
We focus on sparse settings where the average degree per snapshot is of constant order, so that consistent recovery using a single snapshot is impossible, but consistent and even strongly consistent community recovery is possible when the number of snapshots $T$ is large enough (Remark~\ref{rem:consistency_in-constantDegreeRegime}).

\subsubsection{Offline recovery}

\begin{figure}[!ht]
\centering
\captionsetup[subfigure]{justification=centering}
\begin{subfigure}[b]{0.45\textwidth}
\includegraphics[width=\textwidth]{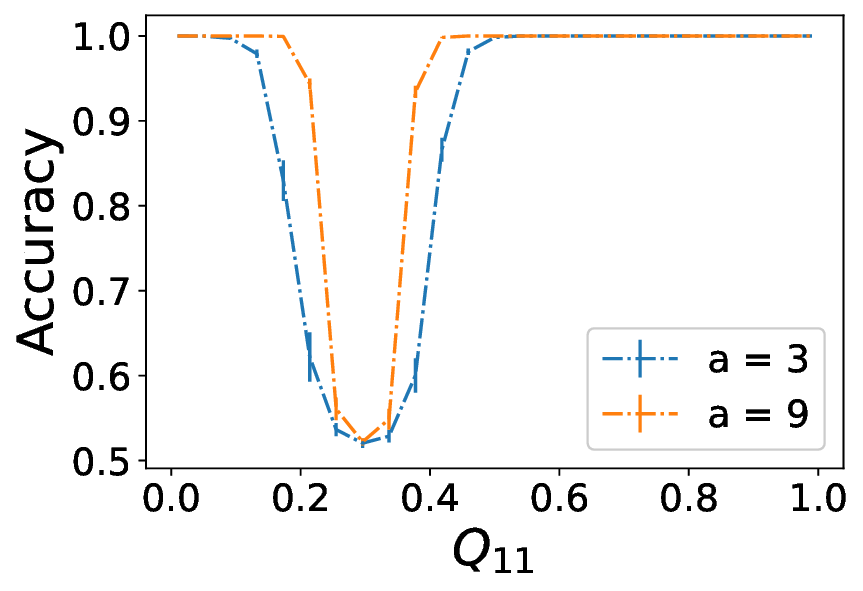}
\caption{$P_{11} = 0.3$.}
\end{subfigure}
\hspace{1mm}
\begin{subfigure}[b]{0.45\textwidth}
\includegraphics[width=\textwidth]{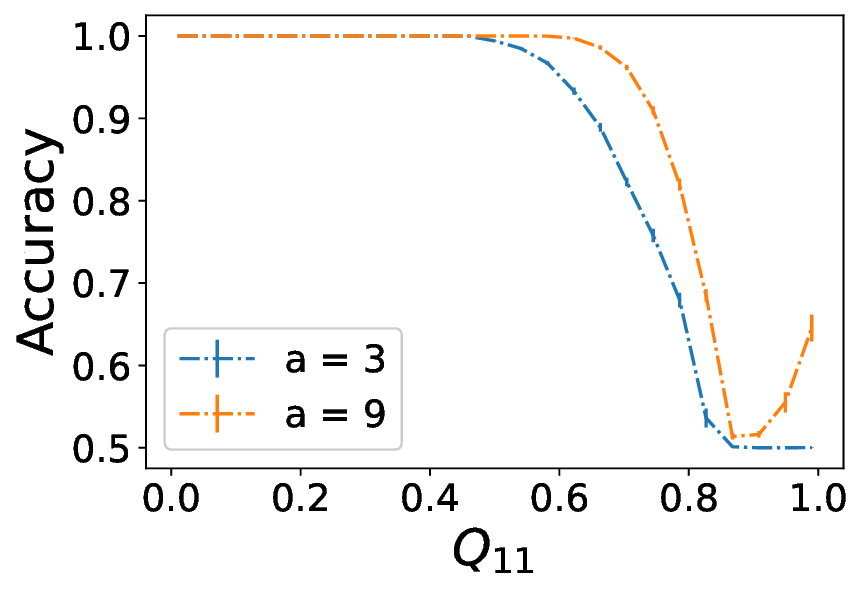}
\caption{$P_{11} = 0.9$.}
\end{subfigure}
\caption{Performance of Algorithm~\ref{algo:likelihood_based_algo} on a temporal network of $N=400$ nodes, $K=2$ blocks, and $T=60$ snapshots, as a function of inter-block link persistence $Q_{11}$. Data are sampled from a stationary Markov SBM with equal intra- and inter-block link densities $\mu_1 = \nu_1 = \frac{a}{N}$, and intra-block link persistence (a) $P_{11}=0.3$, (b) $P_{11}=0.9$.
Results are averaged over 10 samples.}
\label{fig:Algo1_MarkovSBM}
\end{figure}

In an offline setting we apply the generic Algorithm~\ref{algo:likelihood_based_algo} with vector-valued interactions initialised using $\cA = \{ x \in \{0,1\}^T \colon \log \frac{f}{g}(x) \le \log \frac{f}{g}(0) \}$.
Figure~\ref{fig:Algo1_MarkovSBM} presents the algorithm's accuracy on a stationary Markov SBM where the intra- and inter-block interactions are indistinguishable for any single snapshot ($\mu_1=\nu_1$), and communities can be identified only by the different link persistence rates $P_{11}$ and $Q_{11}$ within and between communities.
As expected, the accuracy of community recovery is low for $P_{11} \approx Q_{11}$. Outside such parameter regions, the performance of Algorithm~\ref{algo:likelihood_based_algo} is remarkably high, even though each single snapshot alone carries no information about the community structure.

\subsubsection{Online recovery with known interaction parameters}

\begin{figure}[!ht]
\centering
\captionsetup[subfigure]{justification=centering}
\begin{subfigure}{0.45\textwidth}
  \includegraphics[width=\textwidth]{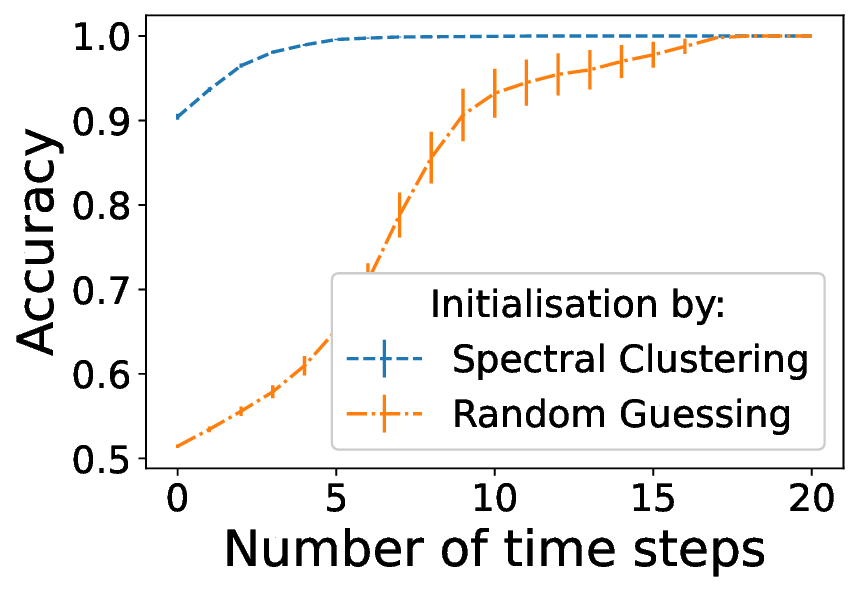}
  \caption{$\mu_1 = \frac{15}{N}$}
  \label{fig:sc_better_than_rg}
\end{subfigure}
\hspace{1mm}
\begin{subfigure}{0.45\textwidth}
  \includegraphics[width=\textwidth]{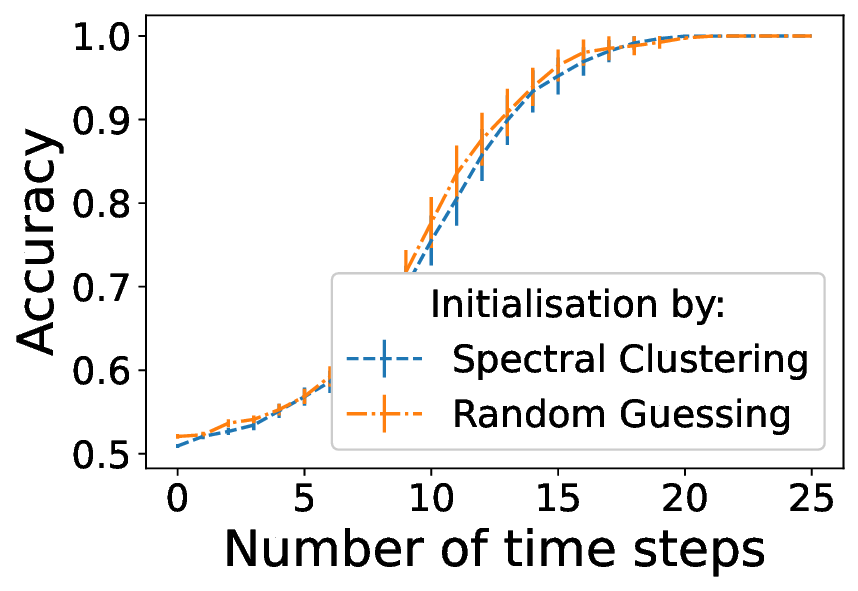}
  \caption{$\mu_1 = \frac{5}{N}$}
  \label{fig:sc_rg_similar}
\end{subfigure}
\caption{Performance of Algorithm~\ref{alg:MarkovClustering} as a function of the number of snapshots $T$ in a temporal network of $N = 500$ nodes and $K=2$ blocks.
Data are sampled from a stationary Markov SBM with intra-block link density
(a) $\mu_1=\frac{15}{N}$,
(b) $\mu_1 = \frac{5}{N}$;
inter-block link density
$\nu_1 = \frac{5}{N}$; and
intra- and inter-block link persistence parameters $P_{11}=0.7$ and $Q_{11}=0.4$. Results are averaged over 25 samples.}
\label{fig:evolution_accuracy_snapshots}
\end{figure}

Figure~\ref{fig:evolution_accuracy_snapshots} illustrates
the number of snapshots needed to recover communities accurately using 
Algorithm~\ref{alg:MarkovClustering}, initiated either by spectral clustering or a blind random guess.
In a sufficiently dense case (Figure~\ref{fig:sc_better_than_rg}), spectral clustering on the first snapshot works well, and even a blind random guess leads to accurate results after a handful of iterations. In a sparser case (Figure~\ref{fig:sc_rg_similar}), 
spectral clustering on the first snapshot performs poorly, but a few online updates rapidly improve accuracy. Remarkably in both cases, a modest number of online updates 
yields a high accuracy, regardless of the quality of the initial clustering.

\subsubsection{Online recovery with unknown interaction parameters}

\begin{figure}[!ht]
\centering
\captionsetup[subfigure]{justification=centering}
\begin{subfigure}[b]{0.45\textwidth}
\includegraphics[width=\textwidth]{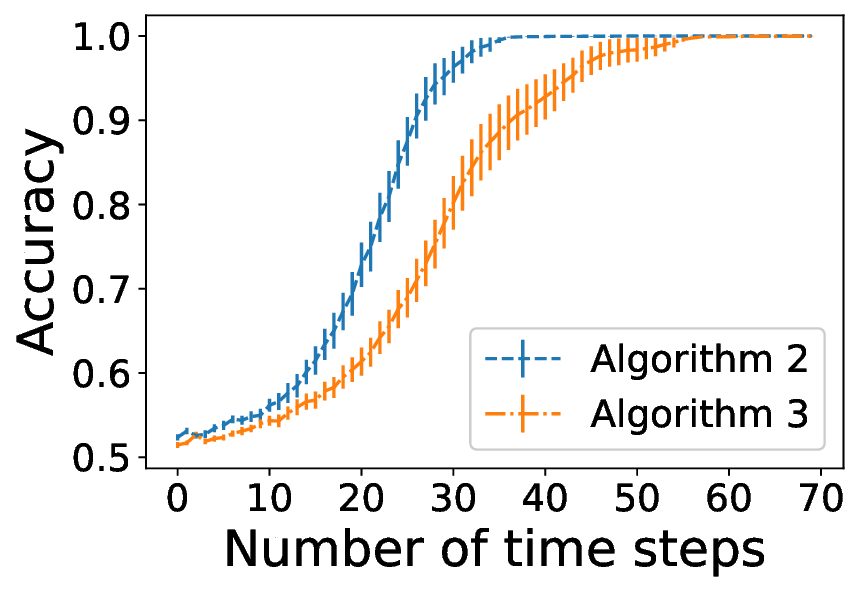}
\caption{$\mu_1 = \nu_1 = \frac{2}{N}$.}
\end{subfigure}
\hspace{1mm}
\begin{subfigure}[b]{0.45\textwidth}
\includegraphics[width=\textwidth]{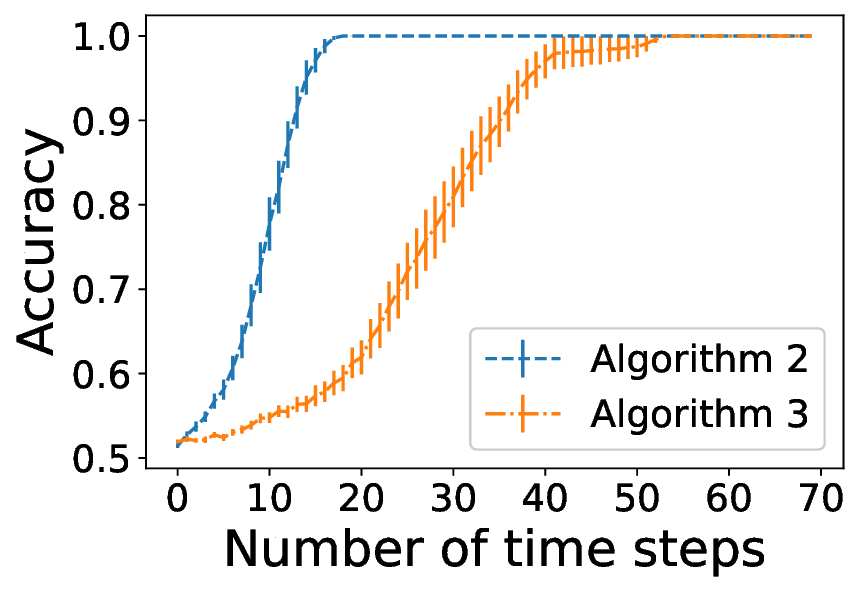}
\caption{$\mu_1 = \frac{8}{N}$, \
$\nu_1 = \frac{2}{N}$.}
\end{subfigure}
\caption{
Performance of
Algorithms~\ref{alg:MarkovClustering} and~\ref{alg:MarkovClustering_parameters_unknown} as a function of the number of snapshots $T$ in a temporal network of $N=400$ nodes and $K=2$ blocks. Data are sampled from a stationary Markov SBM with intra- and inter-block link densities (a) $\mu_1=\nu_1$,
(b) $\mu_1=4\nu_1$; and intra- and inter-block link persistence parameters $P_{11}=0.6$ and $Q_{11}=0.3$. Results are averaged over $25$ samples.}
\label{fig:comparison_accuracy_online_knowingNotKnowing}
\end{figure}

When the interaction parameters are unknown, we replace Algorithm~\ref{alg:MarkovClustering} with Algorithm~\ref{alg:MarkovClustering_parameters_unknown}, which adaptively estimates the interaction parameters jointly with community recovery. Figure~\ref{fig:comparison_accuracy_online_knowingNotKnowing} compares these algorithms in a sparse setting in which spectral clustering on a single snapshot does not provide much more information than a blind random guess. We see that a modest number of additional snapshots suffices to compensate for the need to estimate interaction parameters from data on the fly.

\subsection{Experiments on real data}

\begin{figure}[!ht]
 \centering
 \begin{subfigure}[b]{0.32\textwidth}
 \includegraphics[width=\textwidth]{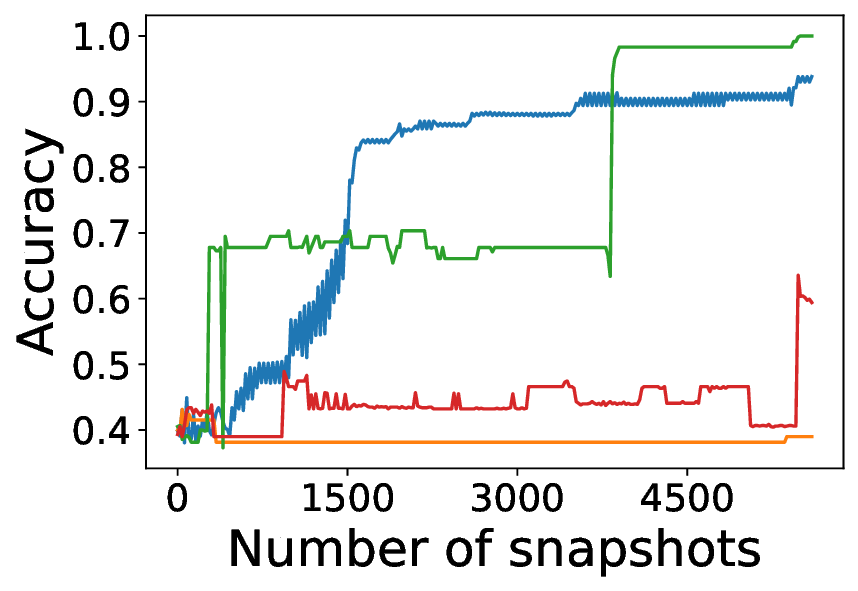}
 \caption{Year $2011$.}
 \end{subfigure}
 \hfill
 \begin{subfigure}[b]{0.32\textwidth}
 \includegraphics[width=\textwidth]{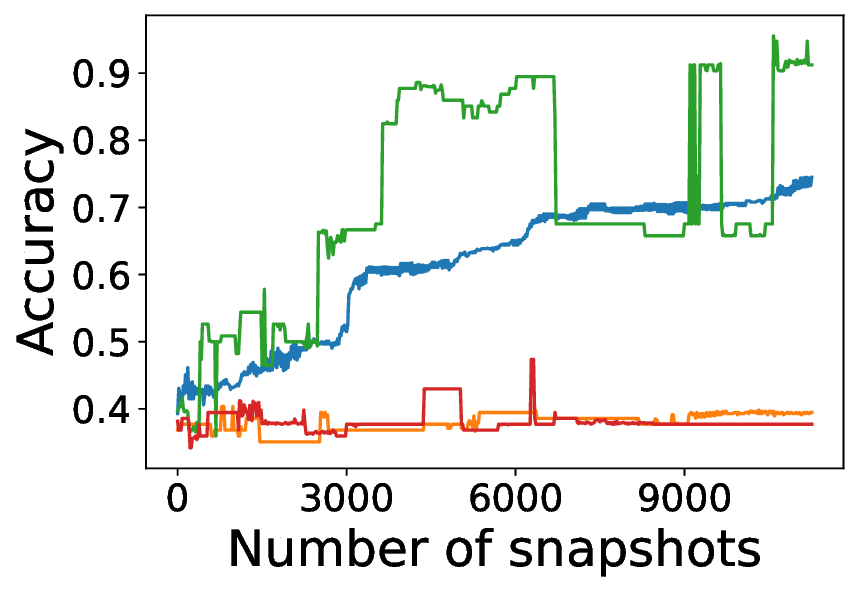}
 \caption{Year 2012.}
 \end{subfigure}
 \hfill
 \begin{subfigure}[b]{0.32\textwidth}
 \includegraphics[width=\textwidth]{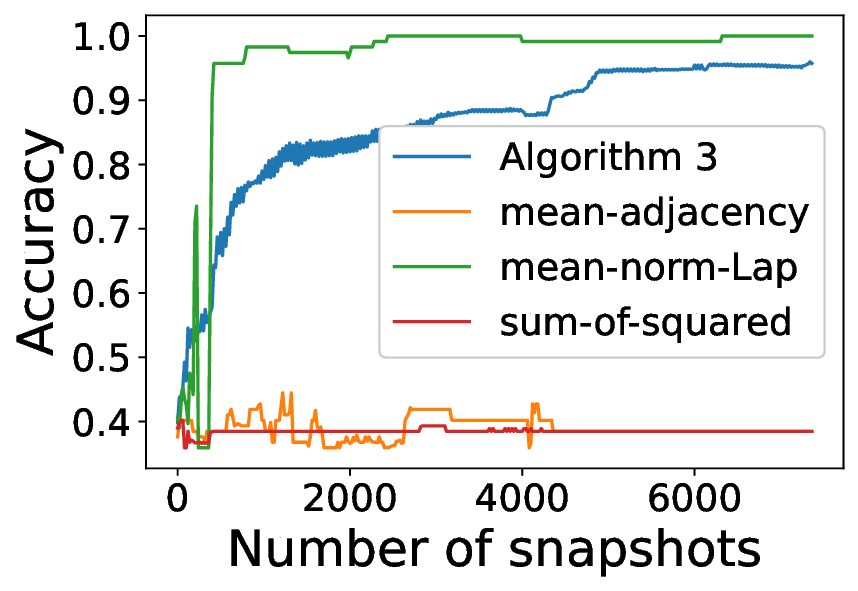}
 \caption{Year 2013.}
 \end{subfigure}
 \caption{Performance of Algorithm~\ref{alg:MarkovClustering_parameters_unknown} vs.\ three reference algorithms on high school data sets. 
 }
 \label{fig:onlinelikehihood_highschool}
\end{figure}

We investigate three data sets collected during three consecutive years from a high school Lycée Thiers in Marseilles, France~\cite{Fournet_Barrat_2014,Mastrandrea_Fournet_Barrat_2015}.
Nodes correspond to
students, interactions to close-proximity encounters, and communities to classes, with about 40 students per class. We restrict to a subset of data corresponding to $K=3$ classes labelled PC, PC$^*$ and PSI$^*$, as they are present in each of the data sets. 
The performance of Algorithm~\ref{alg:MarkovClustering_parameters_unknown} is compared against three reference algorithms:
\begin{itemize}
\item \textit{mean-adjacency} \cite{Paul_Chen_2020}, based on eigenvectors of the time-averaged adjacency matrix $\bar{X} = \frac1T \sum_{t=1}^T X^t$;
\item \textit{mean normalised Laplacian} \cite{Avrachenkov_Dreveton_Leskela_2021}, based on eigenvectors of the normalised Laplacian of $\bar{X} = \frac1T \sum_{t=1}^T X^t$;
\item \textit{sum-of-squared}  \cite{Lei_Lin_2022}, based on eigenvectors of the matrix $\sum_{t=1}^T (( X^t )^2 - D^t)$, where $D^t$ is the diagonal matrix with entries $D^t_{ii} = \sum_{j=1}^N X_{ij}^t$.
\end{itemize}
Figure~\ref{fig:onlinelikehihood_highschool} summarises the results. The \textit{mean normalised Laplacian} algorithm is highly accurate in several cases, but is prone to large fluctuations. In contrast, Algorithm~\ref{alg:MarkovClustering_parameters_unknown} displays a stable performance over time. The \textit{mean-adjacency} and \textit{sum-of-squared} algorithms perform poorly for each of the three data sets. 
We emphasise that the networks are very sparse, with the average degrees per snapshot in the three data sets being $0.04$, $0.02$, and $0.06$. This explains why a large number of snapshots is needed for community recovery.

\section{Technical comparison with related work}
\label{section:comparison_other_work}

Let us discuss our contributions with respect to the most closely related earlier works.

%
Jog and Loh \cite{Jog_Loh_2015} discovered that the Rényi divergence provides a sharp quantity for strong consistency in
homogeneous SBMs with discrete interaction distributions which are sparse in the sense of \eqref{eq:zero_inflated_distributions}; recall Example \ref{exa:SparseCategoricalSBM}.
They analysed the MLE for networks of density $\rho \asymp \frac{\log N}{N}$, assuming that the conditional probability densities $\tf,\tg$ in \eqref{eq:zero_inflated_distributions} do not depend on scale, are strictly positive either with respect to the counting measure on the positive integers or the Lebesgue measure on the real line, and are bounded by
$\supnorm{\log \frac{\tf}{\tg}} = O(1)$.
The latter condition is not satisfied for several cases of interest, for example normal distributions with equal variances but unequal means.
Part \eqref{ite:StronglyConsistent} of Theorem~\ref{the:ConsistencySparse} generalises the setting of \cite{Jog_Loh_2015} to arbitrary probability measures on an arbitrary measurable space, and does not require the condition 
$\supnorm{\log \frac{\tf}{\tg}} = O(1)$.

%
Yun and Proutière \cite{Yun_Proutiere_2016} consider interactions on a finite space and obtain consistency results similar to Theorem~\ref{cor:recovery_conditions} but with additional regularity conditions, which in the homogeneous case correspond to
$\supnorm{\log \frac{f}{g}} = O(1)$ and 
$\sum_{x \ne 0} (f(x) - g(x))^2 = \Omega(\rho^2)$
with $\rho = \max_{x \ne 0} \left( f(x), g(x) \right) \gg N^{-1}$. For the consistency of a spectral clustering algorithm, they also impose $N \min_{x \ne 0} ( f(x) \wedge g(x))
 \wge (N \rho)^{\Omega(1)}
$.
In contrast to \cite{Jog_Loh_2015,Xu_Jog_Loh_2020}, the analysis in \cite{Yun_Proutiere_2016} is valid also for inhomogeneous SBMs with unbalanced block sizes. We believe that Theorems~\ref{thm:asymptotic_exponential_bound_recovery} and~\ref{cor:recovery_conditions} could be extended to similar generality, at the cost of longer and more technical proofs to account for the lack of symmetry.

Xu, Jog, and Loh~\cite{Xu_Jog_Loh_2020} is a major contribution to the study of homogeneous SBMs with unknown interactions, but still relies on several restrictive assumptions. First, their consistency analysis is restricted to interaction distributions having an atom at zero, thereby ruling out purely continuous distributions (e.g.\ Example~\ref{exa:NormalSBM}). Also, the analysis does not extend to discrete probability distributions with infinite support (e.g.\ Examples~\ref{exa:PoissonSBM}, \ref{exa:GeometricSBMAlgo}, and \ref{exa:ZIGeometricSBMAlgo}); nor interactions distributions with finite support of size growing with $N$ (e.g.\ temporal networks with $T \gg 1$). 
%
%
Moreover, some additional technical conditions are needed,
such as the existence of two blocks of sizes $\Nmin$ and $\Nmin+1$ where $\Nmin$ is the minimum block size
(see \cite[Theorem~2]{Xu_Jog_Loh_2020}), as well as some technical smoothness conditions which may be difficult to verify in practice.
Theorems~\ref{thm:asymptotic_exponential_bound_recovery} and Theorem~\ref{cor:recovery_conditions} generalise the framework of \cite{Xu_Jog_Loh_2020}
to a setting which requires neither regularity assumptions on $\fin, \fout$ nor restrictions on the underlying space $\cS$ of interaction types.
%
Theorem~\ref{thm:general_algo_consistency} is similar in spirit to upper bounds in \cite{Xu_Jog_Loh_2020} and \cite{Yun_Proutiere_2016} which perform initial clustering using $\cA = \{0\}$, but is fundamentally different in that it makes no assumptions about truncating the label space $\cS$, nor any assumptions about the regularity of the interaction distributions $f, g$. Moreover, for temporal binary interactions with $\cS = \{0,1\}^T$, the algorithms in \cite{Xu_Jog_Loh_2020} are of exponential complexity in $T$.

Paul and Chen \cite{Paul_Chen_2016} is a key contribution on the recovery thresholds for multilayer SBMs, assuming uncorrelated layers.
Section~\ref{section:temporal_networks} contains both information-theoretic and algorithmic contributions to clustering temporally correlated networks.
Theorems~\ref{the:ShortMarkovSBM} and~\ref{the:LongMarkovSBM} extend the setting of \cite{Paul_Chen_2016} to correlated layers. In addition, the optimal misclassification rate in Theorem~\ref{cor:recovery_conditions} extends \cite[Theorem 6]{Paul_Chen_2016} to non-binary settings (recall Example \ref{exa:MultiplexSBM}). We developed Algorithms \ref{alg:MarkovClustering} and \ref{alg:MarkovClustering_parameters_unknown} for online community recovery in temporal networks.
These algorithms are designed to accurately utilise information related to temporal correlation patterns, and as such
are radically different from the mainstream of methods \cite{Bhattacharyya_Chatterjee_2018-05-27,Bhattacharyya_Chatterjee_2020-04-06,Lei_Lin_2022,Paul_Chen_2016,Paul_Chen_2020} relying on spectral clustering of layer-aggregated adjacency matrices.

\section{Conclusions and future work}
\label{section:conclusion}

In this paper, we studied community recovery in non-binary and dynamic stochastic block models.
Unlike most earlier works, our analysis allows the shape and size of the interaction space to be scale-dependent, which enables us to study correlated interaction patterns over short and long time horizons.
For clarity, most consistency results were stated under the assumption that the number of blocks is bounded, but quantitative bounds in Theorem~\ref{thm:asymptotic_exponential_bound_recovery} allow several generalisations to cases with $K \gg 1$.
We proposed Algorithm~\ref{algo:likelihood_based_algo} that fully utilises the non-binary nature of the observed data for recovering community memberships. 
%
%
Unlike earlier methods, Algorithm~\ref{algo:likelihood_based_algo} is provably consistent for general interaction distributions (not requiring atoms), including standard continuous distributions such as normal and exponential.
Our analysis of consistency essentially requires bounded \Renyi divergences of order 3/2. Investigating whether this condition can be relaxed remains an open problem.

For temporal and multiplex networks, we proposed Algorithms~\ref{alg:MarkovClustering} and \ref{alg:MarkovClustering_parameters_unknown} for community recovery based on fast online likelihood updating, and investigated their performance with numerical experiments on synthetic and real data. We observed that even in sparse or low-information regimes, both algorithms appear to produce accurate results given a reasonable number of temporal snapshots. The theoretical consistency analysis of these algorithms remains an open problem.


%% file: text/text_appendix.tex
\section{Preliminaries}

\subsection{Table of notations}

We keep the same notations as in the main text. Additionally, we define $Z_\alpha(f\|g) = \int f^\alpha g^{1-\alpha}$, so that $\dren_\alpha(f\|g) = (\alpha-1)^{-1} \log Z_\alpha(f\|g)$. We also denote by $\dkl(f\|g) = \int f \log \frac{f}{g}$ the Kullback-Leibler divergence between $f$ and $g$ and we introduce $\vkl(f\|g) = \int f \log^2 \frac{f}{g} - \dkl^2(f\|g)$.
Table \ref{tab:table of notation} summarises commonly used notations in the article.

\begin{table}[!ht]
\centering
\begin{tabular}{ll} \toprule
\bf Symbol & \bf Meaning \\ \midrule
$\delta_x$ & Dirac measure at $x$ \\
$\delta_{ab}$ & Kronecker delta \\
$\rho$ & Overall density parameter \\
$\eta$ & Scale parameter \\
\\
$N$ & Number of nodes \\
$K$ & Number of communities (blocks) \\
$T$ & Number of snapshots (temporal networks) \\
$\cS$ & Space of interaction types $(\cS = \{0,1\}^T$ for temporal networks) \\
$\cX$ & Space of observations \\
$\cZ$ & Space of node labellings (subset of $[K]^{[N]}$) \\
\\
$i, j$ & Node indices \\
$k, \ell$ & Community (block) indices\\
$\sigma$ & Node labelling ($\sigma \in [K]^N$) \\
$X = \left( X_{ij} \right)$ & Data array ($X \in \cS^{N\times N}$) \\
\\
$f(x)$, $g(x)$ & Probability of an interaction of type $x \in \cS$ between two nodes \\
$\mu_a$, $\nu_a$ & Initial intra- and inter-block interaction distributions, $a \in \{0,1\}$ (for Markov dynamics) \\
$P_{ab}$, $Q_{ab}$ & Probability of  transition $a \to b$ ($a, b \in \{0,1\})$ for intra- and inter-block interactions\\
\\
$\dren_\alpha(f \|g)$ & \Renyi divergence of order $\alpha$ \\
$\dren^s_\alpha(f,g)$ & Symmetric \Renyi divergence of order $\alpha$ \\
$Z_\alpha(f\|g)$ & Hellinger integrals \\
$\Hel(f,g)$ & Hellinger distance \\
$\beta_r(f,g)$ & \Renyi divergence ratio (defined in~\eqref{eq:alphaRatio}) \\
\\
$\dham(\sigma_1,\sigma_2)$ & Hamming distance \\
$\dhams(\sigma_1,\sigma_2)$ & Absolute classification error \\
$\Mir(\sigma_1,\sigma_2)$ & Mirkin distance \\
\bottomrule
\end{tabular}
\caption{Common notations.}
\label{tab:table of notation}
\end{table}

\subsection{Multinomial concentration}

Fix integers $N,K \ge 1$, and consider the space $[K]^N$ of mappings $\sigma: [N] \to [K]$.  For any such mapping, we denote the frequencies of output values by $N_k(\sigma) = \sum_{i=1}^N \delta_{\sigma(i)k}$ for $k=1,\dots,K$.  When the space $[K]^N$ is equipped with a probability measure $\pr$, then $\sigma \mapsto (N_1(\sigma), \dots, N_K(\sigma))$ is considered as a random variable. Given $\epsilon > 0$ and $\alpha_1,\dots,\alpha_K \in [0,1]$, we shall be interested in probabilities of events of the form
\begin{align}
 \label{eq:MultinomialConcentration}
 \cA_\epsilon
 &\weq \Big\{ \sigma: \abs{N_k(\sigma) - \alpha_k N} \le \epsilon \alpha_k N \ \text{for all $k \in [K]$} \Big\}, \\
 \label{eq:MultinomialConcentrationPlus}
 \cA_{\epsilon,+}
 &\weq \Big\{ \sigma: N_k(\sigma) \ge (1-\epsilon) \alpha_k N \ \text{for all $k \in [K]$} \Big\}.
\end{align}

\begin{lemma}
\label{the:MultinomialConcentration}
Let $0 < \epsilon \le 1$.
\begin{enumerate}[(i)]
\item If $\pr = \alpha^{\otimes N}$ for a probability measure $\alpha$ on $[K]$, then
$\pr( \cA_\epsilon^c ) \le 2 \sum_{k=1}^K e^{-(\epsilon^2/3) \alpha_k N}$
and $\pr( \cA_{\epsilon,+}^c ) \le \sum_{k=1}^K e^{-(\epsilon^2/2) \alpha_k N}$.
\item If $\pr$ is the uniform distribution on $[K]^N$, then
$\pr( \cA_\epsilon^c ) \le 2 e^{\log K - \epsilon^2 N /(3K)}$
and $\pr( \cA_{\epsilon,+}^c ) \le e^{\log K - \epsilon^2 N /(2K)}$.
\end{enumerate}
\end{lemma}
\begin{proof}
(i) Because $N_k$ is $\Bin(N,\alpha_k)$-distributed, a Chernoff bound \cite[Corollary 2.3]{Janson_Luczak_Rucinski_2000} implies that 
\[
 \pr( \abs{N_k(\sigma) - \alpha_k N} > \epsilon \alpha_k N )
 \wle 2 e^{-(\epsilon^2/3) \alpha_k N}.
\]
Similarly, another Chernoff bound \cite[Theorem 2.1]{Janson_Luczak_Rucinski_2000} implies that 
\[
 \pr( N_k(\sigma) \le (1-\epsilon) \alpha_k N )
 \wle e^{-(\epsilon^2/2) \alpha_k N}.
\]
Hence the first claim follows by the union bound.

(ii) The second claim follows immediately from (i) after noting that the uniform distribution on $[K]^N$ can be represented as $\pi = \alpha^{\otimes N}$ where $\alpha_k = K^{-1}$ for all $k$.
\end{proof}

We shall also be interested in random variables defined by $\Nmin(\sigma) = \min_k N_k(\sigma)$ and $\Delta N(\sigma) = \max_{k,\ell} \abs{N_k(\sigma) - N_\ell(\sigma)}$. The following result implies that for large-scale uniformly distributed settings with $N \gg K \log K$, these random variables are bounded by $\Nmin \ge (1-\epsilon) K^{-1}N$ and $\Delta N \le 2 \epsilon K^{-1} N$ with high probability for $(\frac{K \log K}{N})^{1/2} \ll \epsilon \le 1$. For example, we may select $\epsilon = (\frac{K \log K}{N})^{0.499}$.

\begin{lemma}
\label{the:MultinomialConcentrationMinMaxDiff}
Let $0< \epsilon \le 1$.
(i) If $\pr = \alpha^{\otimes N}$ for a probability measure $\alpha$ on $[K]$, then 
\begin{equation}
\begin{aligned}
 \label{eq:MultinomialConcentrationMinMaxDiff}
 \pr\Big( \Nmin \ge (1-\epsilon) \alphamin N\Big) &\wge 1-\delta_1, \\
 \pr\Big( \Delta N \le (2\epsilon \alphamax + \Delta\alpha) N\Big) &\wge 1-\delta_2,
\end{aligned}
\end{equation}
where $\delta_1 = K e^{-(\epsilon^2/2)\alphamin N}$ and $\delta_2 = 2 K e^{-(\epsilon^2/3)\alphamin N}$, together with $\alphamin = \min_k \alpha_k$, $\alphamax = \max_k \alpha_k$, and $\Delta\alpha = \max_{k,\ell} \abs{\alpha_k-\alpha_\ell}$.

(ii) If $\pr$ is the uniform distribution on $[K]^N$, then
\begin{equation*}
\begin{aligned}
 \pr\Big( \Nmin \ge (1-\epsilon) K^{-1} N \Big) &\wge 1 - \delta_1, \\
 \pr\Big( \Delta N \le 2 \epsilon K^{-1} N \Big) &\wge 1 - \delta_2,
\end{aligned}
\end{equation*}
with $\delta_1 = e^{\log K -\epsilon^2 N/(2K)}$ and
$\delta_2 = 2 e^{\log K -\epsilon^2 N/(3K)}$.
\end{lemma}
\begin{proof}
(i) By Lemma~\ref{the:MultinomialConcentration}, then events $\cA_{\epsilon}$ and $\cA_{\epsilon,+}$ defined by \eqref{eq:MultinomialConcentration}--\eqref{eq:MultinomialConcentrationPlus} satisfy $\pr( \cA_{\epsilon_+}^c ) \le \delta_1$ and $\pr( \cA_{\epsilon}^c ) \le \delta_2$. On the event $\cA_{\epsilon,+}$,
$\Nmin \ge (1-\epsilon)\alphamin N$.
Hence the first inequality in \eqref{eq:MultinomialConcentrationMinMaxDiff} follows. For the second inequality, we note that on the event $\cA_\epsilon$
\begin{align*}
 \abs{N_k - N_\ell}
 &\wle \abs{N_k - \alpha_k N} + \abs{N_\ell-\alpha_\ell N} + \abs{\alpha_k N - \alpha_\ell N} \\
 &\wle \epsilon \alpha_k N + \epsilon \alpha_\ell N + \abs{\alpha_k - \alpha_\ell} N \\
 &\wle 2 \epsilon \alphamax N + \Delta\alpha N
\end{align*}
for all $k,\ell$. This confirms the second inequality in \eqref{eq:MultinomialConcentrationMinMaxDiff}.

(ii) This follows immediately from (i) after noting that the uniform distribution on $[K]^N$ can be represented as $\pi = \alpha^{\otimes N}$ where $\alpha_k = K^{-1}$ for all $k$.
\end{proof}

\subsection{Elementary analysis}

\begin{lemma}
\label{the:GeometricMoments}
For any integer $j \ge 1$ and any real number $0 \le q < 1$,
\[
 \sum\limits_{k=j}^\infty \binom{k}{j} q^{k-j} = (1-q)^{-(j+1)}.
\]
\end{lemma}
\begin{proof}
Denote the falling factorial by $(x)_j = x(x-1)\cdots(x-j+1)$,  and let $f(q) = (1-q)^{-1}$. Then the $j$-th derivative of $f$ equals $f^{(j)}(q) = j! (1-q)^{-(j+1)}$. Because $f(q) = \sum_{k=0}^\infty q^k$, we find that the $j$-th derivative of $f$ also equals $\sum_{k=j}^\infty (k)_j q^{k-j}$. Hence the claim follows.
\end{proof}

\begin{lemma}
\label{the:Log}
(i) For $t \ge 0$, $\log(1+t) = t-\epsilon_1$ where $0 \le \epsilon_1 \le \frac12 t^2$.
(ii) For $0 \le t < 1$, $\log(1-t) = -t - \epsilon_2$ where $0 \le \epsilon_2 \le \frac{t^2}{2(1-t)^2}$, and especially, $0 \le \epsilon_2 \le 2 t^2$ for $0 \le t \le \frac12$.
\end{lemma}
\begin{proof}
(i) By taking two derivatives of $t \mapsto \log(1+t)$, we find that $\log(1+t) = t-\epsilon_1$ with
$\epsilon_1 = \int_0^t \int_0^s (1+u)^{-2} du ds$.

(ii) Similarly, we find that $\log(1-t) = -t-\epsilon_2$ with
$\epsilon_2 = \int_0^t \int_0^s (1-u)^{-2} du ds$.
\end{proof}

\begin{lemma}
\label{the:Power}
For any $0 \le x \le \frac12$ and $a>0$, the error term in the approximation $(1-x)^a = 1 - a x - r(x)$ is bounded by
$\abs{r(x)} \le \frac{2 \abs{a-1}}{2^a} a x^2$. Moreover, $r(x) \ge 0$ when $a \ge 1$.
\end{lemma}
\begin{proof}
The error term in the approximation $f(x) = f(0) + f'(0) x + r(x)$ equals $r(x) = \int_0^x \int_0^t f''(s) ds dt$ and is bounded by $\abs{r(x)} \le \frac12 c x^2$ with $c = \max_{0 \le x \le 1/2} \abs{f''(x)}$. The function $f(x) = (1-x)^a$ satisfies $f(0) = 1$ and $f'(0) = -a$, together with $f''(x) = a(a-1) (1-x)^{a-2}$. The claims follow after noticing that
\[
 \max_{0 \le x \le 1/2} \abs{f''(x)}
 \weq
 \begin{cases}
 \abs{f''(\tfrac12)}
 = \frac{4}{2^a} a \abs{a-1}
 &\quad \text{for $0 < a < 2$}, \\
 f''(0)
 = a (a-1)
 &\quad \text{for $a \ge 2$}.
 \end{cases}
\]
\end{proof}

\begin{lemma}
\label{the:SquareRoot}
Fix $0 \le \delta < 1$. Then the error term in the approximation $\sqrt{1-x} = 1-\frac12 x - \epsilon(x)$ satisfies
$0 \le \epsilon(x) \le c x^2$ for all $0 \le x \le \delta$, where $c = \frac18 (1-\delta)^{-3/2}$.
\end{lemma}
\begin{proof}
Consider Taylor's approximation $f(x) = f(0) + f'(0) x + r(x)$ where $r(x) = \int_0^x \int_0^t f''(s) ds dt$ is bounded by $\frac12 c_1 x^2 \le r(x) \le \frac12 c_2 x^2$ with $c_1 = \min_{0 \le x \le \delta} f''(x)$ and $c_2 = \max_{0 \le x \le \delta} f''(x)$. The function $f(x) = (1-x)^{1/2}$ satisfies $f(0) = 1$ and $f'(0) = -\frac12$, together with $f''(x) = -\frac14 (1-x)^{-3/2}$. Now $c_1 = -\frac14 (1-x)^{-3/2}$ and $c_2 = -\frac14 \le 0$. Hence the claim is true with $\epsilon(x) = -r(x)$.
\end{proof}

\begin{lemma}
\label{the:InclusionExclusion}
For any $0 \le p_1,\dots, p_n \le 1$, $A - B \le 1 - \prod_i (1-p_i) \le A$ with $A = \sum_i p_i$ and $B = \frac12 \sum_i \sum_{j \ne i} p_i p_j$.
\end{lemma}
\begin{proof}
Let $E_1,\dots,E_n$ be independent events with probabilities $p_1,\dots, p_n$. Apply inclusion--exclusion to the probability of the event $E = \cup_i E_i$ having probability $\pr(\cup_i E_i) = 1 - \pr(\cap_i E_i^c) = 1 - \prod_i (1-p_i)$.
\end{proof}

\begin{lemma}
\label{the:PowerSeriesBound}
For any integer $M \ge 1$ and any number $0 \le s < 1$,
\[
 M s^{M}
 \wle \sum_{m=M}^\infty m s^{m}
 \wle (1-s)^{-2} M s^M.
\]
\end{lemma}
\begin{proof}
Denote $S = \sum_{m=M}^\infty m s^{m}$. By differentiating 
$
 \sum_{m=M}^\infty s^{m}
 = (1-s)^{-1} s^M
$, we find that
\begin{align*}
 s^{-1} S
 \weq \sum_{m=M}^\infty m s^{m-1}
 \weq (1-s)^{-2} s^M + (1-s)^{-1} M s^{M-1},
\end{align*}
from which we see that
\begin{align*}
 S
 \weq s (1-s)^{-2} \Big( s^M + (1-s) M s^{M-1} \Big)
 \weq \frac{M s^M}{(1-s)^2} \Big( 1 - s(1-1/M)  \Big)
\end{align*}
The upper bound now follows from $1 - s(1-1/M) \le 1$.
The lower bound is immediate, corresponding to the first term of the nonnegative series.
\end{proof}

\subsection{Hamming distances}

\begin{lemma} 
\label{the:HamBound}
For any node labelling $\sigma: [N] \to [K]$, the number $Z_{\sigma,m}$ of node labellings  
$\sigma': [N] \to [K]$ such that $\Ham(\sigma, \sigma') = m$ satisfies
\[
 Z_{\sigma,m}
 \weq \binom{N}{m} (K-1)^m
 \wle \left( \frac{eN(K-1)}{m} \right)^m.
\]
\end{lemma}
\begin{proof}
Any node labelling $\sigma': [N] \to [K]$ which differs from a particular $\sigma$ at exactly $m$ input values can be constructed as follows. First choose a set of $m$ input values out of $N$; there are $\binom{N}{m}$ ways to do this.  Then for each $i$ of the chosen $m$ input values, select a new output value from the of $K-1$ values excluding $\sigma(i)$; there are $(K-1)^m$ ways to do this.  Hence the equality follows.

To verify the inequality, we note that $\frac{m^m}{m!} \le \sum_{s=0}^\infty \frac{m^s}{s!} = e^m$.  Therefore, we see that $\binom{N}{m} \le \frac{N^m}{m!} \le (\frac{eN}{m})^m$, and the inequality follows.
\end{proof}

\section{Comparing partitions}

\subsection{Classification error}

The absolute classification error between node labellings $\sigma, \sigma' : [N] \to [K]$ is defined by 
\[
 \ace(\sigma, \sigma')
 \weq \min_{\rho \in \Sym(K)} \Ham(\sigma, \rho \circ \sigma'),
\]
where $\Ham(\sigma, \sigma') = \sum_{i=1}^N 1(\sigma(i) \ne \sigma'(i))$ denotes the Hamming distance and $\Sym(K)$ denotes the group of permutations on $[K]$. We note that $\ace(\sigma, \sigma') = \ace(\rho \circ \sigma, \rho' \circ \sigma')$ for all $\rho, \rho' \in \Sym(K)$, which confirms that the classification error depends on its inputs only via the partitions induced by the preimages of the node labellings. The relative error $N^{-1} \ace(\sigma, \sigma')$ is usually called the classification error \cite{Meila_Heckerman_2001,Meila_2007}.

\subsection{Mirkin distance}
\label{sec:Mirkin}

The Mirkin distance is one of the common pair-counting based cluster validity indices \cite{Gosgens_Tikhonov_Prokhorenkova_2021, Lei_etal_2017}. It is defined between two nodes labellings $\sigma, \sigma': [N] \to [K]$ by
\[
 \Mir(\sigma, \sigma')
 \weq 2 \sum_{1 \le i < j \le N} \big( e_{ij}(1-e'_{ij}) + (1-e_{ij})e'_{ij} \big)
\]
where $e_{ij} = 1(\sigma(i) = \sigma(j))$ and $e'_{ij} = 1(\sigma'(i) = \sigma'(j))$.
The Mirkin distance is related to the Rand index by
$\Mir(\sigma, \sigma') = N(N-1) (1-\Rand(\sigma, \sigma'))$.

For any node labelling $\sigma: [N] \to [K]$, we denote by $E(\sigma)$ the set of unordered node pairs $\{i,j\}$ such that $\sigma(i)=\sigma(j)$, by $\Nmin^\sigma = \min_k \abs{C_k}$ and $\Nmax^\sigma = \max_k \abs{C_k}$ where $C_k = \{i: \sigma(i) = k\}$. Then we note that the Mirkin distance can be written as
\begin{equation}
 \label{eq:MirkinSets}
 \Mir(\sigma, \sigma')
 \weq 2 \left( \abs{ E(\sigma) \setminus E(\sigma')} + \abs{ E(\sigma') \setminus E(\sigma)} \right).
\end{equation}
The following result shows that when the Mirkin metric is small, then the maximum set sizes in two partitions cannot differ arbitrarily much.

\begin{lemma}
\label{the:PartitionPairBoundGeneral}
For any node labellings $\sigma, \sigma': [N] \to [K]$,
\[
 \abs{E(\sigma) \setminus E(\sigma')}
 \wge \frac12 ( \Nmax^\sigma -  \Nmax^{\sigma'} ) \Nmax^\sigma.
\]
\end{lemma}
\begin{proof}
For any $k$, denote by $E(C_k)$ the set of unordered node pairs in $C_k$.  Also denote $N_k = |E(C_k)|$, $N'_\ell = |E(C'_\ell)|$, and $N_{k\ell} = \abs{C_k \cap C'_\ell}$. Then we find that
\[
 \abs{E(C_k) \setminus E(\sigma')}
 \weq \binom{N_k}{2} - \sum_\ell  \binom{N_{k\ell}}{2}.
\]
By applying the bound $N_{k\ell} \le \Nmax^{\sigma'}$, we see that
\[
 \sum_{\ell} \binom{N_{k\ell}}{2}
 \wle \frac12 (\Nmax^{\sigma'} - 1) \sum_{\ell} N_{k\ell}
 \weq \frac12 (\Nmax^{\sigma'} - 1) N_{k}.
\]
Therefore,
\[
 \abs{E(C_k) \setminus E(\sigma')}
 \wge \binom{N_k}{2} - \frac12 ( \Nmax^{\sigma'} - 1) N_k
 \wge \frac12 ( N_k - \Nmax^{\sigma'} ) N_k.
\]
The claim now follows after noting that
\begin{align*}
 \abs{E(\sigma ) \setminus E(\sigma')}
 \weq \sum_k \abs{E(C_k) \setminus E(\sigma')} 
 \wge \max_{k} \abs{E(C_k) \setminus E(\sigma')}.
\end{align*}
\end{proof}

\subsection{Optimal alignments}

The confusion matrix of node labellings $\sigma, \sigma': [N] \to [K]$ is
the $K$-by-$K$ matrix having entries
\[
 N_{k\ell}
 \weq \abs{C_k \cap C'_\ell},
\] 
where $C_k = \sigma^{-1}(k)$ and $C'_\ell = (\sigma')^{-1}(\ell)$.  We say that node labellings $\sigma, \sigma': [N] \to [K]$ are optimally aligned if
\begin{equation}
 \label{eq:OptimalAlignment}
 \ace(\sigma, \sigma') \weq \Ham(\sigma, \sigma').
\end{equation} 
The following result provides an entrywise upper bound for the confusion matrix of optimally aligned node labellings.

\begin{lemma}
\label{the:PartitionIntersection}
If $\sigma$ and $\sigma'$ are optimally aligned, then the associated confusion matrix is bounded by
\begin{equation}
 \label{eq:Claim0}
 N_{k\ell} + N_{\ell k} 
 \wle N_{kk} + N_{\ell\ell}
\end{equation}
and
\begin{equation}
 \label{eq:Claim1New}
 N_{k\ell} \wle \frac{1}{3} ( N_k + N'_\ell )
\end{equation}
for all $k \ne \ell$, where
$N_k = \sum_\ell N_{k\ell}$ and $N'_\ell = \sum_k N_{k\ell}$.
\end{lemma}

\begin{proof}
Fix some distinct $k,\ell \in [K]$.
Define $\sigma'' = \tau \circ \sigma'$ where $\tau$ is the $K$-permutation which swaps $k$ and $\ell$ and leaves other elements of $[K]$ intact.  Denote $C''_j = (\sigma'')^{-1}(j)$. Then we see that $C''_j = C'_\ell$ for $j=k$, $C''_j = C'_k$ for $j=\ell$, and $C''_j = C'_j$ otherwise. Using the formulas
\[
 \Ham(\sigma, \sigma')
 = \sum_j \abs{C_j \setminus C'_j}
 \qquad\text{and}\qquad
 \Ham(\sigma, \sigma'')
 = \sum_j \abs{C_j \setminus C''_j}
\]
we find that
\[
 \Ham(\sigma, \sigma'') - \Ham(\sigma, \sigma')
 \weq \abs{C_k \setminus C'_\ell} 
 - \abs{C_k \setminus C'_k}
 + \abs{C_\ell \setminus C'_k} - \abs{C_\ell \setminus C'_\ell}.
\]
Because
\begin{align*}
 \abs{C_k \setminus C'_\ell} - \abs{C_k \setminus C'_k}
 &\weq (N_k - N_{k\ell}) - (N_k - N_{kk})
 \weq N_{kk} - N_{k\ell},
\end{align*}
and the same formula holds also with the roles of $k$ and $\ell$ swapped, it follows that
\[
 \Ham(\sigma, \sigma'') - \Ham(\sigma, \sigma')
 \weq N_{kk} - N_{k\ell} + N_{\ell\ell} - N_{\ell k}.
\]
Because $\sigma$ and $\sigma'$ are optimally aligned, we see that $\Ham(\sigma, \sigma') \le \Ham(\sigma, \sigma'')$.  Therefore, the left side of the above equality is nonnegative, and \eqref{eq:Claim0} follows.

Next, by applying the bounds
$N_{kk} \le N_k - N_{k\ell}$ and $N_{\ell\ell} \le N'_\ell - N_{k\ell}$, we may conclude that
\[
 0
 \wle N_{kk} - N_{k\ell} + N_{\ell\ell} - N_{\ell k}
 \wle N_k + N'_\ell - 3 N_{k\ell} - N_{\ell k}.
\]
The inequality \eqref{eq:Claim1New} now follows by noting that
\[
 N_{k\ell}
 \wle \frac13 ( N_k + N'_\ell - N_{\ell k} )
 \wle \frac13 ( N_k + N'_\ell ).
\]
\end{proof}

\subsection{Relating the classification error and the Mirkin distance}

The next result provides a way to bound the absolute classification error $\ace(\sigma, \sigma')$ using the Mirkin distance $\Mir(\sigma, \sigma')$. 

\begin{lemma}
\label{the:PartitionPairBound}
For any node labellings $\sigma, \sigma': [N] \to [K]$, 
\[
 \abs{E(\sigma) \setminus E(\sigma')}
 \wge \max\left\{ \Nmin^\sigma - \ace(\sigma, \sigma'), \ \frac13 \Nmin^\sigma - \frac16 \Nmax^{\sigma'} \right\} \ace(\sigma, \sigma').
\]
\end{lemma}



\begin{proof}
Let us note that all quantities appearing in the statement of the lemma remain invariant if we replace $\sigma'$ by $\rho \circ \sigma'$, where $\rho \in \Sym(K)$ is an arbitrary permutation. Therefore, we may without loss of generality assume that $\sigma$ and $\sigma'$ are optimally aligned according to \eqref{eq:OptimalAlignment}.

For sets $C,D \subset [N]$, we denote by $E(C,D)$ the collection of unordered pairs which can be written as $e = \{i,j\}$ with $i \in C$ and $j \in D$, and we denote the set of node pairs internal to $C$ by $E(C) = E(C,C)$. 
%
%
We observe that the set $\Gamma = E(\sigma) \setminus E(\sigma')$ can be partitioned into $\Gamma = \cup_k \Gamma_k$, where $\Gamma_k = E(C_k) \setminus E(\sigma')$. We may further split this set according to $\Gamma_k = \Gamma_{k1} \cup \Gamma_{k2}$, where
\begin{align*}
 \Gamma_{k1}
 &\weq E(C_k \cap C'_k, \, C_k \setminus C'_k), \\
 \Gamma_{k2}
 &\weq E(C_k \setminus C'_k) \setminus E(\sigma').
\end{align*}
Therefore, it follows that
$
 \abs{\Gamma} = \sum_k ( \abs{\Gamma_{k1}} + \abs{\Gamma_{k2}} ).
$

To analyse the sizes of $\Gamma_{k1}$ and $\Gamma_{k2}$, denote $N_{k\ell} = \abs{C_k \cap C'_\ell}$ and $D_k = \abs{C_k \setminus C'_k}$. Then we immediately see that 
\begin{equation}
 \label{eq:GammaK1}
 \abs{\Gamma_{k1}}
 \weq N_{kk} D_k.
\end{equation}
Furthermore, we see that
$
 E(C_k \setminus C'_k) \cap E(\sigma')
 = \cup_{\ell \ne k} E(C_k \cap C'_\ell),
$
and it follows that
\begin{equation}
 \label{eq:GammaK2}
 \abs{\Gamma_{k2}}
 \weq \abs{E(C_k \setminus C'_k)} - \sum_{\ell \ne k} \abs{E(C_k \cap C'_\ell)} 
 \weq \binom{D_k}{2} - \sum_{\ell \ne k} \binom{N_{k\ell}}{2}.
\end{equation}
By combining \eqref{eq:GammaK1}--\eqref{eq:GammaK2} we conclude that
\[
 \abs{\Gamma}
 \weq \sum_k ( \abs{\Gamma_{k1}} + \abs{\Gamma_{k2}} )
 \weq \sum_k \bigg\{ N_{kk} D_k + \binom{D_k}{2}  - \sum_{\ell \ne k} \binom{N_{k\ell}}{2} \bigg\}.
\]

Let us derive a lower bound for $\abs{\Gamma}$. Denote $B_k = \max_{\ell \ne k} N_{k\ell}$.  Then by noting that $\sum_{\ell \ne k} N_{k\ell} = D_k$, we see that
\[
 \sum_{\ell \ne k} \binom{N_{k\ell}}{2}
 \weq \frac12 \sum_{\ell \ne k} N_{k\ell}( N_{k\ell}-1)
 \wle \frac12 D_k (B_k-1),
\]
and by applying \eqref{eq:GammaK2}, it follows that
\begin{align*}
 \abs{\Gamma_{k2}}
 &\wge \frac12 D_k (D_k-1) - \frac12 D_k (B_k -1)
 \weq \frac12 D_k(D_k-B_k).
\end{align*}
By applying \eqref{eq:GammaK1} and noting that $N_{kk} = N_k - D_k$, it now follows that
\begin{equation}
 \label{eq:GammaKBoundGen}
 \abs{\Gamma_{k}} \wge D_k(N_k - D_k) + \frac12 D_k(D_k-B_k).
\end{equation}

We shall apply \eqref{eq:GammaKBoundGen} to derive two lower bounds for $\abs{\Gamma}$.
First, by Lemma~\ref{the:PartitionIntersection}, we find that
$
 B_k \le \frac{1}{3} ( N_k + \Nmax^{\sigma'} ),
$
and hence
\begin{align*}
 \abs{\Gamma_{k}}
 &\wge (N_k - D_k)D_k + \frac12 \left( D_k - \frac{1}{3} N_k - \frac{1}{3} \Nmax^{\sigma'} \right) D_k \\
 &\weq \bigg( \frac{5}{6} N_k - \frac12 D_k - \frac16 \Nmax^{\sigma'} \bigg) D_k.
\end{align*}
Because $D_k \le N_k$, we conclude that
\begin{align*}
 \abs{\Gamma_{k}}
 &\wge \bigg( \frac{1}{3} N_k - \frac16 \Nmax^{\sigma'} \bigg) D_k
 \wge \bigg( \frac{1}{3} \Nmin^{\sigma} - \frac16 \Nmax^{\sigma'} \bigg) D_k
\end{align*}
By summing the above inequality over $k$ and noting that $\sum_k D_k = \Ham(\sigma, \sigma') = L$ for optimally aligned $\sigma$ and $\sigma'$, we conclude that
\begin{equation}
 \label{eq:GammaKBound1}
 \abs{\Gamma}
 \wge \bigg( \frac{1}{3} \Nmin^{\sigma} - \frac16 \Nmax^{\sigma'} \bigg) L.
\end{equation}
 
Second, by noting that $B_k \le D_k$, we see that \eqref{eq:GammaKBoundGen} implies
\[
 \abs{\Gamma_{k}}
 \wge D_k(N_k - D_k) 
 \wge D_k(\Nmin^{\sigma} - D_k).
\]
By summing the above inequality over $k$, we find that
\[
 \abs{\Gamma}
 \wge \Nmin^{\sigma} \sum_k D_k - \sum_k D_k^2
 \wge \Nmin^{\sigma} \sum_k D_k - (\sum_k D_k)^2.
\]
By recalling that $\sum_k D_k = L$, we conclude that
\begin{equation}
 \label{eq:GammaKBound2}
 \abs{\Gamma}
 \wge \Nmin^{\sigma} L - L^2.
\end{equation}
By combining \eqref{eq:GammaKBound1}--\eqref{eq:GammaKBound2}, the claim follows.
\end{proof}

\section{Proof of the lower bound of Theorem~\ref{thm:asymptotic_exponential_bound_recovery}}
\label{appendix:lower_bound_proof}

This section is devoted to proving the lower bound of Theorem~\ref{thm:asymptotic_exponential_bound_recovery} and is organised as follows:
Section~\ref{sec:LowerBoundInhomogeneous} describes a lower bound (Theorem~\ref{the:LowerBoundGen}) which is valid for general SBMs, not necessarily homogeneous or binary.
Section~\ref{sec:ProofOfLowerBoundGen} presents the proof of Theorem~\ref{the:LowerBoundGen}.
Section~\ref{subsection:lower_bound_application_homogeneous_model} specialises the lower bound into homogeneous SBMs and leads to Proposition~\ref{the:LowerBoundHom}.

\subsection{A quantitative lower bound}
\label{sec:LowerBoundInhomogeneous}

The following theorem lower bounds the expected loss made by any algorithm in clustering a non-homogeneous SBM.

\newcommand{\alphamaxK}{\alpha_{{\rm max}, \cK}}
\begin{theorem}
\label{the:LowerBoundGen}
Consider a SBM defined by \eqref{eq:PairwiseInteractionModel}--\eqref{eq:JointProbability} where the block membership structure is distributed according to $\pi = \alpha^{\otimes N}$ for some probability distribution $\alpha$ on $[K]$. Fix an arbitrary $\cK \subset [K]$ and probability distributions $f^*_1,\dots, f^*_K$.
Assume that $N \ge 8 \alphamin^{-1} \log (K/\delta)$ for
$\delta = \frac14 \left( \alpha_\cK - \alphamaxK \right)$.  Then for any estimator $\hsigma: \cX \to \cZ$, the error 
is lower bounded in expectation by
\begin{equation}
 \label{eq:LowerBoundGen}
 \E \dhams\left( \hsigma \right) 
 \wge \frac{1}{21} N \alphamin^2 \delta e^{-N I_1 - \alpha_\cK^{1/2}
  \delta^{-1/2} \sqrt{N I_{21} + N^2 I_{22}}}
 - \frac16 N \alphamin K e^{- \frac18 N \alphamin},
\end{equation}
where the quantities $I_1, I_{21}$ and $I_{22}$ are defined by
\begin{equation}
\label{eq:I1I21I22}
\begin{aligned}
I_1 &\weq \sum_k \sum_{\ell} \alpha^*_k \alpha_\ell \dkl(f^*_\ell \| f_{k\ell}), \\
I_{21} &\weq \sum_k \sum_{\ell} \alpha^*_k \alpha_\ell \vkl(f^*_\ell \| f_{k\ell})
+ \sum_k \alpha^*_k B_k, \\
I_{22} &\weq \sum_k \alpha^*_k A_k^2 - \Big(\sum_k \alpha^*_k A_k \Big)^2,
\end{aligned}
\end{equation}
with $A_k = \sum_\ell \alpha_\ell \dkl(f^*_\ell \| f_{k\ell})$ and $B_k = \sum_\ell \alpha_\ell \dkl(f^*_\ell \| f_{k\ell})^2 - ( \sum_\ell \alpha_\ell \dkl(f^*_\ell \| f_{k\ell}))^2$, together with
$\alpha^*_k = 1(k \in \cK) \frac{\alpha_k}{\alpha_\cK}$ and $\alpha_\cK = \sum_{k \in \cK} \alpha_k$.
\end{theorem}

\begin{remark}
\label{rem:LowerBoundGen}
The second term on the right side of \eqref{eq:LowerBoundGen} is $o(1)$ when $\alphamin \ge 9 N^{-1} \log N$ and $2 \le K \le N$.
\end{remark}

\begin{remark}
\label{rem:LowerBoundGenNovelty}
The lower bound of Theorem~\ref{the:LowerBoundGen} is quantitative, and hence valid regardless of any scaling assumptions, and also for all finite models with fixed, not asymptotic, size. This is one of the first explicit quantitative lower bounds in this context.
\end{remark}

\begin{remark}
In homogeneous models with uniform node labels, one can specify the quantities $I_1, I_{21}$ and $I_{12}$ to obtain the lower bound stated in Theorem~\ref{thm:asymptotic_exponential_bound_recovery}. This is done in Section~\ref{subsection:lower_bound_application_homogeneous_model}.  
\end{remark}

\subsection{Proof of Theorem~\ref{the:LowerBoundGen}}
\label{sec:ProofOfLowerBoundGen}

This section is devoted to proving Theorem \ref{the:LowerBoundGen} step by step.

\subsubsection{Key result on block permutations}
The following key result implies that when $L(\sigma_1, \sigma_2) = \min_\tau \dham(\sigma_1, \tau \circ \sigma_2) < \frac12 \Nmin(\sigma_1)$, then the minimum Hamming distance is attained by a unique block permutation.
\begin{lemma}
\label{the:UniquePermutation}
Let $\sigma_1, \sigma_2: [N] \to [K]$ be such that $\dham(\sigma_1, \tau^* \circ \sigma_2) < \frac12 \Nmin$ for some $K$-permutation $\tau^*$, where $\Nmin = \min_k \abs{\sigma_1^{-1}(k)}$. Then $\tau^*$ is the unique minimiser of $\tau \mapsto \dham(\sigma_1, \tau \circ \sigma_2)$.
\end{lemma}
 This corresponds to \cite[Lemma B.6]{Xu_Jog_Loh_2020}.

 \begin{proof}
 Assume that $\tau \in \Sym(K)$ satisfies $\dham(\tau \circ \sigma_1, \sigma_2) < \frac{s}{2}$, where $s = \Nmin$. Fix $k \in [K]$ and let $U_k = \{i: \sigma_1(i) = k, \sigma_2(i) \ne \tau(k)\}$. Then every node $i$ in $U_k$ satisfies $\tau\circ\sigma_1(i) \ne \sigma_2(i)$, and therefore
 $ \abs{U_k} \le \dham(\tau \circ \sigma_1, \sigma_2) < \frac{s}{2}.$
 Hence for any $\ell \ne \tau(k)$,
    \[
     \abs{\sigma_1^{-1}(k) \cap \sigma_2^{-1}(\ell)}
     \wle \abs{U_k}
     \ < \ \frac{s}{2}.
    \] 
 On the other hand,
    \[
     \abs{\sigma_1^{-1}(k) \cap \sigma_2^{-1}(\tau(k)) }
     \weq \abs{\sigma_1^{-1}(k)} - \abs{U_k}
     \wge s - \frac{s}{2}
     \wge \frac{s}{2}.
    \]
 Hence $\tau(k)$ is the unique value which maximizes $\ell \mapsto \abs{\sigma_1^{-1}(k) \cap \sigma_2^{-1}(\ell)}$. Because this conclusion holds for all $k$, it follows that $\tau$ is uniquely defined.
 \end{proof}

\subsubsection{Lower bounding by critical node count}
\label{sec:LowerBoundingByCriticalNodeCountGen}

This method apparently originates from \cite{Zhang_Zhou_2016}. Let $\Opt(\sigma_1,\sigma_2)$ be the set of $K$-permutations $\tau$ for which $\dham(\sigma_1, \tau\circ\sigma_2)$ is minimised. Given an estimated node labelling $\hsigma_x$, we define a set of critical nodes by
\[
 \Crit(\sigma,\hsigma_x)
 \weq \{ j \in [N]: \sigma(j) \ne \tau\circ\hsigma_x(j) \ \text{for some $\tau \in \Opt(\sigma,\hsigma_x)$}\}.
\]
We denote the number of critical nodes by
\[
 \nCrit(\sigma,x)
 \weq \abs{\Crit(\sigma,\hsigma_x)}.
\]

\begin{lemma}
\label{the:MeanPermutedHammingLowerBoundGen}
For any estimate $\hsigma_x$ obtained as a deterministic function of observed data, let $L = L(\sigma, \hsigma_x) = \Ham^*(\sigma, \hsigma_x)$. Then
\begin{equation}
 \label{eq:MeanPermutedHammingLowerBoundGen}
 \E L
 \wge \frac{\alphamin}{6} \left( \E \nCrit - N K e^{- \frac18 N \alphamin} \right).
\end{equation}
\end{lemma}

\begin{proof}
We shall consider $\nCrit = \nCrit(\sigma, \hsigma_x)$, and $\Nmin = \Nmin(\sigma)$ as random variables defined on $\cZ \times \cX$. By Lemma~\ref{the:UniquePermutation}, $L = \nCrit$ on the event $L < c$ where $c = \frac12 \Nmin$. Given a node labelling $\sigma \in \cZ$, we consider $x \mapsto \hsigma_x$, $x \mapsto L(\sigma,x)$ and $x \mapsto \nCrit(\sigma,x)$ as random variables on $\cX$.
Consider the following two cases:
\begin{enumerate}[(i)]
\item 
If $P_\sigma(L \ge c) \ge \frac{1}{N+c} E_\sigma \nCrit$, then 
\begin{align*}
 E_\sigma L 1( L \ge c )
 \wge c P_\sigma( L \ge c )
 \wge \frac{c}{N+c} E_\sigma \nCrit.
\end{align*}

\item If $P_\sigma(L \ge c) \le \frac{1}{N+c} E_\sigma \nCrit$, then 
\begin{align*}
 E_\sigma \nCrit 1( L \ge c )
 \wle N P_\sigma( L \ge c )
 \wle \frac{N}{N+c} E_\sigma \nCrit,
\end{align*}
so that
\begin{align*}
 E_\sigma L 1( L < c )
 \weq E_\sigma \nCrit 1( L < c )
 \weq E_\sigma \nCrit - E_\sigma \nCrit 1( L \ge c )
 \wge \frac{c}{N+c} E_\sigma \nCrit.
\end{align*}
\end{enumerate}
In both cases, $E_\sigma L \ge \frac{c}{N+c} E_\sigma \nCrit$, so that
\[
 E_\sigma L
 \wge \frac{\Nmin}{2N + \Nmin} E_\sigma \nCrit
 \wge \frac{\Nmin}{3N} E_\sigma \nCrit.
\]
By taking expectations with respect to the prior, we find that
\begin{equation}
 \label{eq:MeanPermutedHammingLowerBoundPreGen}
 \E L
 \wge \frac{1}{3N} \E \Nmin Y,
\end{equation}
where $Y= E_\sigma \nCrit$ is viewed as a random variable on probability space $\cS$ equipped with probability measure $\pi$. Let $t = \frac12 N\alphamin$. We note that $0 \le Y \le N$ surely, and that $\Nmin > t$ with high probability. Observe that
\begin{align*}
 \E \Nmin Y
 \wge \E \Nmin Y \, 1(\Nmin>t)
 \wge t \, \E Y 1(\Nmin > t).
\end{align*}
and, due to $Y \le N$,
\begin{align*}
 \E Y 1(\Nmin > t)
 \weq \E Y - \E Y 1(\Nmin \le t)
 \wge \E Y - N \pr(\Nmin \le t).
\end{align*}
By noting that $\E Y = \E \nCrit$ and applying Lemma~\ref{the:MultinomialConcentrationMinMaxDiff}, we find that
\begin{align*}
 \E \Nmin Y
 &\wge t \left( \E Y - N \pr(\Nmin \le t) \right) \\
 &\weq \frac12 N\alphamin \left( \E \nCrit - N \pr(\Nmin \le \frac12 N\alphamin) \right) \\
 &\wge \frac12 N\alphamin \left( \E \nCrit - N K e^{- \frac18 N \alphamin} \right).
\end{align*}
Together with \eqref{eq:MeanPermutedHammingLowerBoundPreGen}, the claim now follows.
\end{proof}

\subsubsection{Change of measure}

Fix a reference node $i$, a set $\cK \subset [K]$, and some probability distributions $f^*_1,\dots, f^*_K$ on the interaction space $S$. We define an alternative statistical model for $\sigma$ and $x$ by modifying $P_{\sigma}(x)$ defined in~\eqref{eq:PairwiseInteractionModel} according to
\begin{equation}
 \label{eq:KernelModifiedGen}
 P^{*i}_{\sigma}(x)
 \weq \left( 1_{\cK}(\sigma(i)) \prod_{j \ne i} \frac{f^*_{\sigma(j)}(x_{ij})}{f_{\sigma(i) \sigma(j)}(x_{ij})}
 + 1_{\cK^c}(\sigma(i)) \right) P_{\sigma}(x),
\end{equation}
and defining a modified probability measure on $\cZ \times \cX$ by
\begin{equation}
 \label{eq:JointModifiedGen}
 \pr^{*i}(\sigma,x)
 \weq \pi(\sigma) P^{*i}_{\sigma}(x).
\end{equation}

In the modified model, node labels are sampled independently as before, and all interactions not involving node $i$ are sampled just as in the original model. If the label of node $i$ belongs to $\cK$, then we sample all $i$-interactions from $f^*_1, \dots, f^*_K$.  The following lemma confirms that under the alternative model, $\sigma_i$ is conditionally independent of observed data $x$ and other labels $\sigma_{-i}$ given $\sigma_i \in \cK$.

\begin{lemma}
\label{the:ConditionalDistributionAlt}
For $(\sigma, x)$ sampled from model \eqref{eq:JointModifiedGen}, the conditional distribution of the label $\sigma_i$ given that $\sigma_i \in \cK$, the other labels are $\sigma_{-i}$, and the observed interactions are $x$, equals
\[
 \pr^{*i}( \sigma_i = k \cond \sigma_i \in \cK, \sigma_{-i}, x )
 \weq \alpha^*_k
 \qquad \text{for all $\sigma_{-i}, x$},
\]
where $\alpha^*_k = 1(k \in \cK) \frac{\alpha_k}{\alpha_\cK}$.
\end{lemma}
\begin{proof}
Observe that $P^{*i}_{\sigma}(x) = Q_{\sigma_{-i}}(x)$ for all $\sigma$ such that $\sigma_i \in \cK$, where
\[
 Q_{\sigma_{-i}}(x)
 \weq 
 \bigg( \prod_{j \ne i} f^*_{\sigma_{-i}(j)}(x_{ij}) \bigg)
 \bigg( \prod_{uv \in E_{-i}} f_{\sigma_{-i}(u) \sigma_{-i}(v)}(x_{uv}) \bigg),
\]
and $E_{-i}$ is the set of unordered node pairs not incident to $i$. Especially, $\pr^{*i}( \sigma, x ) = \alpha(\sigma_i) \, \pi_{-i}(\sigma_{-i}) \, Q_{\sigma_{-i}}(x)$ whenever $\sigma_i \in \cK$.
Hence the conditional probability distribution of $\sigma_i$ given $(\sigma_{-i}, x)$ satisfies
$\pr^{*i}( \sigma_i \cond \sigma_{-i}, x ) = \alpha(\sigma_i)$ for all $\sigma_i \in \cK$. The claim follows by summing this equality with respect to $\sigma_i \in \cK$.
\end{proof}


To analyse how much the alternative model differs from the original model, we will investigate the associated log-likelihood ratio 
\[
 \LLR_i(\sigma,x)
 \weq \log \frac{\pr^{*i}(\sigma,x)}{\pr(\sigma,x)}.
\]

\begin{lemma}
\label{the:LLRMeanVarGen}
The mean and variance of the log-likelihood ratio given $\sigma(i) \in \cK$ are equal to
%
%
$\E^{*i}( \LLR_i \cond \sigma_i \in \cK) = (N-1) I_1$ and $\V^{*i}( \LLR_i \cond \sigma_i \in \cK) = (N-1) I_{21} + (N-1)^2 I_{22}$, where $I_1,I_{21},I_{22}$ are given by \eqref{eq:I1I21I22}.
\end{lemma}
\begin{proof}
The conditional distribution of $(\sigma,x)$ sampled from $\pr^{*i}$ given $\sigma(i) \in \cK$ can be represented as
\[
 \tilde\pr^{*i}(\sigma, x)
 \weq \tilde\pi^{*i}(\sigma) P_\sigma^{*i}(x),
\]
where $\tilde\pi^{*i}(\sigma) = \alpha^*_{\sigma(i)} \prod_{j \ne i} \alpha_{\sigma(j)}$
and $\alpha^*_k = 1(k \in \cK) \frac{\alpha_k}{\alpha_\cK}$, and 
$P_\sigma^{*i}$ is defined by \eqref{eq:KernelModifiedGen}.
Furthermore, the log-likelihood ratio can be written as
\[
 \LLR_i(\sigma, x)
 \weq 1(\sigma(i) \in \cK) \sum_{j \ne i} \log \frac{f^*_{\sigma(j)}(x_{ij})}{f_{\sigma(i) \sigma(j)}(x_{ij})}.
\]
The conditional expectation $A(\sigma) = E^{*i}_\sigma \LLR_i$ of the log-likelihood ratio given $\sigma$ hence equals
\[
 A(\sigma) \weq 1(\sigma(i) \in \cK) \sum_{j \ne i} m_{\sigma(i) \sigma(j)}
\]
where $m_{k\ell} = \dkl(f^*_\ell \| f_{k\ell})$. Hence, treating
$(\sigma,x) \mapsto \sigma(i)$, $(\sigma,x) \mapsto \sigma(j)$, and $(\sigma,x) \mapsto A(\sigma)$,
as random variables on probability space $(\cZ \times \cX, \, \tilde\pr^{*i})$, and noting that $\sigma(i) \in \cK$ with $\tilde\pr^{*i}$-probability one, we find that
\[
 \tilde\E^{*i} \Lambda_i
 \weq \tilde\E^{*i} A
 \weq \sum_{j \ne i} \tilde\E^{*i} m_{\sigma(i) \sigma(j)}
 \weq (N-1) \sum_k \sum_\ell m_{k\ell} \alpha^*_k \alpha_\ell,
\]
which implies the first claim.


To compute the variance, we observe that
\begin{equation}
 \label{eq:VarDecomp}
 \tilde\V^{*i} \Lambda_i
 \weq \tilde\E^{*i} B + \tilde\V^{*i} A,
\end{equation}
where $B = V^{*i}_\sigma \LLR_i$. We note that by the conditional independence of $x_{ij}$, $j \ne i$, given $\sigma$, it follows that
\[
 B
 \weq 1(\sigma(i) \in \cK) \sum_{j \ne i} v_{\sigma(i) \sigma(j)},
\]
where $v_{k\ell} = \vkl(f^*_\ell \| f_{k\ell})$. By taking expectations, we find that
\begin{equation}
 \label{eq:MeanB}
 \tilde\E^{*i} B
 \weq (N-1) \sum_{k} \sum_{\ell} \alpha^*_k \alpha_\ell v_{k\ell}.
\end{equation}
We still need to compute the variance of $A$. To do this, we condition on the label of node $i$ and observe that on the event $\sigma(i) \in \cK$ of $\tilde\pr^{*i}$-probability one, 
\begin{align*}
 \tilde\E^{*i} (A \cond \sigma(i) \, )
 &\weq (N-1) A_{\sigma(i)}, \\
 \tilde\V^{*i} (A \cond \sigma(i) \, )
 &\weq (N-1) B_{\sigma(i)},
\end{align*}
where $A_k = \sum_\ell \alpha_\ell m_{k\ell}$ and $B_k = \sum_\ell \alpha_\ell m_{k\ell}^2 - ( \sum_\ell \alpha_\ell m_{k\ell})^2$. Therefore,
\begin{align*}
 \tilde\V^{*i} A
 &\weq \tilde\E^{*i} \tilde\V^{*i} (A \cond \sigma(i) \, ) + \V \tilde\E^{*i} (A \cond \sigma(i) \, ) \\
 &\weq (N-1) \tilde\E^{*i} B_{\sigma(i)} + (N-1)^2 \tilde\V^{*i} A_{\sigma(i)} \\
 &\weq (N-1) \sum_{k} \alpha^*_k B_k
 + (N-1)^2 \bigg\{ \sum_k \alpha^*_k A_k^2 - \Big(\sum_k \alpha^*_k A_k \Big)^2 \bigg\}.
\end{align*}
By combining this with \eqref{eq:VarDecomp} and \eqref{eq:MeanB}, we find that
\begin{align*}
 \tilde\V^{*i} \Lambda_i
 &\weq (N-1) \sum_{k} \sum_{\ell}  \alpha^*_k \alpha_\ell v_{k\ell}
 + (N-1) \sum_{k} \alpha^*_k B_k \\
 &\qquad + (N-1)^2 \bigg\{ \sum_k \alpha^*_k A_k^2
 - \Big(\sum_k \alpha^*_k A_k \Big)^2 \bigg\},
\end{align*}
and the second claim follows.
\end{proof}

\subsubsection{Lower bound of critical node count} The following is key to proving the lower bound, and rigorously handling stochastic dependencies implied by optimal $K$-permutations in the definition of $L$.  
Recall that $\alpha_\cK = \sum_{k \in \cK} \alpha_k$ together with
$\alphamin = \min_{k \in [K]} \alpha_k$ and
$\alphamaxK = \max_{k \in \cK} \alpha_k$.

\begin{lemma}
\label{the:LowerBoundCriticalNodeCountGen}
Assume that
$N \ge 8 \alphamin^{-1} \log (K/\delta)$ for
$\delta = \frac14 \left( \alpha_\cK - \alphamaxK \right)$.
Then for any estimator $x \mapsto \hsigma_x$, the expected number of critical nodes is bounded by
\begin{equation}
 \label{eq:MidJulyBound4HighlightedGen}
 \E \nCrit
 \wge \frac{2}{7} \alphamin \delta Ne^{-t}
\end{equation}
for $t = \max_i \Big( \E^{*i} (\LLR_i | \sigma_i \in \cK ) + \alpha_\cK^{1/2} \delta^{-1/2} \sqrt{\V^{*i}(\LLR_i | \sigma_i \in \cK) } \Big)$.
\end{lemma}

\begin{proof}
Denote $\epsilon = \frac16 \alphamin$.  The proof contains four steps which are treated one by one in what follows. 

(i) Denote the event that node $i$ is critical by
\[
 \cC_i
 \weq \big\{ (\sigma,x): \sigma(i) \ne \tau(\hsigma_x(i))
 \text{\ for some $\tau \in \Opt(\sigma, \hsigma_x)$} \big\},
\]
and let
\[
 \cE_i
 \weq \cC_i \cup \big\{ (\sigma,x): L^+(\sigma,\hsigma_x) > \epsilon N \big\}.
\]
Recall that $\E L^+ = \sum_i \pr(\cC_i)$. Markov's inequality then implies that
\begin{align*}
 \sum_i \pr( \cE_i )
 &\wle \sum_i \left( \pr(\cC_i) + (\epsilon N)^{-1} \E \nCrit \right).
\end{align*}
By noting that the right side above equals $(1+\epsilon^{-1}) \E \nCrit$, we obtain a lower bound
\begin{equation}
 \label{eq:MidJulyBound1GenNew}
 \E \nCrit
 \wge \frac{\epsilon}{1+\epsilon} \sum_i \pr( \cE_i ).
\end{equation}

(ii) We will now focus on a particular node $i$, and derive a lower bound for the probability of event $\cE_i$ under the perturbed model $\pr^{*i}$ defined by \eqref{eq:KernelModifiedGen}. We start by deriving an upper bound for the probability of the event
\[
 \pr^{*i}( \cE_i^c, \, \Nmin > 3 \epsilon N, \, \sigma(i) \in \cK)
 \weq \pr^{*i}( \cC_i^c, \, \cB, \, \sigma(i) \in \cK),
\]
where
\[
 \cB \weq \{ (\sigma, x): L^+(\sigma, \hsigma_x) \le \epsilon N, \, \Nmin(\sigma) > 3 \epsilon N\}
\]
and $\Nmin(\sigma) = \min_k \abs{\sigma^{-1}(k)}$. On the event $\cB$, we see that $L^+(\sigma, \hsigma_x) < \frac13 \Nmin(\sigma)$, and Lemma~\ref{the:UniquePermutation} implies that $L^+(\sigma, \hsigma_x) = \min_\tau \Ham(\sigma, \tau\circ \hsigma_x)$ is attained by a unique $K$-permutation $\tau$.  This is why we may split the above probability into
\begin{equation}
 \label{eq:CondDistrSplit1}
 \pr^{*i}( \cC_i^c, \, \cB, \, \sigma(i) \in \cK)
 \weq \sum_\tau \pr^{*i}( \cC_i^c, \, \cB_\tau, \, \sigma(i) \in \cK)
\end{equation}
where
\[
 \cB_\tau
 \weq \left\{ (\sigma,x): \Ham(\sigma, \tau\circ\hsigma_x) \le \epsilon N, \, \Nmin(\sigma) > 3 \epsilon N \right\}.
\]
To analyse events associated with $\cB_\tau$, let us introduce some more notation. We define $\Hami(\sigma_1,\sigma_2) = \sum_{j \ne i} 1(\sigma_1(j) \ne \sigma_2(j))$ and denote $\Nmini(\sigma) = \min_k \abs{\sigma^{-1}(k) \setminus \{i\}}$, and consider an event
\begin{align*}
 \cB^{-i}_\tau 
 &\weq \left\{ (\sigma,x): \Hami(\sigma, \tau\circ\hsigma_x) \le \epsilon N, \, \Nmini(\sigma) > 3 \epsilon N - 1\right\}.
\end{align*}
Then we find that
\begin{align*}
 \cC_i^c \cap \cB_\tau
 &\weq \{\sigma(i) = \tau(\hsigma_x(i))\} \cap \cB_\tau \\
 &\wsubset \{\sigma(i) = \tau(\hsigma_x(i))\} \cap \cB^{-i}_\tau,
\end{align*}
so that, under the conditional distribution $\tilde \pr^{*i}(\cdot) = \pr^{*i}(\cdot \cond \sigma(i) \in \cK)$,
\begin{equation}
 \label{eq:CondDistrBound1}
 \tilde\pr^{*i}( \cC_i^c, \, \cB_\tau)
 \wle \tilde\pr^{*i}( \sigma(i) = \tau(\hsigma_x(i)), \, \cB^{-i}_\tau).
\end{equation}
We note that the event $\cB^{-i}_\tau$ is completely determined by $(\sigma_{-i}, x)$, 
%
and according to Lemma~\ref{the:ConditionalDistributionAlt}, we know that when $(\sigma,x)$ is sampled from $\pr^{*i}$, then $\sigma(i)$ is $\alpha^*$-distributed and conditionally independent of $(\sigma_{-i},x)$ given $\sigma(i) \in \cK$.
Therefore, under the conditional distribution $\tilde \pr^{*i}(\cdot) = \pr^{*i}(\cdot \cond \sigma(i) \in \cK)$, we find that
\begin{align*}
 \tilde \pr^{*i} \big(\sigma(i) = \tau(\hsigma_x(i)), \, \cB^{-i}_\tau \big)
 &\weq \sum_{k \in \cK} \tilde\pr^{*i} \big(\sigma(i) = k, \, \tau(\hsigma_x(i)) = k , \, \cB^{-i}_\tau \big) \\
 &\weq \sum_{k \in \cK} \alpha^*_k
 \, \tilde\pr^{*i}\big( \tau(\hsigma_x(i)) = k , \, \cB^{-i}_\tau \big),
\end{align*}
from which we conclude together with \eqref{eq:CondDistrBound1} that
\begin{align*}
 \tilde\pr^{*i}( \cC_i^c, \, \cB_\tau)
 &\wle  \frac{\alphamaxK}{\alpha_\cK} \, 
 \tilde\pr^{*i}\big( \cB^{-i}_\tau \big).
\end{align*}
Because $N \ge \epsilon^{-1}$ due to $\log (K/\delta) \ge \log(4K) \ge 1$ and $N \ge 8 \alphamin^{-1} \log (K/\delta)$, we see that $\epsilon N < \frac12 (3 \epsilon N - 1)$.  Therefore, $\Hami(\sigma, \tau\circ\hsigma_x) < \frac12 \Nmini(\sigma)$ on the event $\cB^{-i}_\tau$.
Then again by Lemma~\ref{the:UniquePermutation}, the events $\cB^{-i}_\tau$ are mutually exclusive,
%
and in light of \eqref{eq:CondDistrSplit1} it follows that 
\begin{align*}
 \tilde\pr^{*i}\big( \cC_i^c, \, \cB \big)
 \wle \frac{\alphamaxK}{\alpha_\cK} \, \tilde\pr^{*i}(\cup_\tau \cB^{-i}_\tau)
 \wle \frac{\alphamaxK}{\alpha_\cK}.
\end{align*}
By recalling the definitions of $\cC_i, \cE_i$, we now conclude that
\begin{align*}
 \pr^{*i}\big( \cE_i^c, \, \Nmin > 3 \epsilon N, \, \sigma(i) \in \cK \big)
 &\weq \pr^{*i}\big( \cC_i^c, L^+ \le \epsilon N, \, \Nmin > 3 \epsilon N, \, \sigma(i) \in \cK \big) \\
 &\weq \pr^{*i}\big( \cC_i^c, \, \cB, \, \sigma(i) \in \cK \big) \\
 &\wle \alphamaxK,
\end{align*}
and therefore,
\begin{equation}
 \label{eq:MidJulyBound2Gen}
 \pr^{*i}(\cE_i^c, \, \sigma(i) \in \cK)
 \wle \alphamaxK + \pr^{*i}(\Nmin \le 3 \epsilon N).
\end{equation}

(iii) Next, by recalling our choice of $\epsilon = \frac16 \alphamin$ and applying Lemma~\ref{the:MultinomialConcentrationMinMaxDiff}, we see that
$\pr^{*i}(\Nmin \le 3 \epsilon N) = \pr(\Nmin \le \frac12 N \alphamin) \le K e^{- \frac18 N \alphamin} \le \delta$ due to
$N \ge 8 \alphamin^{-1} \log (K/\delta)$.
By combining this with \eqref{eq:MidJulyBound2Gen}, we see that
$\pr^{*i}(\cE_i^c, \, \sigma(i) \in \cK) \le \alphamaxK + \delta$. Hence,
by our choice of $\delta$, it follows that
\begin{equation}
 \label{eq:MidJulyBound3GenNew}
 \begin{aligned}
 \pr^{*i}( \cE_i, \, \sigma(i) \in \cK)
 &\wge \pr( \sigma(i) \in \cK) - \alphamaxK - \delta \\
 &\weq \alpha_\cK - \alphamaxK - \delta \\
 &\weq 3 \delta.
 \end{aligned}
\end{equation}


(iv) Finally, we will transform the lower bound \eqref{eq:MidJulyBound3GenNew} into one involving the original probability distribution $\pr$ instead of $\pr^{*i}$. By writing
\begin{align*}
 \pr( \cE_i, \, \sigma(i) \in \cK)
 \weq \E^{*i} e^{-\Lambda_i} 1( \cE_i, \, \sigma(i) \in \cK),
\end{align*}
and noting that $e^{-\LLR_i} 1(\cE_i, \sigma_i \in \cK) \ge e^{-t} 1(\cE_i, \sigma_i \in \cK, \LLR_i \le t)$, it follows that 
\begin{align*}
 \pr( \cE_i, \, \sigma(i) \in \cK)
 &\wge e^{-t} \pr^{*i}( \cE_i, \, \sigma(i) \in \cK, \, \LLR_i \le t) \\
 &\wge e^{-t} \bigg( \pr^{*i}( \cE_i, \, \sigma(i) \in \cK) - \pr^{*i}( \cE_i, \, \sigma(i) \in \cK, \, \LLR_i > t) \bigg) \\
 &\wge e^{-t} \bigg( \pr^{*i}( \cE_i, \, \sigma(i) \in \cK) - \pr^{*i}( \sigma(i) \in \cK, \, \LLR_i > t) \bigg).
\end{align*}
For $t \ge \tilde\E^{*i} \left( \LLR_i \right) + \left( \frac{\alpha_\cK}{\delta} \tilde\V^{*i}(\LLR_i) \right)^{1/2}$, Chebyshev's inequality implies that $\tilde\pr^{*i}(\LLR_i > t) \le \frac{\delta}{\alpha_\cK}$, and hence $\pr^{*i}(\sigma(i) \in \cK, \, \LLR_i > t) \le \delta.$
By substituting this bound and the bound \eqref{eq:MidJulyBound3GenNew} to the right side above, we see that
\begin{align*}
 \pr( \cE_i, \, \sigma(i) \in \cK)
 &\wge e^{-t} ( 3\delta  - \delta )
 \weq 2\delta e^{-t}.
\end{align*}
%
By \eqref{eq:MidJulyBound1GenNew} it now follows that
\[
 \E \nCrit
 \wge \frac{\epsilon}{1+\epsilon} \sum_i \pr( \cE_i )
 \wge \frac{\epsilon}{1+\epsilon} \sum_i \pr( \cE_i, \, \sigma(i) \in \cK ),
\]
so that
\[
 \E \nCrit
 \wge \frac{2 N \delta e^{-t}}{1+\epsilon^{-1}}.
\]
Because $1+\epsilon^{-1} \le \frac{7}{6} \epsilon^{-1} = 7 \alphamin^{-1}$, the claim follows.
\end{proof}

\subsubsection{Concluding the proof of Theorem \ref{the:LowerBoundGen}}

By Lemma~\ref{the:MeanPermutedHammingLowerBoundGen}, we find that
\[
 \E L
 \wge \frac{\alphamin}{6} \left( \E \nCrit - N K e^{- \frac18 N \alphamin} \right).
\]
By Lemma~\ref{the:LowerBoundCriticalNodeCountGen}, 
\[
 \E \nCrit
 \wge \frac{2}{7} \alphamin \delta Ne^{-t}
\]
for $t = \max_i \Big( \tilde\E^{*i} \LLR_i + \alpha_\cK^{1/2} \delta^{-1/2} \sqrt{\tilde\V^{*i}(\LLR_i)} \Big)$.
By Lemma~\ref{the:LLRMeanVarGen}, 
$\tilde \E^{*i} \LLR_i \le N I_1$ and $\tilde \V^{*i} (\LLR_i) \le N I_{21} + N^2 I_{22}$, so that
$t \le N I_1 + \alpha_\cK^{1/2} \delta^{-1/2} \sqrt{N I_{21} + N^2 I_{22}}$. By combining these facts, it follows that
\begin{align*}
 \E L
 &\wge \frac{\alphamin}{6} \left( \E \nCrit - N K e^{- \frac18 N \alphamin} \right) \\
 &\wge \frac{\alphamin}{6} \left( \frac{2}{7} \alphamin \delta Ne^{-t}
 - N K e^{- \frac18 N \alphamin} \right) \\
 &\wge \frac{\alphamin}{6}
 \left( \frac{2}{7} \alphamin \delta N e^{-N I_1 - \alpha_\cK^{1/2} \delta^{-1/2} \sqrt{N I_{21} + N^2 I_{22}}}
 - N K e^{- \frac18 N \alphamin} \right).
\end{align*}
Hence the claim of Theorem \ref{the:LowerBoundGen} is valid.
\qed

\subsection{Application to homogeneous models}
\label{subsection:lower_bound_application_homogeneous_model}

\subsubsection{Log-likelihood ratio in homogeneous models}

The expected log-likelihood ratio equals $(N-1)I_1$ where $I_1$ is given in~\eqref{eq:I1I21I22}. The following result shows how to minimise this in the homogeneous case with intra-block and inter-block interaction distributions $f$ and~$g$.

\begin{lemma}
\label{the:MinimumI1}
For any homogeneous SBM and for any $\cK \subset [K]$ of size at least two such that $\alpha_k > 0$ for all $k \in \cK$,
\begin{equation}
 \label{eq:MinExpectedLLRHomNew}
 \min_{f^*_1,\dots, f^*_K} I_1
 \weq \sum_{k \in \cK} \alpha^*_k \alpha_k D_{1-\alpha_k^*}( g \| f ),
\end{equation}
with $\alpha_k^* = \alpha_k / (\sum_{k \in \cK} \alpha_k)$, and the minimum is attained by setting
\begin{equation}
 \label{eq:OptimalDistributions}
 f^*_k
 \weq
 \begin{cases}
  Z_{\alpha_k^*}^{-1} f^{\alpha_k^*} g^{1-\alpha_k^*} &\quad \text{for $k \in \cK$}, \\
  g, &\quad \text{otherwise}.
 \end{cases}
\end{equation}
Furthermore, when $\alpha$ is the uniform distribution on $[K]$,
\begin{equation}
 \label{eq:MinExpectedLLRHom}
 \min_{\cK: \abs{\cK} \ge 2} \min_{f^*_1,\dots, f^*_K} I_1
 \weq K^{-1} D_{1/2}(f \| g ).
\end{equation}
\end{lemma}
\begin{proof}
Observe that $I_1 = I_{11} + I_{12}$ where
\[
 I_{11} = \sum_{\ell \in \cK} \alpha_\ell \sum_{k \in \cK} \alpha^*_k \dkl(f^*_\ell \| f_{k\ell}) 
 \quad\text{and}\quad
 I_{12} = \sum_{\ell \in \cK^c} \alpha_\ell \sum_{k \in \cK} \alpha^*_k \dkl(f^*_\ell \| f_{k\ell}).
\]
We see that
\begin{align*}
 I_{11}
 &\weq \sum_{\ell \in \cK} \alpha_\ell \Big( \alpha^*_\ell \dkl(f^*_\ell \| f)
   + ( 1 - \alpha^*_\ell ) \dkl(f^*_\ell \| g) \Big)
\end{align*}
and
\[
 I_{12} \weq \sum_{\ell \in \cK^c} \alpha_\ell \dkl(f^*_\ell \| g).
\]
Because each $f^*_\ell$ appears only once in the sums above, we minimise $I_{11}$ and $I_{12}$ separately. To minimise $I_{12}$, we set $f^*_\ell = g$ for all $\ell \in \cK^c$, leading to $I_{12} = 0$. To minimise $I_{11}$, we see by applying \cite[Theorem 30]{vanErven_Harremoes_2014} that for all $\ell \in \cK$,
\[
 \min_{f^*_\ell} \Big( \alpha^*_\ell \dkl(f^*_\ell \| f) + ( 1 -  \alpha^*_\ell) \dkl(f^*_\ell \| g) \Big)
 \weq (1-\alpha^*_\ell) D_{\alpha^*_\ell}( f \| g ),
\]
and the minimum is attained by setting $f^*_\ell$ as in \eqref{eq:OptimalDistributions}. Hence the minimum value of $I_1$ equals
\[
 I_1 \weq \sum_{\ell \in \cK} \alpha_\ell (1-\alpha^*_\ell) D_{\alpha^*_\ell}( f \| g ).
\]
Finally, by skew symmetry of \Renyi divergences, we know that $(1-\alpha^*_\ell) D_{\alpha^*_\ell}( f \| g ) = \alpha^*_\ell D_{1-\alpha^*_\ell}( g \| f )$, so that we can also write the minimum as
\[
 I_1
 \weq \sum_{\ell \in \cK} \alpha_\ell \alpha^*_\ell D_{1-\alpha^*_{\ell}}( g \| f )
 \weq \alpha_\cK^{-1} \sum_{\ell \in \cK} \alpha_\ell^2 D_{1-\alpha^*_{\ell}}( g \| f ).
\]

Assume now that $\alpha$ is the uniform distribution on $[K]$.  Then the minimum above equals
$I_1 = (K/r) K^{-2} r D_{1-1/r}( g \| f )
= K^{-1} D_{1-1/r}( g \| f )$
for $r = \abs{\cK}$. Because $r \mapsto D_{1 - \frac{1}{r}}(g \| f)$ is increasing in $r$, we see that $I_1$ is increasing as a function of $\abs{\cK}$.  The minimum with respect to $\cK$ is hence attained at an arbitrary $\cK$ with $\abs{\cK} = 2$, confirming \eqref{eq:MinExpectedLLRHom}.
\end{proof}

The following result describes the variance terms $I_{21}$ and $I_{22}$ given by \eqref{eq:I1I21I22} for a uniform homogeneous SBM, when the reference distributions $f^*_1,\dots, f^*_K$ are selected to minimise $I_1$ according to Lemma~\ref{the:MinimumI1}.

\begin{lemma}
\label{the:MinimumI2}
Consider a homogeneous SBM with intra-block and inter-block interaction distributions $f$ and $g$, and uniform $\alpha$ on $[K]$. Fix $\cK \subset [K]$ of size 2, and define
$f^*_\ell$ as in \eqref{eq:OptimalDistributions}.
Then
\begin{align*}
 I_{21} &\weq \left(\frac12-K^{-1} \right) K^{-1} I^2 + \frac12 K^{-1} J, \\
 I_{22} &\weq 0,
\end{align*}
where $I = D_{1/2}(f \| g)$ and $J = \int h \log^2 \frac{f}{g}$ with $h = Z_{1/2}^{-1} (f g)^{1/2}$.
\end{lemma}
\begin{proof}
When $\alpha$ is uniform on $[K]$ and $\abs{\cK} = 2$, we see that the distributions in \eqref{eq:OptimalDistributions} are given by $f^*_\ell = h$ for $\ell \in \cK$, $f^*_\ell = g$ otherwise. Recall that
\begin{align*}
 I_{21} &\weq \sum_{k \in \cK} \sum_{\ell} \alpha^*_k \alpha_\ell \vkl(f^*_\ell \| f_{k\ell})
 + \sum_{k \in \cK} \alpha^*_k B_k, \\
 I_{22} &\weq \sum_{k \in \cK} \alpha^*_k A_k^2 - \Big(\sum_{k \in \cK} \alpha^*_k A_k \Big)^2,
\end{align*}
with $A_k = \sum_\ell \alpha_\ell \dkl(f^*_\ell \| f_{k\ell})$ and $B_k = \sum_\ell \alpha_\ell \dkl(f^*_\ell \| f_{k\ell})^2 - ( \sum_\ell \alpha_\ell \dkl(f^*_\ell \| f_{k\ell}))^2$. Now for any $k \in \cK$, we have by a direct computation (or using~\cite[Theorem~30]{vanErven_Harremoes_2014})
\[
 A_k
 \weq K^{-1} \Big( \dkl(h \| f) + \dkl(h \| g) \Big)
 \weq K^{-1} I.
\]
This implies that $I_{22}=0$.

Observe next that for $k \in \cK$,
\begin{align*}
 B_k
 &\weq \sum_\ell \alpha_\ell \dkl(f^*_\ell \| f_{k\ell})^2 - A_k^2 \\
 &\weq K^{-1} \Big( \dkl(h \| f)^2 + \dkl(h \| g)^2 \Big) - K^{-2} I^2.
\end{align*}
Because $\log Z = - \frac12 I$, we find that $\log \frac{h}{f} = \frac12 I - \frac12 \log \frac{f}{g}$ and
$\log \frac{h}{g} = \frac12 I + \frac12 \log \frac{f}{g}$. By squaring these equalities and integrating against $h$, we find that
\[
 \vkl(h \| f) + \vkl(h \| g)
 \weq \frac12 I^2 + \frac12 J - \dkl(h \| f)^2 - \dkl(h \| g)^2.
\]
It follows that
\begin{align*}
 \sum_{k \in \cK} \sum_{\ell} \alpha^*_k \alpha_\ell \vkl( f^*_\ell \| f_{k\ell})
 &\weq \sum_{k \in \cK} \sum_{\ell \in \cK} \alpha^*_k \alpha_\ell \vkl(h \| f_{k\ell}) \\
 &\weq \frac12 K^{-1} \sum_{k \in \cK} ( \vkl(h \| f) + \vkl(h \| g) ) \\
 &\weq K^{-1} \Big( \vkl(h \| f) +  \vkl(h \| g)  \Big) \\
 &\weq K^{-1} \Big( \frac12 I^2 + \frac12 J - \dkl(h \| f)^2 - \dkl(h \| g)^2 \Big).
\end{align*}
Therefore,
\begin{align*}
 I_{21}
 &\weq \sum_{k \in \cK} \sum_{\ell} \alpha^*_k \alpha_\ell \vkl( f^*_\ell \| f_{k\ell})
 + \sum_{k \in \cK} \alpha^*_k B_k \\
 &\weq K^{-1} \Big( \frac12 I^2 + \frac12 J - \dkl(h \| f)^2 - \dkl(h \| g)^2 \Big) \\
 &\qquad + K^{-1} \Big( \dkl(h \| f)^2 + \dkl(h \| g)^2 \Big) - K^{-2} I^2 \\
 &\weq \left( \frac12-K^{-1} \right) K^{-1} I^2 + \frac12 K^{-1} J.
\end{align*}
\end{proof}

 \begin{lemma}
 \label{lemma:bounding_J_over_I}
  Let  $I = \drenh(f,g) = -2 \log Z$ and $J = Z^{-1} \int \log^2(f/g) \sqrt{fg}$, where $Z = \int \sqrt{fg}$.
  Assume that $f,g > 0$ on $S$, and that $Z > 0$. Then
    \[
     J
     \wle 8 ( e^{I/2} - 1 ).
    \]
    Especially, $J \le 14 I$ whenever $I \le 1$.
    \end{lemma}

    \begin{proof}
    Let us fix some $x \in S$ for which $f(x) \ne g(x)$. At this point, for $t = \sqrt{f / g}$,
    \[
     \frac{(\log f - \log g)^2}{(\sqrt{f} - \sqrt{g})^2} \sqrt{fg}
     \weq 4 \frac{(\log \sqrt{f} - \log \sqrt{g})^2}{(\sqrt{f} - \sqrt{g})^2} \sqrt{fg}
     \weq 4 \phi(t)
    \]
    where $\phi(t) = \frac{(\log t)^2}{(t-1)^2} \, t$. Assume that $t > 1$, and let $u = \frac12 \log t$. Then $t = e^{2u}$ and
    \[
     \phi(t)
     \weq \left(\frac{2u}{e^{2u}-1} \right)^2 e^{2u}
     \weq \left( \frac{2u}{e^{u}-e^{-u}} \right)^2
     \weq \left( \frac{u}{\sinh u} \right)^2,
    \]
    where
    \[
     \sinh u
     \weq \frac12 ( e^u - e^{-u} )
     \weq \sum_{k > 0, \rm{odd}} \frac{u^k}{k!}
     \wge u.
    \]
    Hence $\phi(t) \le 1$ for all $t > 1$. Next, by noting that $\phi(t) = \phi(1/t)$ for all $0 < t$, we conclude that $\phi(t) \le 1$ for all $t > 0$ such that $t \ne 1$. We conclude that 
    \[
     (\log f - \log g)^2 \sqrt{fg}
     \wle 4 (\sqrt{f} - \sqrt{g})^2
    \]
    whenever $f \ne g$. Obviously the same inequality holds also when $f=g$. By integrating both sides, it follows that
    \[
     ZJ
     \wle 4 \int (\sqrt{f} - \sqrt{g})^2
     \weq 4 ( 2 - 2 Z )
     \weq 8 (1-Z).
    \]
    Hence $J \le 8( Z^{-1} - 1)$. The first claim follows because $Z = e^{-I/2}$. The second claim follows by noting that $e^{t/2}-1 = \int_0^{t/2} e^s ds \le e^{1/2} t$ for $t \le 1$, and $8 e^{1/2} \le 14$.
    \end{proof}
    
\subsubsection{Lower bound for homogeneous models}

\begin{proposition}
\label{the:LowerBoundHom}
Consider a stochastic block model defined by \eqref{eq:PairwiseInteractionModel}--\eqref{eq:JointProbability}. Suppose that $\alpha$ is the uniform distribution over $[K]$, and that the interactions are homogeneous. Then for any estimator $\hsigma: \cX \to \cZ$, the error 
is bounded in expectation by
\[
 \E \left( \frac{\dhams(\hsigma)}{N} \right)
 \wge \frac{1}{84} K^{-3} e^{ - \frac{N}{K} I - \sqrt{ 8 N I_{21} } } - \frac16 e^{ - \frac{N}{8K} }
\]
where $I_{21} = \left(\frac12-K^{-1} \right) K^{-1} I^2 + \frac12 K^{-1} J$.
\end{proposition}

\begin{proof}
Theorem~\ref{the:LowerBoundGen} states that
\begin{equation}
 \E \Ham^*(\hsigma) \wge \frac{1}{21} N \alphamin^2 \delta e^{-N I_1 - \alpha_\cK^{1/2} \delta^{-1/2} \sqrt{N I_{21} + N^2 I_{22}}} - \frac16 N \alphamin K e^{- \frac18 N \alphamin}.
\end{equation}
Lemma~\ref{the:MinimumI1} implies that $\min_{\cK: \abs{\cK} \ge 2} \min_{f^*_1,\dots, f^*_K} I_1 \weq K^{-1} D_{1/2}(f \| g )$. 
When the minimum is achieved,  Lemmas~\ref{the:MinimumI2} and~\ref{lemma:bounding_J_over_I} ensure that $I_{22} =0$ and 
 $I_{21} = \left(\frac12-K^{-1} \right) K^{-1} I^2 + \frac12 K^{-1} J$. Furthermore, we have $\alpha_{\cK} = \frac2K$ and $\delta = \frac14 \left(\frac2K - \frac1K \right) = \frac{1}{4K}$ since $\alpha$ is uniform.
\end{proof}

\section{Upper bound on ML estimation error}
\label{appendix:mle_consistency}

This section is devoted to analysing the accuracy of maximum-likelihood estimators.
Section~\ref{sec:MLE} describes how ML estimation error probabilities are characterised by the Mirkin distance.
Section~\ref{sec:ML_UB_Balanced} provides an upper bound on a worst-case ML estimation error among balanced block structures.
Section~\ref{sec:ML_UB_Averaged} provides an upper bound (Proposition~\ref{the:MLEUpperAverage}) on an average ML estimation error among all block structures, which confirms the upper bound of Theorem~\ref{thm:asymptotic_exponential_bound_recovery}, and also shows that any maximum-likelihood estimator achieves the upper bound.
Section~\ref{sec:ML_Consistency} analyses the upper bound of Theorem~\ref{thm:asymptotic_exponential_bound_recovery} in a large-scale setting and yields a proof of the existence part of  Theorem~\ref{cor:recovery_conditions}, summarised as Proposition~\ref{the:MLEConsistent}.

\subsection{Maximum likelihood estimators}
\label{sec:MLE}

A maximum likelihood estimator of $\sigma$ is a map $\hsigma: \cX \to \cZ$ such that
\begin{equation}
 \label{eq:MLE}
 P_{\hsigma_x}(x) \ge P_{\sigma'}(x) \qquad \text{for all $\sigma' \in \cZ$ and $x \in \cX$}.
\end{equation}
%
%
%
%
The following results help us to analyse situations in which a maximum likelihood estimator produces outputs diverging from the correct value. The result is stated using the Mirkin distance $\Mir(\sigma,\sigma')$ defined in Section~\ref{sec:Mirkin}.

\begin{lemma}
\label{the:MirkinBound}
For a homogeneous SBM with $N$ nodes, $K$ blocks, and interaction distributions $f$ and $g$ with 
$I = \drenh(f,g)$
\[
 P_\sigma\{x: P_{\sigma'}(x) \ge P_\sigma(x) \}
 \wle e^{-\frac14 \Mir(\sigma, \sigma') I},
\]
for all node labellings $\sigma, \sigma'$.
\end{lemma}
\begin{proof}
Observe that $P_\sigma\{x: P_{\sigma'}(x) \ge P_\sigma(x) \} = P_\sigma(\ell \ge 0)$, where the log-likelihood ratio $\ell(x) = \log \frac{P_{\sigma'}(x)}{P_\sigma(x)}$ is viewed as a random variable on probability space $(\cX, P_\sigma)$. Also denote by $E$ (resp.\ $E'$) the set of node pairs $\{i,j\}$ for which $\sigma(i) = \sigma(j)$
(resp.\ $\sigma'(i) = \sigma'(j)$). Then we find that
\[
 \ell(x)
 \weq \sum_{ij \in E' \setminus E} \log \frac{f}{g}(x_{ij}) 
 \ \ - \sum_{ij \in E \setminus E'} \log \frac{f}{g}(x_{ij}).
\]
Therefore, the distribution of $x \mapsto \ell(x)$
on the probability space $(\cX, P_\sigma)$ is the same as the law of
\[
 \sum_{j=1}^{\abs{E'\setminus E}} \log\frac{f}{g}(Y_j)
 - \sum_{i=1}^{\abs{E\setminus E'}} \log\frac{f}{g}(X_i),
\]
in which the random variables $X_i,Y_j$ are mutually independent and distributed according to $\law(X_i) = f$ and $\law(Y_j) = g$.  By applying Markov's inequality and the above representation, we find that
\[
 P_\sigma(\ell \ge 0)
 \weq P_\sigma( e^{\frac12 \ell} \ge 1)
 \wle E_\sigma e^{\frac12 \ell}
 \weq e^{-\frac12 ( \abs{E' \setminus E} + \abs{E \setminus E'} ) I },
\]
where $I = \drenh(f,g)$.  Hence the claim follows. 
\end{proof}

\subsection{Upper bound on worst-case error among balanced node labellings}
\label{sec:ML_UB_Balanced}

The following result is key minimax upper bound characterising the worst-case estimation accuracy
among block structures which are balanced according to $\sigma \in \cZ_{1-\epsilon, 1+\epsilon}$, where
\begin{equation}
 \label{eq:Balanced}
 \cZ_{a,b} 
 \weq \left\{ \sigma \in \cZ: \ a \frac{N}{K} \le \abs{\sigma^{-1}(k)} \le b \frac{N}{K} \right\},
\end{equation}
and we recall that $\cZ = [K]^{[N]}$. Similar upper bounds in the context of binary SBMs have been derived in \cite{Zhang_Zhou_2016}.

\begin{proposition}
\label{the:MLEUpperMinimax}
For a homogeneous SBM with $N$ nodes and $K$ blocks, any estimator $\hsigma: \cX \to \cZ$ satisfying the MLE property \eqref{eq:MLE} has classification error bounded by
\[
 \max_{\sigma \in \cZ_{1-\epsilon,1+\epsilon}} E_\sigma \ace(\sigma, \hsigma)
 \wle 8 e N (K-1) e^{-(1 - \zeta - \kappa) \nik}
 + N K^N e^{-\frac14 (\frac{\zeta}{K-1} - \epsilon) (N/K)^2 I }
\]
for all $0 \le \epsilon \le \zeta \le \frac{1}{21}$, where
$\kappa = 56 \max\{ K^2 e^{-\frac18 \nik}, \, K N^{-1} \}$
and $I = \drenh(f,g)$.
\end{proposition}

\begin{proof}
We note that due to homogeneity, $P_\sigma = P_{[\sigma]}$ depends on $\sigma$ only via the partition $[\sigma] = \{\sigma^{-1}(k): k \in [K]\}$.  A similar observation also holds for the absolute classification error $\ace(\sigma_1, \sigma_2) = \ace([\sigma_1], [\sigma_2])$.  In the proof we denote by $\cP_{1-\epsilon, 1+\epsilon} = \{ [\sigma]: \sigma \in \cZ_{1-\epsilon, 1+\epsilon}\}$ the collection of partitions corresponding to node labellings in $\cZ_{1-\epsilon, 1+\epsilon}$. We select a node labelling $\sigma \in \cZ_{1-\epsilon,1+\epsilon}$, and split the error according to
\begin{equation}
 \label{eq:AverageErrorRateSplitNew}
 E_\sigma L
 \weq
 E_\sigma L 1( \hsigma \in \cZ_{1-\zeta, 1+\zeta} )
 + E_\sigma L 1( \hsigma \notin \cZ_{1-\zeta, 1+\zeta} ).
\end{equation}
The remainder of the proof consists of two parts, where we derive upper bounds for both terms on the right side above.

(i) For analysing the first term on the right side of \eqref{eq:AverageErrorRateSplitNew},
we note that $\hsigma \in \cZ_{1-\zeta, 1+\zeta}$ if and only if $[\hsigma] \in \cP_{1-\zeta, 1+\zeta}$,
and therefore,
\begin{equation}
 \label{eq:AverageErrorRateSplit1}
 E_\sigma L 1( \hsigma \in \cZ_{1-\zeta, 1+\zeta} )
 \weq \sum_{m=1}^N m p_m
\end{equation}
where $p_m = P_\sigma\{x: [\hsigma_x] \in \cP_{1-\zeta, 1+\zeta}(\sigma,m)\}$ is the probability of the event that the partition associated to $\hsigma_x$ belongs to the set
\[
 \cP_{1-\zeta, 1+\zeta}(\sigma,m)
 \weq \{ \theta \in \cP_{1-\zeta, 1+\zeta}: \ace([\sigma], \theta) = m\}.
\]
On such event  
%
there exists a partition $\theta \in \cP_{1-\zeta, 1+\zeta}(\sigma, m)$ such that
$P_{\theta}(x) \ge P_{[\sigma]}(x)$.  Hence by the union bound,
\[
 p_m
 \wle \sum_{\theta \in \cP_{1-\zeta, 1+\zeta}(\sigma, m)}
 P_\sigma\{ x: P_{\theta}(x) \ge P_{[\sigma]}(x) \}.
\]
Observe next that to every partition $\theta \in \cP_{1-\zeta, 1+\zeta}(\sigma, m)$ there corresponds exactly $K!$ node labellings $\sigma'$ belonging to the set
\[
 \cZ_{1-\zeta, 1+\zeta}(\sigma,m)
 \weq \{ \sigma' \in \cZ_{1-\zeta, 1+\zeta}: \ace(\sigma, \sigma') = m\}.
\]
Therefore, the above upper bound can be rewritten as
\begin{equation}
 \label{eq:AmUpperBound}
 p_m
 \wle (K!)^{-1}
 \nhquad \sum_{\sigma' \in \cZ_{1-\zeta, 1+\zeta}(\sigma,m)} \nhquad
 P_\sigma\{ x: P_{\sigma'}(x) \ge P_{\sigma}(x) \}.
\end{equation}

Let us next analyse the probabilities on the right side of \eqref{eq:AmUpperBound}.
By Lemma~\ref{the:MirkinBound}, we find that
\[
 P_\sigma\{ x: P_{\sigma'}(x) \ge P_\sigma(x) \}
 \wle e^{-\frac14 \Mir(\sigma, \sigma') I}.
\]
Because $\epsilon \le \zeta$, it follows that $\cZ_{1-\epsilon, 1+\epsilon} \subset \cZ_{1-\zeta, 1+\zeta}$. We note that
$\frac14 \Mir(\sigma, \sigma')
= \frac12 ( \abs{E\setminus E'} + \abs{E' \setminus E} )
\ge \min\{ \abs{E\setminus E'},  \abs{E' \setminus E}\}$,
where $E$ (resp., $E'$) denotes the set of node pairs for which $\sigma$ (resp., $\sigma'$) assigns the same label. With the help of Lemma~\ref{the:PartitionPairBound} we then find that for all $\sigma, \sigma' \in \cZ_{1-\zeta, 1+\zeta}$, such that $\ace(\sigma, \sigma') = m$,
\[
 \frac14 \Mir(\sigma, \sigma')
 \wge \max\left\{ (1-\zeta) \fracnk - m, \
 \frac13 (1-\zeta) \fracnk - \frac16 (1+\zeta) \fracnk \right\}m.
\]
We note that
$\frac13 (1-\zeta) - \frac16 (1+\zeta) = \frac16 - \frac12 \zeta \ge \frac17$
when $\zeta \le \frac{1}{21}$. Hence,
\[
 \frac14 \Mir(\sigma, \sigma')
 \wge \max\left\{ (1-\zeta) \fracnk - m, \ \frac17 \fracnk \right\}m,
\]
and we conclude that for all $\sigma \in \cZ_{1-\epsilon, 1+\epsilon}$ and $\sigma' \in \cZ_{1-\zeta, 1+\zeta}$,
\begin{equation}
 \label{eq:MLEBoundHam}
 P_\sigma\{ x: P_{\sigma'}(x) \ge P_\sigma(x)\}
 \wle
 \min \left\{ e^{- (1-\zeta) \nik+ m I}, \ e^{-\frac17 \nik} \right\}^m.
\end{equation}
Furthermore, let us analyse the cardinality of the sum on \eqref{eq:AmUpperBound}.
Because $\ace(\sigma, \sigma') = m$ if and only if $\Ham( \tau\circ\sigma, \sigma')=m$ for some $\tau\in \Sym(K)$, a union bound combined with Lemma~\ref{the:HamBound} implies that
\begin{align*}
 \abs{\cZ_{1-\zeta, 1+\zeta}(\sigma,m)}
 \wle K! \, \abs{ \{ \sigma' \in \cZ: \Ham(\sigma, \sigma') = m \} }
 \wle K! \left( \frac{eN(K-1)}{m} \right)^m.
\end{align*}
By combining this bound with \eqref{eq:AmUpperBound} and \eqref{eq:MLEBoundHam}, 
we may now conclude that
\begin{equation}
 \label{eq:ProbBoundNew}
 p_m
 \wle \min \left\{ \frac{eN(K-1)}{m} e^{- (1-\zeta) \nik + m I}, \ 
 \frac{eN(K-1)}{m} e^{- \frac17 \nik} \right\}^m.
\end{equation}

We will now apply the bounds in \eqref{eq:ProbBoundNew} to derive an upper bound for the sum in \eqref{eq:AverageErrorRateSplit1} which we will split according to
\begin{equation}
 \label{eq:AverageErrorRateSplit1Detail}
 \sum_{m=1}^N m p_m
 \weq
 \sum_{m \le m_1} m p_m
 \ + \sum_{m_1 < m \le N} \nhquad m p_m
\end{equation}
using a threshold parameter $m_1$.  We will also select another threshold parameter $0 < m_0 \le m_1$.  Using these, the probabilities $p_m$ are bounded by $p_m \le s_1^m$ for $m_0 \le m \le m_1$, and $p_m \le s_2^m$ for $m \ge m_1$,
where
\begin{align*}
 s_1 = \frac{eN(K-1)}{m_0} e^{- (1-\zeta) \nik + m_1 I}
 \qquad\text{and}\qquad
 s_2 = \frac{eN(K-1)}{m_1} e^{- \frac17 \nik }.
\end{align*}
To obtain a good upper bound, $m_1$ should be small enough to keep the exponent in $s_1$ small, and large enough so that $s_2 < 1$.
From the latter point of view, we see that $s_2 \le \frac12$ when 
$m_1 \ge 2 eN (K-1) e^{-\frac17 \nik}$. To leave some headroom, we set a slightly larger $m_1$ corresponding to $\frac17$ replaced by $\frac18$. For later purposes, we also require that $m_1 \ge 56$ which guarantees that
$\frac{1}{m_1} \le \frac17 - \frac18$. Therefore, we set
\[
 m_1 \weq 2 eN (K-1) e^{-\frac18 \nik} \vee 56.
\]
With this choice, we find that $s_2 \le \frac12 e^{-\frac{1}{56} \nik} \le \frac12$. Hence,
\[
 \sum_{m_1 < m \le N} m p_m
 \wle N \sum_{m \ge m_1} s_2^m
 \weq N \frac{s_2^{\ceil{m_1}}}{1-s_2}
 \wle N \frac{s_2^{m_1}}{1-s_2}
 \wle 2 N s_2^{m_1}.
\]
Furthermore, $m_1 \ge 56$ implies that $s_2^{m_1} \le e^{-\frac{m_1}{56} \nik} \le e^{-\nik}$. It follows that the second term on the right side of \eqref{eq:AverageErrorRateSplit1Detail} is bounded by
\begin{equation}
 \label{eq:TailBoundNew}
 \sum_{m_1 < m \le N} m p_m
 \wle 2 N e^{-\nik}.
\end{equation}

Let us next derive an upper bound for the first term on the right side of \eqref{eq:AverageErrorRateSplit1Detail}. We define $B = eN(K-1) e^{-(1-\zeta) \nik + m_1 I}$, and consider the following two cases. \\
(a) If $B \le \frac12$, we set $m_0=1$, which implies that $s_1 = B$, and we find that
\begin{equation}
 \label{eq:HiInformationBound}
 \sum_{1 \le m \le m_1} m p_m
 \wle \sum_{m=1}^\infty m s_1^m
 \weq \sum_{m=1}^\infty m B^m
 \weq \frac{B}{(1-B)^2}
 \wle 4 B.
\end{equation}
(b) If $B > \frac12$, we set $m_0 = 2B$, so that $s_1 = \frac12$, and we find that
\begin{align*}
 \sum_{1 \le m \le m_1} \nhquad m p_m
 \weq \sum_{1 \le m \le m_0} \nhquad m p_m + \sum_{m_0 < m \le m_1} \nhquad m p_m 
 \wle m_0 \ + \sum_{m > m_0} \nhquad m s_1^m.
\end{align*}
By noting that $m_0 > 1$, we find that
$2 \le \floor{m_0}+1 \le 2m_0$.  Then by applying Lemma~\ref{the:PowerSeriesBound} it follows that
\[
 \sum_{m > m_0} m s_1^m
 \weq \sum_{m=\floor{m_0}+1}^\infty m 2^{-m}
 \wle 4 (\floor{m_0}+1) 2^{-(\floor{m_0}+1)}
 \wle 2 m_0.
\]
Hence, $\sum_{1 \le m \le m_1} m p_m \le 3 m_0 = 6B$.  In light of \eqref{eq:HiInformationBound}, we conclude that the latter conclusion holds for both $B \le \frac12$ and $B > \frac12$.  By combining these observations with \eqref{eq:TailBoundNew}, and noting that $B \ge N e^{-\nik}$, it follows that
\begin{align*}
 \sum_{1 \le m \le N} m p_m
 \wle 2 N e^{-\nik} + 6B
 \wle 8 B
 \weq 8 eN(K-1) e^{-(1-\zeta) \nik + m_1 I}.
\end{align*}
After noting that $m_1 I = \nik \max\{ 2 e K(K-1) e^{-\frac18 \nik}, \, 56 \frac{K}{N}\}$, we see that $m_1 I \le \kappa \nik$ for $\kappa = 56 \max\{ K^2 e^{-\frac18 \nik}, \, K N^{-1} \}$. Then we conclude that
the first term on the right side of \eqref{eq:AverageErrorRateSplitNew} is bounded by
\begin{equation}
 \label{eq:MLEUpperBoundPart1}
 E_\sigma L 1( \hsigma \in \cZ_{1-\zeta, 1+\zeta} )
 \wle 8 eN(K-1) e^{-(1-\zeta- \eta) \nik}.
\end{equation}

(ii) Finally, it remains to derive an upper bound for the second term on the right side of \eqref{eq:AverageErrorRateSplitNew}.  Denote $\gamma = (K-1)^{-1}\zeta$. Then the generic bound $N \le \Nmin(\sigma') +  (K-1) \Nmax(\sigma')$ implies that
$\Nmin(\sigma') \ge N - (K-1)(1+\gamma)\frac{N}{K} = (1-\zeta) \frac{N}{K}$ for all $\sigma' \in \cZ_{0,1+\gamma}$. Therefore,
$\cZ_{0,1+\gamma} \subset \cZ_{1-\zeta, 1+\gamma} \subset \cZ_{1-\zeta, 1+\zeta}$.
Especially,
\[
 P_\sigma( \hsigma \not\in \cZ_{1-\zeta, 1+\zeta})
 \wle P_\sigma( \hsigma \not\in \cZ_{0, 1+\gamma} ).
\]
On the event that $\hsigma \not\in \cZ_{0, 1+\gamma}$, the MLE property \eqref{eq:MLE} implies that there exists $\sigma'$ with $\Nmax(\sigma') > (1+\gamma)\fracnk$
for which $P_{\sigma'}(x) \ge P_\sigma(x)$. For any such $\sigma'$, 
$\Nmax(\sigma') -  \Nmax(\sigma) \ge (\gamma-\epsilon) \frac{N}{K}$, so that by
Lemma~\ref{the:PartitionPairBoundGeneral}, we see that
\[
 \Mir(\sigma, \sigma')
 \weq 2 ( \abs{E \setminus E'} + \abs{E' \setminus E} ) 
 \wge 2 \abs{E' \setminus E}
 \wge (\gamma-\epsilon) (N/K)^2.
\]
By Lemma~\ref{the:MirkinBound}, we conclude that for all $\sigma' \notin \cZ_{0,1+\gamma}$,
\[
 P_\sigma\{ x: P_{\sigma'}(x) \ge P_\sigma(x) \}
 \wle e^{-\frac14 (\gamma-\epsilon) (N/K)^2 I }.
\]
Hence, by the union bound it follows that
\begin{equation}
 \label{eq:AverageErrorRate2}
 P_\sigma( \hsigma \not\in \cZ_{1-\zeta, 1+\zeta})
 \wle P_\sigma( \hsigma \not\in \cZ_{0, 1+\gamma} )
 \wle K^N e^{-\frac14 (\gamma-\epsilon) (N/K)^2 I },
\end{equation}
and we conclude that second term on the right side of \eqref{eq:AverageErrorRateSplitNew} is bounded by
\begin{equation}
 \label{eq:MLEUpperBoundPart2}
 E_\sigma L 1( \hsigma \notin \cZ_{1-\zeta, 1+\zeta} )
 \wle N P_\sigma( \hsigma \not\in \cZ_{1-\zeta, 1+\zeta})
 \wle N K^N e^{-\frac14 (\gamma-\epsilon) (N/K)^2 I }.
\end{equation}
The claim now follows by combining \eqref{eq:MLEUpperBoundPart1}--\eqref{eq:MLEUpperBoundPart2}.
\end{proof}

\subsection{Upper bound on average error among all node labellings}
\label{sec:ML_UB_Averaged}

The following result is the upper bound of Theorem~\ref{thm:asymptotic_exponential_bound_recovery}.

\begin{proposition}
\label{the:MLEUpperAverage}
For a homogeneous SBM with $N$ nodes and $K$ blocks, any estimator $\hsigma: \cX \to \cZ$ satisfying the MLE property \eqref{eq:MLE} has classification error bounded by
\[
 \E \ace(\sigma, \hsigma)
 \wle 8 e N (K-1) e^{-(1 - \zeta - \kappa) \nik}
 + N K^N e^{-\frac14 (\frac{\zeta}{K-1} - \epsilon) (N/K)^2 I }
 + 2 N K e^{-\frac13 \epsilon^2 \frac{N}{K}},
\]
for all $0 \le \epsilon \le \zeta \le \frac{1}{21}$, where
$\kappa = 56 \max\{ K^2 e^{-\frac18 \nik}, \, K N^{-1} \}$
and $I = \drenh(f,g)$.
\end{proposition}

\begin{proof}
Denote $L = \ace(\sigma, \hsigma)$. By noting that the classification error is bounded by $L \le N$ with probability one, it follows that
\begin{align*}
 \E L
 &\wle \sum_{\sigma \in \cZ_{1-\epsilon,1+\epsilon}} \nhquad \pi_\sigma E_\sigma L
 \ + \sum_{\sigma \in \cZ_{1-\epsilon,1+\epsilon}^c} \nhquad \pi_\sigma E_\sigma L \\
 & \wle \max_{\sigma \in \cZ_{1-\epsilon,1+\epsilon}} E_\sigma L
 \ + N \pi( \cZ_{1-\epsilon,1+\epsilon}^c ).
\end{align*}
For a random node labelling $\sigma = (\sigma_1, \dots, \sigma_N)$ sampled from the uniform distribution $\pi$ on~$\cZ$, we see that coordinates are mutually independent and uniformly distributed on in $[K]$. A multinomial concentration inequality (Lemma~\ref{the:MultinomialConcentrationMinMaxDiff}) then implies that
\[
 \pi( \cZ_{1-\epsilon,1+\epsilon}^c )
 \wle 2 K e^{-\frac13 \epsilon^2 \frac{N}{K}}.
\]
The claim follows by Proposition~\ref{the:MLEUpperMinimax}.
\end{proof}

\subsection{Upper bound for large-scale settings}
\label{sec:ML_Consistency}

The following result implies the existence statements of Theorem~\ref{cor:recovery_conditions}.

\begin{proposition}
\label{the:MLEConsistent}
Consider a large-scale homogeneous SBM with $N \gg 1$ nodes and $K \asymp 1$ blocks, and interaction distributions $f,g$ such that $I = \drenh(f,g)$, and let $\hsigma$ be any estimator having the MLE property \eqref{eq:MLE}. 
\begin{enumerate}[(i)]
\item If $I \gg N^{-1}$, then the estimator $\hsigma$ is consistent in the sense that $\E \ace(\hsigma) = o(N)$.
\item If $I \ge (1+\Omega(1)) \frac{K \log N}{N}$, then the estimator $\hsigma$ is strongly consistent in the sense that $\E \ace(\hsigma) = o(1)$.
\end{enumerate}

\end{proposition}

\begin{proof}
Denote $L =  \ace(\hsigma)$. By Proposition~\ref{the:MLEUpperAverage}, we see that
\begin{equation}
 \label{eq:MLEAverageBound}
 \E L
 \wle 8 e N K e^{-(1 - \zeta - \eta) \nik}
 + N K^N e^{-\frac14 (\frac{\zeta}{K-1} - \epsilon) (N/K)^2 I }
 + 2 N K e^{-\frac13 \epsilon^2 \frac{N}{K}},
\end{equation}
where $\eta = 56 \max\{ K^2 e^{-\frac18 \nik}, \, K N^{-1} \}$
and $I = \drenh(f,g)$, and where we are 
free to choose any $0 \le \epsilon \le \zeta \le \frac{1}{21}$, 

(i) Suppose $I \gg N^{-1}$. Let us define $\epsilon = 3 ( \frac{K \log N}{N} )^{1/2} \ll 1$ and $\zeta = \epsilon K + 5 \frac{K^3 \log K}{NI} \ll 1$. 
We have $e^{-\frac13 \epsilon^2 \frac{N}{K}} = N^{-3}$, and the last term on the right side of \eqref{eq:MLEAverageBound} equals $2 KN^{-2}$. We also find that
\begin{align*}
 \left( \frac{\zeta}{K-1} - \epsilon \right) (N/K)^2 I 
 &\wge \left( \frac{\zeta}{K} - \epsilon \right) (N/K)^2 I 
 \weq 5 N \log K,
\end{align*}
so that the middle term on the right side of \eqref{eq:MLEAverageBound} is bounded by
\[
 N K^N e^{-\frac14 (\frac{\zeta}{K-1} - \epsilon) (N/K)^2 I }
 \wle N K^N K^{-\frac{5}{4}N}
 \weq N K^{-\frac{1}{4}N}.
\]
We conclude that
\[
 \E L
 \wle 8 e N K e^{-(1 - \zeta - \eta) \nik} + N K^{-\frac{1}{4}N} + 2 N^{-2} K.
\]
We note that $\eta \ll 1$ and $\log(8 e K) \ll \nik$. We note that $N K^{-\frac{1}{4}N} \le N 2^{-\frac{1}{4}N} = o(1)$. Hence we conclude that
\[
 \E L
 \wle N e^{-(1 - o(1)) \nik} + o(1).
\]

(ii) The condition for strong consistency follows immediately from the above bounds.
\end{proof}

\section{Consistency of Algorithm~\ref{algo:likelihood_based_algo}}
\label{appendix:consistency_likelihood_based_algo}

This section is devoted to the proof of Theorem~\ref{thm:general_algo_consistency} characterising the accuracy of Algorithm~\ref{algo:likelihood_based_algo}.
Section~\ref{the:SingleLabel} presents an upper bound for the estimation error of a conditional ML estimator.
Section~\ref{sec:Refinement_and_consensus} describes the analysis of refinement and consensus steps in Algorithm~\ref{algo:likelihood_based_algo}.
Section~\ref{sec:general_algo_consistency_proof} concludes the proof of Theorem~\ref{thm:general_algo_consistency}.

\subsection{Single node label estimation}
\label{the:SingleLabel}

Given a reference node $i$ and a node labelling\footnote{In this section we assume that $\tilde\sigma_i$ is nonrandom.} $\tilde\sigma_i$ on $[N] \setminus \{i\}$, define an estimator for the label of~$i$ by $\hsigma_i(i) = \argmax\limits_{k \in [K]} h_i(k)$, with arbitrary tie breaks, where
\begin{equation}
 \label{eq:SingleNodeLogLikelihoodRatio}
 h_i(k)
 \weq \sum_{j \colon \tilde\sigma_i(j) = k} \log \frac{ \fin(X_{ij})} {\fout(X_{ij})}.
\end{equation}
This is a maximum likelihood estimator in the special case where $\tilde\sigma_i$ assigns a correct label to all $j \ne i$.  When this is not the case, we need to account for errors caused by corrupted likelihoods due to misclassified nodes in $\tilde\sigma_i$. The error in such a setting is given by the following lemma. For $r > 0$ we define a ratio between symmetrised \Renyi divergences by
\begin{equation} 
 \label{eq:alphaRatio}
 \beta_r(f,g) \weq \frac{ \dren_{1+r}^s \left( f , g \right)  } { \dren_{r}^s \left( f, g \right) }.
\end{equation}

\begin{lemma}
\label{lemma:individual_error_renyi}
Let $\sigma: [N] \to [K]$ and assume that $X_{ij}$, $j \ne i$, are mutually independent $\cS$-valued random variables such that $\law(X_{ij}) = \fin$ for $\sigma(i) = \sigma(j)$ and $\law(X_{ij}) = \fout$ otherwise.
The error probability when estimating the label of node $i$ as a maximiser of \eqref{eq:SingleNodeLogLikelihoodRatio} is bounded by
\[
 \pr \left( \tau \circ \hat\sigma_i(i) \ne \sigma(i) \right)	
 \wle
 K e^{-( N_{\min} - 1 - ( 2 + \frac{r}{1-r} \beta_r) \dhamsi ) 2 (1-r) \dren_r^s(\fin, \fout)} \qquad \text{for all $r \in (0, 1)$},
\]
where $\beta_r = \beta_r(f,g)$ is defined by \eqref{eq:alphaRatio}, $\dhamsi = \dhams(\tilde\sigma_i, \sigma_{-i})$ is the symmetrised Hamming distance from $\tilde\sigma_i$ to the restriction $\sigma_{-i}$ of the true node labeling $\sigma$ to $[N] \setminus \{i\}$, and $\tau$ is an arbitrary $K$-permutation such that $\dham(\tau \circ \tilde\sigma_i, \sigma_{-i}) = \dhamsi$.
\end{lemma}

\begin{proof}
Denote $k^* = \tau^{-1}(\sigma(i))$. Observe that $\tau \circ \hat\sigma_i(i) \ne \sigma(i)$ if and only if $\hat\sigma_i(i) \ne k^*$, and the latter is possible only if $L_k = h_i(k) - h_i(k^*) \ge 0$ for some $k \ne k^*$.
Let us fix some $0 < r < 1$.
After noting that
$\pr(L_k \ge 0) = \pr\left( e^{r L_k} \ge 1 \right) \le \E e^{r L_k}$, it follows that
\begin{equation}
\label{eq:ExpMarkov}
\pr( \tau \circ \hat\sigma_i(i) \ne \sigma(i) )
\wle \sum_{k \ne k^*} \pr( L_k \ge 0)
\wle \sum_{k \ne k^*} \E e^{r L_k}.
\end{equation}

Denote by $C_k = \{j \ne i: \sigma(j) = k\}$ the peers of $i$ with true label $k$, and by $\tilde C_k = \{j \ne i: \tilde\sigma_i(j) = k\}$ the set of peers labelled $k$ by $\tilde\sigma_i$. Denote $Z_\alpha(f \| g) = \int f^\alpha g^{1-\alpha}$. By noting that for any $j \ne i$,
\[
\E \left( \frac{ \fin(X_{ij})} {\fout(X_{ij})} \right)^{r} 
\weq
\begin{cases}
Z_{1+r}(\fin \| \fout), &\quad \sigma(j) = \sigma(i), \\
Z_{r}(\fin \| \fout), &\quad \text{else}, \\
\end{cases}
\]
and
\[
\E \left( \frac{ \fin(X_{ij})} {\fout(X_{ij})} \right)^{-r} 
\weq
\begin{cases}
Z_{r}(\fout \| \fin), &\quad \sigma(j) = \sigma(i), \\
Z_{1+r}(\fout \| \fin), &\quad \text{else}, \\
\end{cases}
\]
we find that for all $k$, the log-likelihood ratio $h_i(k)$ defined in \eqref{eq:SingleNodeLogLikelihoodRatio} satisfies
\begin{align*}
\E e^{ r h_i(k)} &\weq Z_{ 1+r}(\fin \| \fout)^{\vin_k} 
Z_{r}(\fin \| \fout)^{\vout_k}, \\
\E e^{- r h_i(k)} &\weq Z_{1+r}(\fout \| \fin)^{\vout_k} Z_{r}(\fout \| \fin)^{\vin_k},
\end{align*}
where $\vin_k = \abs{\tilde C_k \cap C_{\sigma(i)}}$ and $\vout_k = \abs{\tilde C_k \setminus C_{\sigma(i)}}$. Because $h_i(k)$ and $h_i(\ell)$ are mutually independent for $k \ne \ell$, it follows that $L_k $ for $k \ne k^*$ satisfies
\begin{align*}
\E e^{r L_k}
&\weq Z_{1+r}(\fin \| \fout)^{\vin_k}
Z_{1+r}(\fout \| \fin)^{\vout_{k^*}} 
Z_{r}(\fin \| \fout)^{\vout_k}
Z_{r}(\fout \| \fin)^{\vin_{k^*} }.
\end{align*}
Because $Z_{r} = e^{-(1-r) \dren_r}$ and $Z_{1+r} = e^{r \dren_{1+r}}$, we may rephrase the above equality as $\E e^{r L_k} = e^{t}$, where
\[
t \weq s_1 \vin_k + s_2 \vout_{k^*} - u_1 \vout_k - u_2 \vin_{k^*},
\]
with 
$s_1 = r D_{1+r}(\fin \| \fout)$, 
$s_2 = r D_{1+r}(\fout \| \fin)$,
$u_1 = (1-r)\dren_r (\fin \| \fout ) $,
and
$u_2 = (1-r)\dren_r (\fout \| \fin)$.
By noting that $\vin_k + \vout_k = \abs{\tilde C_k}$, we see that
\begin{align*}
t & \weq \left( u_1 + s_1  \right) \vin_k + 
\left( u_2 + s_2 \right) \vout_{k^*} - u_1 \abs{\tilde C_k} - u_2 \abs{\tilde C_{k^*}}.
\end{align*}
One may verify that $\tau \circ \tilde\sigma_i(j) \ne \sigma(j)$ for all $j \in \tilde C_k \cap C_{\sigma(i)}$ and all $k \ne k^*$. Therefore, $\vin_k = \abs{C_k \cap C_{\sigma(i)}} \le \dham( \tau \circ \tilde\sigma_i, \sigma_{-i} )$.  Similarly, $\tau \circ \tilde\sigma_i(j) \ne \sigma(j)$ for all $j \in \tilde C_{k^*} \setminus C_{\sigma(i)}$ implies that $\vout_{k^*} = \abs{\tilde C_{k^*} \setminus C_{\sigma(i)}} \le  \dham( \tau \circ \tilde\sigma_i, \sigma_{-i} )$.
Next, by noting that $\tau \circ \tilde\sigma_i(j) \ne \tau(k)$ and $\sigma(j) = \tau(k)$ for $j \in C_{\tau(k)} \setminus \tilde C_k$, it follows that $\abs{ C_{\tau(k)} \setminus \tilde C_k } \le \dham( \tau \circ \tilde\sigma_i, \sigma_{-i})$. Therefore,
\begin{align*}
\abs{\tilde C_k}
\wge \abs{ \tilde C_k \cap C_{\tau(k)}} 
\weq \abs{ C_{\tau(k)}} - \abs{ C_{\tau(k)} \setminus \tilde C_k} 
\wge \Nmin - 1 - \dham( \tau \circ \tilde\sigma_i, \sigma_{-i}),
\end{align*}
and the above inequality also holds for $k = k^*$. By collecting the above inequalities and recalling that $\dham( \tau \circ \tilde\sigma_i, \sigma_{-i}) = \dhamsi$, we conclude that
\begin{align*}
t & \wle 
\dhamsi \left( u_1 + u_2 + s_1 + s_2  \right) - \left( N_{\min} - 1 - \dhamsi \right) (u_1 + u_2 ) \\
& \wle - (u_1 + u_2) \left( N_{\min} - 1 - 2 \dhamsi - \dhamsi \frac{s_1 + s_2}{u_1 + u_2} \right).
\end{align*}
The claim follows by observing that $u_1 + u_2 = 2 (1-r) \dren_r^{s}(\fin, \fout)$ and $s_1 + s_2 = 2 r \dren_{1+r}^{s} \left(\fin, \fout \right) $.
\end{proof}

\subsection{Analysis of refinement and consensus procedures}
\label{sec:Refinement_and_consensus}

Let us start with a lemma bounding difference between the block sizes given by two node labeling $\sigma_1,\sigma_2$ as a function of the Hamming distance.

\begin{lemma}
\label{the:HammingInverseImages}
For any $\sigma_1, \sigma_2: [N] \to [K]$, 
(i) $\big| \abs{\sigma_1^{-1}(k)} - \abs{\sigma_2^{-1}(k)} \big| \le \dham(\sigma_1,\sigma_2)$ for all $k$, and (ii) $\abs{ \Nmin(\sigma_1) - \Nmin(\sigma_2) } \le \dhams(\sigma_1,\sigma_2)$, where $\Nmin(\sigma_1) = \min_k \abs{\sigma_1^{-1}(k)}$ and $\Nmin(\sigma_2) = \min_k \abs{\sigma_2^{-1}(k)}$.
\end{lemma}
\begin{proof}
(i) Because $\abs{ \sigma_1^{-1}(k) \setminus \sigma_2^{-1}(k) } \le \dham(\sigma_1,\sigma_2)$, we find that
\begin{align*}
 \abs{\sigma_1^{-1}(k)}
 &\weq \abs{ \sigma_1^{-1}(k) \cap \sigma_2^{-1}(k) } + \abs{ \sigma_1^{-1}(k) \setminus \sigma_2^{-1}(k) } \\
 &\wle \abs{ \sigma_2^{-1}(k) } + \dham(\sigma_1,\sigma_2).
\end{align*}
By symmetry, the same inequality is true also with $\sigma_1,\sigma_2$ swapped.

(ii) Let $\tau$ be a $K$-permutation for which $\dham(\tau \circ \sigma_1,\sigma_2) = \dhams(\sigma_1,\sigma_2)$. Then by (i),
\[
 \abs{\sigma_2^{-1}(k)}
 \wge \abs{(\tau \circ \sigma_1)^{-1}(k)} - \dham(\tau \circ \sigma_1,\sigma_2)
 \wge \Nmin(\sigma_1) - \dhams(\sigma_1,\sigma_2).
\]
This implies that $\Nmin(\sigma_2) \ge \Nmin(\sigma_1) - \dhams(\sigma_1,\sigma_2)$. The second claim hence follows by symmetry.
\end{proof}

The following result describes the behaviour of Steps 2 and 3 in Algorithm~\ref{algo:likelihood_based_algo} on the event that Step 1 achieves moderate accuracy.

\begin{lemma}
\label{the:OptimalPermutations}
Assume that the outputs $\tilde\sigma_i$ of Step 1 in Algorithm~\ref{algo:likelihood_based_algo} satisfy $\dhams(\tilde\sigma_i, \sigma_{-i}) < \frac15 \Nmin - 1$ for all $i$. Then there exist unique $K$-permutations $\tau_1,\dots,\tau_N$ such that for all $i$:
\begin{enumerate}[(i)]
\item the outputs $\tilde\sigma_i$ of Step 1 satisfy $\dhams(\tilde\sigma_i, \sigma_{-i}) = \dham(\tau_i \circ \tilde\sigma_i, \sigma_{-i})$;
\item the outputs $\hat\sigma_i$ of Step 2 satisfy $\dhams(\hat\sigma_i, \sigma) = \dham ( \tau_i \circ \hat\sigma_i, \sigma)$;
\item the final output $\hat\sigma$ from Step 3 satisfies $\hat\sigma(i) = (\tau_1^{-1} \circ \tau_i)(\hat\sigma_i(i))$.
\end{enumerate}
\end{lemma}
\begin{proof}
(i) Denote $\epsilon = \max_i \dhams(\tilde\sigma_i, \sigma_{-i})$. Because the smallest block size of $\sigma_{-i}$ is bounded by $\frac12 \Nmin(\sigma_{-i}) \ge \frac12 (\Nmin - 1) > \epsilon$, it follows by Lemma~\ref{the:UniquePermutation} that for every $i$ there exists a unique $K$-permutation $\tau_i$ such that $\dham(\tau_i \circ \tilde\sigma_i, \sigma_{-i}) = \dhams(\tilde\sigma_i, \sigma_{-i})$.

(ii) Observe that 
\begin{equation}
 \label{eq:HalfBlockSize}
 \dham ( \tau_i \circ \hat\sigma_i, \sigma)
 \weq \dham ( \tau_i \circ \tilde\sigma_i, \sigma_{-i} ) + 1(\tau_i\circ\hat\sigma_i(i) = \sigma(i))
 \wle \epsilon + 1.
\end{equation}
Because $\epsilon + 1 < \frac12 \Nmin$, Lemma~\ref{the:UniquePermutation} implies that
$\dhams(\hat\sigma_i, \sigma) = \dham ( \tau_i \circ \hat\sigma_i, \sigma)$.
	
(iii) By Lemma~\ref{the:HammingInverseImages} and \eqref{eq:HalfBlockSize}, the minimum block size of $\hat\sigma_1$ is bounded by $\Nmin(\hat\sigma_1) \ge\Nmin - \dhams(\hat\sigma_1, \sigma) \ge \Nmin - (\epsilon+1)$. Inequality \eqref{eq:HalfBlockSize} also implies that
\begin{align*}
 \dham ( \tau_i \circ  \hat\sigma_i, \tau_1 \circ \hat\sigma_1 )
 \wle \dham ( \tau_i \circ \hat\sigma_i, \sigma ) + \dham ( \sigma, \tau_1 \circ \hat\sigma_1 )
 \wle 2(\epsilon+1).
\end{align*}
Therefore, $\dham ( \tau_1^{-1} \circ \tau_i \circ \hat\sigma_i, \hat\sigma_1 ) \le 2(\epsilon+1)$ as well. Furthermore, because $2(\epsilon+1) < \frac12 (\Nmin - (\epsilon+1)) \le \frac12 \Nmin(\hat\sigma_1)$, 
we conclude by Lemma~\ref{the:UniquePermutation} that $\tau_1^{-1} \circ \tau_i$ is the unique minimiser of $\tau \mapsto \dham(\tau \circ \hat\sigma_i, \hat\sigma_1)$, and
\[
\tau_1^{-1} \circ \tau_i(k)
\weq \argmax_\ell \abs{ \hat\sigma_i^{-1}(k) \cap \hat{\sigma}_1^{-1}( \ell ) }
\quad \text{for all $k$}.
\]
Hence, the output value $\hat\sigma(i)$ satisfies $\hat{\sigma}(i) = (\tau_1^{-1} \circ \tau_i)( \hat\sigma_i(i))$.
\end{proof}

\subsection{Proof of Theorem~\ref{thm:general_algo_consistency}}
\label{sec:general_algo_consistency_proof}

Let us first fix some measurable set $A \subset \cS$, and let $p=f(A)$ and $q = g(A)$. Denote $a = \frac{(p-q)^2}{p \vee q}$ and $b = (p^{1/2}-q^{1/2})^2$. A simple computation shows that $b \le a \le 4b$. Then it follows, by choosing suitable scale-dependent sets $A$, that $a \ge b \gg N^{-1} \beta_r$. In light of inequalities $(p-q)^2 \le (p \vee q)^2$
and $\dren_{1+r}^s(f,g) \ge \dren_{r}^s(f,g)$,
the assumption 
$\frac{(p-q)^2}{p \vee q}
\gg N^{-1} \frac{\dren_{1+r}^s(f,g)}{\dren_{r}^s(f,g)}$
implies that 
$p \vee q \ge \frac{(p-q)^2}{p \vee q} \gg N^{-1}$.

In Algorithm~\ref{algo:likelihood_based_algo}, the outputs of Step~1 are denoted by $\tilde \sigma_1, \dots, \tilde \sigma_N$, the outputs of Step~2 by $\hat \sigma_1, \dots, \hat \sigma_N$, and the final output from Step~3 by $\hat\sigma$.  Recall that $\sigma$ denotes the unknown true node labelling, and $\sigma_{-i}$ its restriction to $[N] \setminus \{i\}$.  As a standard graph clustering algorithm for Step~1, we will employ a spectral clustering algorithm described in \cite[Algorithm 4]{Xu_Jog_Loh_2020} with tuning parameter 
%
%
$\mu=8$ and trim threshold $\tau = 40 K \bar d$, where $\bar d$ is the average degree of $\tilde X_{-i}$, which is a modified version of \cite[Algorithm 2]{Gao_Ma_Zhang_Zhou_2017} with explicitly known error bounds.

Denote by $B$ the event that
$\frac{\Nmin}{N} \ge K^{-1} \Big( 1 - \sqrt{ \frac{ 8 K \log N}{N} } \Big)$.
Lemma~\ref{the:MultinomialConcentrationMinMaxDiff} shows that
\begin{align*}
 \pr \left( B^c \right) \wle K N^{-4}.
\end{align*}
Then $\frac{\Nmin}{N} \ge \frac12 K^{-1}$ on the event $B$, for large values of the scale parameter.
%
%
%
Denote by $E_i$ the event that Step~1 for node $i$ succeeds with accuracy $\dhams(\tilde\sigma_i, \sigma_{-i}) \le \epsilon N$, where 
$\epsilon = 2^{30} K^4 N^{-1} J^{-1}$ with $J = \frac{(p-q)^2}{p \vee q}$.
The matrix $\tilde X_{-i}$ computed in Step~1 of Algorithm~\ref{algo:likelihood_based_algo} is the adjacency matrix of a standard binary SBM with intra-block link probability $\pin$, inter-block link probability $\pout$, and node labelling $\sigma_{-i}$.
Because $p \vee q \ge N^{-1}$, and $N \ge 32K^2 \vee 2000$ and $J \ge 2^{36} K^6 N^{-1}$ for large values of the scale parameter,
by applying \cite[Proposition B.3]{Xu_Jog_Loh_2020} it
follows\footnote{
The statement of \cite[Proposition B.3]{Xu_Jog_Loh_2020} requires $p \wedge q \ge N^{-1}$ but the proof is valid also for $p \vee q \ge N^{-1}$.
}
that
\begin{align}
 \label{eq:in_proof_proba_initialisation_good}
 \pr\left( E_i^c \cap B \right) \wle N^{-5}.
\end{align}

The inequality $J \ge 2^{36} K^6 N^{-1}$ implies that $\epsilon N \le 2^{-4} K^{-1} \Nmin \le \frac{1}{32} \Nmin$, and the event $B$ implies $\Nmin \ge 135$. Therefore, $\epsilon N \le \frac{1}{32} \Nmin < \frac15 \Nmin - 1$,
and we see by applying Lemma~\ref{the:OptimalPermutations} that on the event $E \cap B$ where $E = \cap_i E_i$ there exist unique $K$-permutations $\tau_1,\dots,\tau_N$ such that
$\hat\sigma(i) = (\tau_1^{-1} \circ \tau_i)(\hat\sigma_i(i))$ for all~$i$. Especially,
\[
 \dham^*(\hat\sigma, \sigma)
 \wle \dham( \tau_1 \circ \hat\sigma, \sigma)
 \weq \sum_{i} 1( \tau_i\circ\hat\sigma_i(i) \ne \sigma(i) )
 \qquad \text{on $E \cap B$},
\]
so it follows that $\E \dham^*(\hat\sigma, \sigma) 1_{E\cap B} \le \sum_{i} \pr( \tau_i \circ \hat\sigma_i(i) \ne \sigma(i), E_i , B)$. In light of \eqref{eq:in_proof_proba_initialisation_good}, 
by using $\dhams = \left( 1_E 1_B + 1_{E^c} 1_B + 1_{B^C} \right) \dhams$ and applying the bounds $\dham^*(\hat\sigma, \sigma) \le N$ and $\pr(E^c \cap B ) \le \sum_i \pr(E_i^c \cap B ) \le N^{-4}$, we conclude that
\begin{equation}
 \label{eq:ProofAlgo1Bound1}
 \E \dham^*(\hat\sigma)
 \wle \sum_{i} \pr\left( \tau_i\circ\hat\sigma_i(i) \ne \sigma(i), E_i, B \right) + (K+1) N^{-3}.
\end{equation}

Let us analyse the sum on the right side of \eqref{eq:ProofAlgo1Bound1}.  Note that $\tilde\sigma_i$ and the $K$-permutation $\tau_i$ are fully determined by the entries of the sub-array $X_{-i} = (X_{j,j'}^t: j,j' \in [N] \setminus \{i\}, t \in [T] )$. Conditionally on $X_{-i}$, we may hence treat $\tilde\sigma_i$ and $\tau_i$ as nonrandom, and apply Lemma~\ref{lemma:individual_error_renyi} to conclude that on the event $E_i \cap B$,
\begin{align*}
 \pr \left( \tau_i \circ \hat\sigma_i(i) \ne \sigma(i) \cond X_{-i} \right)	
 &\wle
 K e^{-( N_{\min} - 1 - ( 2 + \frac{r}{1-r} \beta_r) \dhamsi ) 2 (1-r) \dren_r^s(\fin, \fout)} \\
 &\wle
 K e^{-\left( K^{-1} N \left( 1 - \sqrt{\frac{8 K \log N}{N}}\right) - 1 - ( 2 + \frac{r}{1-r} \beta_r) \epsilon N \right) 2 (1-r) \dren_r^s(\fin, \fout)},
\end{align*}
where the latter inequality is due to the definition of event $B$, and fact that
$\dhamsi = \dhams(\tilde\sigma_i, \sigma_{-i} ) \le \epsilon N$ on $E_i$.
%
%
Because the event $E_i \cap B$ is measurable with respect to the sigma-algebra generated by $X_{-i}$, we conclude that
the right side of the above inequality is also an upper bound for $\pr \left( \tau_i \circ \hat\sigma_i(i) \ne \sigma(i), E_i, B \right)$. We see by combining this with \eqref{eq:ProofAlgo1Bound1} that
\begin{align*}
 \E \frac{\dhams(\hsigma)}{N}
 \wle
 K e^{-\left( K^{-1} N \left( 1 - \sqrt{\frac{8 K \log N}{N}}\right) - 1 - ( 2 + \frac{r}{1-r} \beta_r) \epsilon N \right) 2 (1-r) \dren_r^s(\fin, \fout)}
 + (K+1) N^{-4}.
\end{align*}
The proof of Theorem~\ref{thm:general_algo_consistency} now follows by the inequality $(1-r) \dren_r \ge r \drenh$, valid for any $r \in \left(0, \frac12 \right]$ ~\cite[Theorem~16]{vanErven_Harremoes_2014},
and noting that the assumption $J \gg N^{-1} \beta_r$ implies that $\epsilon \beta_r \ll 1$.
\qed

\section{Information-theoretic divergences of sparse binary Markov chains}
\label{appendix:information_theoretic_divergences_sparse_markov_chains}

This section discusses binary Markov chains with initial distributions $\mu,\nu$ and transition probability matrices $P,Q$. In this case the \Renyi divergence of order $\alpha \in (0,\infty) \setminus \{1\}$ for the associated path probability distributions $f,g$ on $\{0,1\}^T$ equals
\begin{equation}
 \label{eq:MarkovRenyiGeneral}
 D_\alpha(f || g)
 \weq \frac{1}{\alpha-1}
 \log
 \bigg(
 \sum_{x \in \{0,1\}^T} \mu^\alpha_{x_1} \nu_{x_1}^{1-\alpha}
 \prod_{t=2}^T P_{x_{t-1} x_t}^\alpha Q_{x_{t-1} x_t}^{1-\alpha}
 \bigg).
\end{equation}
Such divergences will be analysed using weighted geometric and arithmetic averages of transition parameters defined by
\begin{equation}
 \label{eq:MarkovMean}
 \begin{aligned}
 r_a &= \mu_a^\alpha \nu_a^{1-\alpha},    \qquad\qquad  & R_{ab} &= P_{ab}^\alpha Q_{ab}^{1-\alpha}, \\
 \hr_a &= \alpha \mu_a + (1-\alpha) \nu_a,  & \qquad  \hR_{ab} &= \alpha P_{ab} + (1-\alpha) Q_{ab}.
 \end{aligned}
\end{equation}
We note that $D_\alpha(f || g) = \frac{1}{\alpha-1} \log Z$, where
$
 Z
 = \sum_{x \in \{0,1\}^T} r_{x_1}
 \prod_{t=2}^T R_{x_{t-1} x_t}. 
$
Moreover, $r_1 = 1-\hr_1+O(\rho^2)$ and $R_{01} = 1-\hR_{01} + O(\rho^2)$ when $\mu_1,\nu_1,P_{01},Q_{01} \lesim \rho$.

Section~\ref{sec:Proof_Thm_F1} presents the proof of Proposition~\ref{thm:renyi_divergence_order_less_1_sparse_markov_chains}. 
Section~\ref{sec:HighRenyiUpper} discusses high-order \Renyi divergences.

\subsection{Proof of Proposition~\ref{thm:renyi_divergence_order_less_1_sparse_markov_chains}}
\label{sec:Proof_Thm_F1}

Proposition~\ref{thm:renyi_divergence_order_less_1_sparse_markov_chains} follows by substituting $\alpha=\frac12$ in the following result.

\begin{proposition}
\label{the:RenyiSparseMCGeneral}
Consider binary Markov chains with initial distributions $\mu,\nu$ and
transition probability matrices $P,Q$. Assume that 
$\mu_1, \nu_1, P_{01}, Q_{01} \le \rho$ for some $\rho$ such that $\rho T \le 0.01$. Then the \Renyi divergence of order $\alpha \in (0,1)$ between the associated path probability distributions defined by~\eqref{eq:MarkovRenyiGeneral} is approximated by
\begin{equation}
 \label{eq:RenyiSparseMCGeneral}
 \dren_\alpha(f \| g)
 \weq \frac{1}{1-\alpha}
 \bigg( \hr_1 - r_1 + \sum_{t=2}^T J_t + \epsilon \bigg),
\end{equation}
where the error term satisfies
$\abs{\epsilon} \le 46 (\rho T)^2$,
\[
 J_t
 \weq
 \begin{cases}
 \hR_{01} - R_{01} +
 \Big(1 - \frac{R_{10}}{1 - R_{11}} \Big)
 \bigg( R_{01} + \Big( r_1 (1- R_{11}) - R_{01} \Big) R_{11}^{t-2} \bigg),
 &\quad R_{11} < 1, \\
 \hR_{01} - R_{01}, &\quad R_{11}=1,
 \end{cases}
\]
and the parameters $r_a,\hr_a,R_{ab},\hR_{ab}$ are given by \eqref{eq:MarkovMean}.
\end{proposition}

The rest of Section \ref{sec:Proof_Thm_F1} is devoted to proving Proposition~\ref{the:RenyiSparseMCGeneral}. 

\subsubsection{Basic results on binary sequences}

For a path $x = (x_1,\dots,x_T)$ in $\{0,1\}^T$, denote by $x_{ij} = \sum_{t=2}^T 1(x_{t-1}=i, x_t=j)$ the number of $ij$-transitions. Then, the path probability of a binary Markov chain with initial distribution $\mu$ and transition matrix $P$ can be written as $f(x) = \mu_{x_1} \prod_{ij} P_{ij}^{x_{ij}}$. For sparse Markov chains, we will analyse path probabilities by focusing on the total number of 1's
$\norm{x} = \sum_t x_t$,
and the number of on-periods
\[
\xon \weq x_1+x_{01} \weq x_{10}+x_T.
\]
The quantity $x_1+x_{01}$ counts the number of on-period start times, and the quantity $x_{10}+x_T$ counts the number of on-period end times. We also note that $x_{01} + x_{11} = \sum_{t=2}^T x_t$ implies that $\norm{x} = \xon + x_{11}$. The data $(\xon, \norm{x}, x_1, x_T)$ suffices to determine the path probability of $x$ because the transition counts can be recovered using the formulas $x_{01}  = \xon - x_1$, $x_{10} = \xon - x_T$, $x_{11} = \norm{x}-\xon$, together with $x_{00}+x_{01}+x_{10}+x_{11}=T-1$. Especially, the probability of a path with $(\xon, \norm{x}, x_1, x_T) = (j,t,a,b)$ equals
\begin{equation}
 \label{eq:MarkovPathOnPeriod}
 f(x)
 \weq \mu_0^{1-a} \mu_1^a P_{00}^{T-1-(t+j-a-b)}
 P_{01}^{j-a} P_{10}^{j-b} P_{11}^{t-j}.
\end{equation}

The number of such paths is summarised in the next result.

\begin{lemma}
\label{lemma:cardinals_transitions_binary_markov_chain}
Denote by $c_{jt}(ab)$ the number of paths $x \in \{0,1\}^T$ such that $\xon = j$, $\norm{x}=t$, $x_1=a$, and $x_T=b$.
%
Then the nonzero values of $c_{jt}(ab)$ are given by $c_{00}(00) = 1$,
\[ 
 c_{1t}(ab)
 \weq
 \begin{cases}
 T-t-1, &\quad (a,b)=(0,0), \ 1 \le t \le T-2, \\
 1,
 &\quad (a,b) = (0,1),(1,0), \ 1 \le t \le T-1, \\
 1, &\quad (a,b) = (1,1), \ t=T, 
\end{cases}
\]
and $c_{jt}(ab) = \binom{t-1}{j-1}\binom{T-t-1}{j-a-b}$ for $2 \le j \le \ceil{T/2}$, and $j \le t \le T-1-j+a+b$.
\end{lemma}

\begin{proof}
We compute the cardinalities separately for the three cases in which the number of on-periods equals $j=0$, $j=1$, and $j \ge 2$.
	
(i) Case $j=0$. The only path with no on-periods is the path of all zeros. Therefore, ${c_{0t}(ab)} = 1$ for $t=0$ and $(a,b) = (0,0)$, and ${c_{0t}(a,b)} = 0$ otherwise.
	
(ii) Case $j=1$. In this case ${c_{1t}(00)} = T-t-1$ for $1 \le t \le T-2$ and zero otherwise.
Furthermore, ${c_{1t}(01)} = {c_{1t}(10)} = 1$ for $1 \le t \le T-1$, and both are zero otherwise.
Finally, ${c_{1t}(11)} = 1$ for $t=T$ and zero otherwise.
	
(iii) Case $j \ge 2$. Now we proceed as follows. First, given a series of $t$ ones, we choose $j-1$ places to break the series: there are $\binom{t-1}{j-1}$ ways of doing so. Then, we need to fill those breaks with zeros chosen among the $T-t$ zeros of the chain. Note that when $a=b=0$, we also need to put zeros before and after the chain of ones. There are $j-1+(1-a)+(1-b) = j+1-a-b$ places to fill with $T-t$ zeros, and we need to put at least one zero in each place: there are $\binom{T-t-1}{j-a-b}$ ways of doing so.\footnote{A combinatorial fact, often referred as the \textit{stars and bars} method, is that the number of ways in which $n$ identical balls can be divided into $m$ distinct bins is
$
\begin{pmatrix}
n+m-1 \\ m-1
\end{pmatrix}
$, and
$\begin{pmatrix}
n-1 \\ m-1
\end{pmatrix}$ if bins cannot be empty.
}
Therefore, we conclude that
\[
c_{jt}(ab)
\weq \binom{t-1}{j-1}\binom{T-t-1}{j-a-b}.
\]
\end{proof}

\subsubsection{Useful Taylor expansions}


\begin{lemma}
\label{the:UsefulTaylor} 
Assume that $\alpha \in (0,1)$ and $\max\{\mu_1, \nu_1, P_{01}, Q_{01}\} \le \rho$ for some $\rho \le \frac13$. Then the geometric and arithmetic means defined by \eqref{eq:MarkovMean} are related according to
\begin{align*}
 r_0 &\weq 1 - \hr_1 + \epsilon_1, \\
 R_{00} &\weq 1 - \hat R_{01} + \epsilon_2, \\
 R_{00}^{T-1} &\weq 1 - (T-1) \hat R_{01} + \epsilon_3, \\
 r_0 R_{00}^{T-1} &\weq 1 - \hr_1 - (T-1) \hR_{01} + \epsilon_4,
\end{align*}
where the error terms are bounded by
$\abs{\epsilon_1}, \abs{\epsilon_2} \le (1+\rho)\rho^2$,
$\abs{\epsilon_3} \le 2(1+\rho) (\rho T)^2$, and $\abs{\epsilon_4} \le 4(1+2\rho) (\rho T)^2$. 
\end{lemma}

\begin{proof}
Note that $\hr_1 \le \rho$ and $\hat R_{01} \le \rho$.
Taylor's approximation (Lemma~\ref{the:Power}) implies that $(1-\mu_1)^\alpha = 1 - \alpha \mu_1 + \epsilon_{11}$ and $(1-\nu_1)^{1-\alpha} = 1 - (1-\alpha) \nu_1 + \epsilon_{12}$ for $\abs{\epsilon_{11}}, \abs{\epsilon_{12}} \le \frac12 \rho^2$. By multiplying these, we find that
\[
 r_0
 \weq (1-\mu_1)^\alpha(1-\nu_1)^{1-\alpha}
 \weq 1 - \hr_1 + \epsilon_1,
\]
where the error term is bounded by
$\abs{\epsilon_1}
\le (1+\frac14 \rho^2) \rho^2
\le (1+\rho) \rho^2$.
Because $R_{00} = (1-P_{01})^\alpha(1-Q_{01})^{1-\alpha}$, repeating the same argument yields $\abs{\epsilon_2} \le (1+\rho) \rho^2$.
	
Assume next that $T \ge 2$ (otherwise the third claim is trivial). Note that
$0 \le 1 - R_{00}
= \hat R_{01} - \epsilon_2
\le \rho + (1+\rho) \rho^2
\le \frac12$ due to
$\hR_{01} \le \rho$ and $\rho \le \frac13$.
By applying Lemma~\ref{the:Power}, we then see that
\[
 R_{00}^{T-1}
 \weq (1 - \hR_{01} + \epsilon_2 )^{T-1}
 \weq 1 - (T-1) (\hR_{01}-\epsilon_2) + \epsilon_{31},
\]
where
$\abs{\epsilon_{31}}
\le T^2 (\hR_{01}-\epsilon_2)^2
\le 2 T^2(\hR_{01}^2 + \epsilon_2^2)$.
It follows that
$R_{00}^{T-1} = 1 - (T-1) \hR_{01} + \epsilon_3$
with
$\epsilon_3 = (T-1) \epsilon_2 + \epsilon_{31}$
bounded by
$\abs{\epsilon_3}
\le T \abs{\epsilon_2} + \abs{\epsilon_{31}}
\le T \abs{\epsilon_2} + 2 (T \abs{\epsilon_2})^2 + 2 (\rho T)^2,
$
so that $\abs{\epsilon_3} \le 2(1+\rho) \rho^2 T^2$.
	
Finally, by multiplying the approximation formulas of $r_0$ and $R_{00}^{T-1}$, we find that
\begin{align*}
 \epsilon_4
 \weq \epsilon_1 ( 1 - (T-1) \hR_{01} )
 + \epsilon_3 (1 - \hr_1 )
 + \epsilon_1 \epsilon_3
 + (T-1) \hr_1 \hR_{01}.
\end{align*}
By the triangle inequality, we find that for $T \ge 2$, 
$\abs{\epsilon_4}
\le (1+\rho T) \abs{\epsilon_1} + \abs{\epsilon_3} + \abs{\epsilon_1 \epsilon_3} + \rho^2 T$,
from which one may check that $\abs{\epsilon_4} \le 4(1+2\rho)(\rho T)^2$.
\end{proof}

\subsubsection{Analysing paths with two or more on-periods}

\begin{lemma}
\label{lemma:FirstOrderRenyi_more_than_two_onperiods}
For any $\alpha \in (0,1)$ and any Markov chain path distributions $f,g$ with transition matrices $P,Q$ satisfying $P_{11} Q_{11} < 1$,
\[
 \sum_{x:\xon \ge 2} f_x^\alpha g_x^{1-\alpha}
 \wle R_{01}(r_1 + R_{01}) T^2 \Zgeo (1+\Zgeo) e^{\Zgeo R_{01}T},
\]
where $\Zgeo = \frac{R_{10}}{1-R_{11}}$ and the weighted geometric means $r_a,R_{ab}$ are defined by \eqref{eq:MarkovMean}.
\end{lemma}
\begin{proof}
Fix an integer $2 \le j \le \ceil{T/2}$, and denote $Z_j = \sum_{x:\xon=j} f_x^\alpha g_x^{1-\alpha}$. By \eqref{eq:MarkovPathOnPeriod}, we see that for any path $x$ with $j$ on-periods, $t$ ones, initial state $a$, and final state $b$,
\[
 f_x^\alpha g_x^{1-\alpha}
 \weq r_0^{1-a}r_1^a
 R_{00}^{T-1-(t+j-a-b)} R_{01}^{j-a} R_{10}^{j-b} R_{11}^{t-j}
 \wle r_1^a R_{01}^{j-a} R_{10}^{j-b} R_{11}^{t-j}.
\]
By Lemma~\ref{lemma:cardinals_transitions_binary_markov_chain}, the number of such paths equals
\[
 c_{jt}(ab)
 \weq \binom{t-1}{j-1} \binom{T-t-1}{j-a-b} .
\]
To obtain an upper bound for the path count, we note that
$\binom{T-t-1}{j-a-b} \le \frac{T^{j-a-b}}{(j-a-b)!}$. Furthermore, we also see that
$\binom{t-1}{j-1} = \frac{t-1}{j-1} \binom{t-2}{j-2} \le T \binom{t-2}{j-2}$.
The latter bound implies that $\binom{t-1}{j-1} \le T^b \binom{t-b-1}{j-b-1}$ for all $b \in \{0,1\}$. As a consequence, we conclude that
\[
 c_{jt}(ab)
 \wle \frac{T^{j}}{(j-2)!} \binom{t-b-1}{j-b-1}
\]
holds for all $a,b \in \{0,1\}$. 
Hence,
\begin{align*}
 Z_j
 &\wle \frac{T^{j}}{(j-2)!}
 \sum_{a,b=0}^1 \sum_{t \ge j}
 \binom{t-b-1}{j-b-1}
 r_1^a R_{01}^{j-a} R_{10}^{j-b} R_{11}^{t-j}.
\end{align*}
Using a geometric moment formula (Lemma~\ref{the:GeometricMoments}), we find that
\[
 \sum_{t=j}^\infty \binom{t-b-1}{j-b-1} R_{11}^{t-j}
 \weq (1-R_{11})^{-(j-b)}
 \weq R_{10}^{b-j} \Zgeo^{j-b},
\]
and it follows that
\begin{align*}
 Z_j
 \wle \frac{T^{j}}{(j-2)!}
 \sum_{a,b=0}^1 r_1^a R_{01}^{j-a} \Zgeo^{j-b}
 \weq T^2 \frac{(R_{01} \Zgeo T)^{j-2}}{(j-2)!}
 \sum_{a,b=0}^1 r_1^a R_{01}^{2-a} \Zgeo^{2-b}.
\end{align*}
By noting that $\sum_{a,b=0}^1 r_1^a R_{01}^{2-a} \Zgeo^{2-b} = R_{01}(r_1 + R_{01}) \Zgeo (1+\Zgeo)$ and summing the above inequality with respect to $j \ge 2$, the claim follows.
\end{proof}

\subsubsection{Proof of Proposition~\ref{the:RenyiSparseMCGeneral}}

\begin{proof}
By definition,
$\dren_\alpha(f\|g) = \frac{1}{\alpha-1} \log Z$, where $Z = \sum_x f_x^\alpha g_x^{1-\alpha}$. We may split the latter sum as
\begin{align*}
 Z
 \weq Z_0 + Z_1 + \sum_{j=2}^{ \ceil{T/2} } Z_j,
\end{align*}
where $Z_j = \sum_{x:\xon = j} f_x^{\alpha} g_x^{1-\alpha}$ indicates a Hellinger sum over paths with $j$ on-periods. We will approximate the first two terms on the right by $Z_0 = \hZ_0 + \epsilon_0$, $Z_1 = \hZ_1 + \epsilon_1$, where
\[
 \hZ_0 \weq 1 - \hr_1 - (T-1)\hR_{01}
\]
and
\[
 \hZ_1
 \weq R_{01} R_{10}
 \sum_{t=1}^{T-2}(T-t-1) R_{11}^{t-1} 
 + \left( R_{01} + r_1 R_{10} \right)
 \sum_{t=1}^{T-1} R_{11}^{t-1} 
 + r_1 R_{11}^{T-1}.
\]
Then, it follows that
\begin{equation}
 \label{eq:ApproximationZ}
 Z
 \weq \hZ_0 + \hZ_1 + \epsilon_0 + \epsilon_1 + \epsilon_2,
\end{equation}
where $\epsilon_2 = \sum_{j=2}^{ \ceil{T/2} } Z_j$. 
When $R_{11} < 1$, by applying formulas $\sum_{t=1}^{T-2} (T-t-1) R_{11}^{t-1} = (1-R_{11})^{-1} \left( (T-1) - \sum_{t=2}^T R_{11}^{t-2} \right)$ and
$R_{11}^{T-1} = 1 - (1-R_{11}) \sum_{t=2}^T R_{11}^{t-2}$ we find that
\[
 \hZ_1
 \weq r_1 + \Zgeo R_{01} (T-1) - (1 - \Zgeo)
 \left( r_1 (1-R_{11}) - R_{01} \right)
 \sum_{t=2}^T R_{11}^{t-2},
\]
where $\Zgeo = \frac{R_{10}}{1-R_{11}}$.
Hence,
\begin{equation}
 \label{eq:HatZ01}
 \hZ_0 + \hZ_1
 \weq 1 -
 \left( \hr_1 - r_1 + \sum_{t=2}^T J_t \right),
\end{equation}
where the expression of $J_t$ coincides with the one in the statement of the proposition. When $R_{11}=1$, we find that $J_t = \hat R_{01} - R_{01}$.


Let us next derive upper bounds for the error terms in \eqref{eq:ApproximationZ}.  We start with $\epsilon_0$. Because the only path with $\xon = 0$ is the identically zero path, we find that $Z_0 = r_0 R_{00}^{T-1}$. By Lemma~\ref{the:UsefulTaylor} we have
$\abs{\epsilon_0} \le 4(1+2\rho)(\rho T)^2
\le 5 (\rho T)^2$.

For the error term $\epsilon_1$,
with the help of formula \eqref{eq:MarkovPathOnPeriod} and Lemma~\ref{lemma:cardinals_transitions_binary_markov_chain}, we see that
\begin{align*}
 Z_1
 &\weq r_0 R_{01} R_{10} \sum_{t=1}^{T-2}(T-t-1) R_{11}^{t-1} R_{00}^{T-2-t} 
+ \left( r_0 R_{01} + r_1 R_{10}  \right) \sum_{t=1}^{T-1} R_{11}^{t-1} R_{00}^{T-1-t} \\ 
 &\qquad\quad + r_1 R_{11}^{T-1}.
\end{align*}
Because $r_0, R_{00} \le 1$, it follows that $Z_1 \le \hZ_1$, and hence $\epsilon_1 \le 0$. Furthermore, Lemma~\ref{the:UsefulTaylor} implies that
$r_0, R_{00} \ge 1 - 2\rho$. By noting that $R_{00}^{T-t} \ge R_{00}^{T-1}$ for $t \ge 1$, it follows that
\[
 Z_1
 \wge (1 - 2 \rho )^T \hZ_1
 \wge (1 - 2 \rho T) \hZ_1.
\]
For $R_{11}=1$ we have $R_{10}=0$ and $\hZ_1 = r_1 + (T-1)R_{01}$.
For $R_{11} < 1$, we observe that
\begin{align*}
 \hZ_1
 &\wle (T-1) R_{01} \frac{R_{10}}{1-R_{11}}
 + (T-1) R_{01}
 + r_1 \frac{R_{10}}{1 - R_{11}} + r_1 \\
 &\wle (1+\Zgeo) \Big(r_1 + (T-1) R_{01} \Big),
\end{align*}
where $\Zgeo = \frac{R_{10}}{1-R_{11}}$. We note that $\Zgeo = Z_\alpha( \Geo(P_{11}) \| \Geo(Q_{11}))$ equals the Hellinger sum of two geometric distributions, and therefore, $\Zgeo \in (0,1]$. Hence, $\hZ_1 \le 2\rho T$, and it follows that
\[
 Z_1
 \wge \hZ_1 - 2 \rho T \hZ_1
 \wge \hZ_1 - 4 (\rho T)^2.
\]
Thus, $\abs{\epsilon_1} \le 4 (\rho T)^2$ for both $R_{11}<1$ and $R_{11}=1$.

For the last error term in \eqref{eq:ApproximationZ},
we see that $\epsilon_2=0$ for $R_{11}=1$, whereas for $R_{11}<1$, Lemma~\ref{lemma:FirstOrderRenyi_more_than_two_onperiods} shows that
$0 \le \epsilon_2 \le 4 (\rho T)^2 e^{\rho T}
\le 5 (\rho T)^2$. By combining the error bounds for $\epsilon_0, \epsilon_1, \epsilon_2$, we may now conclude that
\begin{equation}
 \label{eq:ApproximationZPrime}
 Z
 \weq \hZ_0 + \hZ_1 + \epsilon',
\end{equation}
where
$
 \abs{\epsilon'}
 \le 14 (\rho T)^2.
$

Finally, Taylor's approximation (Lemma~\ref{the:Log}) shows that $\log(1-t) = -t - \epsilon''$ where $0 \le \epsilon'' \le 2t^2$ for $0 \le t \le \frac12$. By applying this with $t = 1-Z$, and noting that $\abs{J_t} \le \rho$ implies $0 \le t \le 3 \rho T + \abs{\epsilon'} \le 4 \rho T$,
we find that
\begin{align*}
 D_\alpha(f\|g)
 \weq \frac{1}{1-\alpha}
 \left( 1 - Z + \epsilon'' \right)
 \weq \frac{1}{1-\alpha}
 \left( 1-\hZ_0 - \hZ_1 - \epsilon' + \epsilon'' \right).
\end{align*}
The error bound of formula \eqref{eq:RenyiSparseMCGeneral} now follows from \eqref{eq:HatZ01} after noting that
\[
 \abs{\epsilon'} + \abs{\epsilon''}
 \wle 14 (\rho T)^2 + 2 ( 4 \rho T)^2
 \wle 46 (\rho T)^2.
\]

%
\end{proof}

\subsection{High-order \Renyi divergences}
\label{sec:HighRenyiUpper}

The following result provides an upper bound on the \Renyi divergence of order $\alpha > 1$ between path probability distributions of binary Markov chains defined by \eqref{eq:MarkovRenyiGeneral}.

\begin{proposition}
\label{prop:renyi_divergence_bound_order_more_than_one_markov_chains}
Assume that $\frac{\mu_1}{\nu_1}, \frac{P_{01}}{Q_{01}}, \frac{P_{10}}{Q_{10}} \le M$ for some $M \ge 1$, $Q_{11}>0$, and $\nu_1, Q_{01} \le \rho$ for some $\rho \le \frac12$. Then the \Renyi divergence of order $1 < \alpha < \infty $ is bounded by
\begin{equation}
 \label{eq:MarkovRenyiLargeBound0}
 \begin{aligned}
 D_\alpha(f || g)
 &\wle \frac{2\alpha}{\alpha-1} \rho T
 + \frac{M^{2\alpha}}{\alpha-1} \rho T \sum_{t=0}^{T-1} \Lambda^t \\
 &\qquad + \frac{4}{\alpha-1}
 \sum_{j=2}^{\lceil{T/2} \rceil} \frac{(M^{2\alpha} \rho T )^j}{(j-2)!}
 \sum_{t=j}^T \binom{t-1}{j-1} \Lambda^{t-j},
 \end{aligned}
\end{equation}
where $\Lambda = P_{11}^\alpha Q_{11}^{1-\alpha}$. 
Furthermore, when $\Lambda < 1$,
\begin{equation}
 \label{eq:MarkovRenyiLargeBound01}
 D_\alpha(f || g)
 \wle \frac{2 \alpha +1}{\alpha-1} C \rho T e^{5 C \rho T}
 \qquad \text{with $C = \frac{M^{2\alpha}}{1-\Lambda}$}.
\end{equation}
\end{proposition}
\begin{proof}
Recall that $D_\alpha(f\|g) = \frac{1}{\alpha-1} \log Z$ where
$Z = \sum_x g_x (f_x/g_x)^\alpha$. Because $\nu_1 \le \rho$ with $\rho \le \frac12$, we find that $\frac{\mu_0}{\nu_0} \le \frac{1}{1-\nu_1} = 1 + \frac{\nu_1}{1-\nu_1} \le 1 + 2 \rho$. Because $Q_{01} \le \rho$, the same argument shows that $\frac{P_{00}}{Q_{00}} \le 1 + 2 \rho$. Because $1-x_1 + x_{00} \le T$, it follows that
\begin{align*}
 \frac{f_x}{g_x}
 &\weq \left( \frac{\mu_0}{\nu_0} \right)^{1-x_1}
 \left( \frac{\mu_1}{\nu_1} \right)^{x_1}
 \left( \frac{P_{00}}{Q_{00}} \right)^{x_{00}}
 \left( \frac{P_{01}}{Q_{01}} \right)^{x_{01}}
 \left( \frac{P_{10}}{Q_{10}} \right)^{x_{10}}
 \left( \frac{P_{11}}{Q_{11}} \right)^{x_{11}} \\
 &\wle (1+2\rho)^T
 M^{x_1 + x_{01} + x_{10}}
 \left( \frac{P_{11}}{Q_{11}} \right)^{x_{11}}.
\end{align*}
Observe also that $g_x \le \nu_1^{x_1} Q_{01}^{x_{01}} Q_{11}^{x_{11}}
\le \rho^{x_1 + x_{01}} Q_{11}^{x_{11}}$. Therefore,
\[
 Z
 \wle (1+2\rho)^{\alpha T}
 \sum_x \rho^{x_1 + x_{01}} M^{\alpha (x_1 + x_{01} + x_{10})}
 \Lambda^{x_{11}},
\]
where $\Lambda = P_{11}^\alpha Q_{11}^{1-\alpha}$. By recalling that $\xon = x_1 + x_{01} = x_{10} + x_T$ and $\norm{x} = x_1 + x_{01} + x_{11} = x_{10} + x_{11} + x_T$, we find that
$x_1 + x_{01} + x_{10} = 2 \xon - x_T \le 2 \xon$ and $x_{11} = \norm{x} - \xon$.
Hence
\begin{equation}
 \label{eq:MarkovRenyiLargeBound1}
 Z
 \wle (1+2\rho)^{\alpha T}
 \sum_{x} \rho^{\xon} M^{2 \alpha \xon}
 \Lambda^{\norm{x} - \xon}
 \weq (1+2\rho)^{\alpha T}
 \sum_{j=0}^{\lceil{T/2} \rceil}  S_j
\end{equation}
where
\[
 S_j
 \weq (M^{2\alpha} \rho)^j  \sum_{t=j}^T c_{jt} \Lambda^{t-j},
\]
and $c_{jt}$ is the number of paths $x \in \{0,1\}^T$ containing $\xon = j$ on-periods and $\norm{x} = t$ ones. Because there is only one path containing no ones, and this path has no on-periods, we find that $S_0=1$. By noting that $\log(1+t) \le t$, it follows from \eqref{eq:MarkovRenyiLargeBound1} that
\begin{equation}
 \label{eq:MarkovRenyiLargeBound2}
 D_\alpha(f || g)
 \wle \frac{2\alpha}{\alpha-1} \rho T + \frac{1}{\alpha-1}  \sum_{j=1}^{\lceil{T/2} \rceil}  S_j.
\end{equation}
Because $c_{1t} \le T$ for all $t$, we se that
\begin{equation}
 \label{eq:S1Bound}
 S_1
 \wle M^{2\alpha} \rho T \sum_{t=0}^{T-1} \Lambda^t.
\end{equation}
For $j \ge 2$, Lemma~\ref{lemma:cardinals_transitions_binary_markov_chain} implies that
$c_{jt} = \sum_{a,b=0}^1 \binom{t-1}{j-1} \binom{T-t-1}{j-a-b} \le 4 \frac{T^j}{(j-2)!} \binom{t-1}{j-1}$, and we find that
\begin{equation}
 \label{eq:S2Bound}
 S_j
 \wle 4 \frac{(M^{2\alpha} \rho T )^j}{(j-2)!}
 \sum_{t=j}^T \binom{t-1}{j-1} \Lambda^{t-j}.
\end{equation}
Inequality \eqref{eq:MarkovRenyiLargeBound0} follows by substituting \eqref{eq:S1Bound}--\eqref{eq:S2Bound} into \eqref{eq:MarkovRenyiLargeBound2}.

Assume next that $\Lambda < 1$, and denote $C = \frac{M^{2\alpha}}{1-\Lambda}$. By replacing $T-1$ by infinity  on the right side of \eqref{eq:S1Bound}, it follows that $S_1 \le C \rho T$. By a geometric moment formula (Lemma~\ref{the:GeometricMoments}), we find that
\[
 \sum_{t=j}^T \binom{t-1}{j-1} \Lambda^{t-j}
 \wle \sum_{t=j}^\infty \binom{t-1}{j-1} \Lambda^{t-j}
 \weq (1-\Lambda)^{-j}.
\]
Then \eqref{eq:S2Bound} implies that
\begin{align*}
 \sum_{j=2}^{\lceil{T/2} \rceil} S_j
 \wle 4\sum_{j=2}^{\lceil{T/2} \rceil} \frac{(C\rho T)^j}{(j-2)!}
 \wle 4\sum_{j=2}^\infty \frac{(C\rho T)^j}{(j-2)!}
 \weq 4 (C \rho T)^2 e^{C \rho T}.
\end{align*}
Now it follows by \eqref{eq:MarkovRenyiLargeBound2} that
\[
 D_\alpha(f || g)
 \wle \frac{2\alpha}{\alpha-1} \rho T + \frac{C\rho T}{\alpha-1}
 + \frac{4 (C \rho T)^2}{\alpha-1} e^{C \rho T}.
\]
Therefore,
\[
 \frac{(\alpha-1) D_\alpha(f || g)}{C\rho T}
 \wle \frac{2\alpha}{C}+1+4C\rho T e^{C\rho T} 
 \wle \left( \frac{2\alpha}{C}+1+4C\rho T \right) e^{C\rho T}.
\]
Because 
$
 \frac{2\alpha}{C}+1+4C\rho 
 \le \left( \frac{2\alpha}{C}+1 \right) (1 + 4C\rho T)
 \le \left( \frac{2\alpha}{C}+1 \right) e^{4C\rho T},
$
we conclude that 
\[
 \frac{(\alpha-1) D_\alpha(f || g)}{C\rho T}
 \wle \left( \frac{2\alpha}{C}+1 \right) e^{5C\rho T}
\]
Because $C \ge 1$, we see that $\frac{2\alpha}{C}+1 \le 2\alpha+1$,
and \eqref{eq:MarkovRenyiLargeBound01} follows.
\end{proof}

%% file: main.bbl
\newcommand{\SortNoop}[1]{}\def\cprime{$'$}
\begin{thebibliography}{59}

\bibitem{Abbe_2018_JMLR}
\begin{barticle}[author]
\bauthor{\bsnm{Abbe},~\bfnm{Emmanuel}\binits{E.}}
(\byear{2018}).
\btitle{Community detection and stochastic block models: {Recent}
  developments}.
\bjournal{Journal of Machine Learning Research}
\bvolume{18}
\bpages{1--86}.
\end{barticle}
\endbibitem

\bibitem{Abbe_Bandeira_Hall_2016}
\begin{barticle}[author]
\bauthor{\bsnm{Abbe},~\bfnm{Emmanuel}\binits{E.}},
  \bauthor{\bsnm{Bandeira},~\bfnm{Afonso~S.}\binits{A.~S.}} \AND
  \bauthor{\bsnm{Hall},~\bfnm{Georgina}\binits{G.}}
(\byear{2016}).
\btitle{Exact recovery in the stochastic block model}.
\bjournal{IEEE Transactions on Information Theory}
\bvolume{62}
\bpages{471--487}.
\end{barticle}
\endbibitem

\bibitem{Alaluusua_Leskela_2022}
\begin{binproceedings}[author]
\bauthor{\bsnm{Alaluusua},~\bfnm{Kalle}\binits{K.}} \AND
  \bauthor{\bsnm{{Leskel{\"a}}},~\bfnm{Lasse}\binits{L.}}
(\byear{2022}).
\btitle{Consistent Bayesian community recovery in multilayer networks}.
In \bbooktitle{2022 IEEE International Symposium on Information Theory (ISIT)}.
\bnote{\url{https://arxiv.org/abs/2202.05823}}.
\end{binproceedings}
\endbibitem

\bibitem{Avrachenkov_Dreveton_Leskela_2021}
\begin{binproceedings}[author]
\bauthor{\bsnm{Avrachenkov},~\bfnm{Konstantin}\binits{K.}},
  \bauthor{\bsnm{Dreveton},~\bfnm{Maximilien}\binits{M.}} \AND
  \bauthor{\bsnm{Leskel{\"a}},~\bfnm{Lasse}\binits{L.}}
(\byear{2021}).
\btitle{Recovering communities in temporal networks using persistent edges}.
In \bbooktitle{Computational Data and Social Networks (CSoNet 2021)}
(\beditor{\bfnm{David}\binits{D.}~\bsnm{Mohaisen}} \AND
  \beditor{\bfnm{Ruoming}\binits{R.}~\bsnm{Jin}}, eds.).
\bseries{Lecture Notes in Computer Science}
\bvolume{13116}
\bpages{243--254}.
\bpublisher{Springer}.
\end{binproceedings}
\endbibitem

\bibitem{Barucca_Lillo_Mazzarisi_Tantari_2018}
\begin{barticle}[author]
\bauthor{\bsnm{Barucca},~\bfnm{Paolo}\binits{P.}},
  \bauthor{\bsnm{Lillo},~\bfnm{Fabrizio}\binits{F.}},
  \bauthor{\bsnm{Mazzarisi},~\bfnm{Piero}\binits{P.}} \AND
  \bauthor{\bsnm{Tantari},~\bfnm{Daniele}\binits{D.}}
(\byear{2018}).
\btitle{Disentangling group and link persistence in dynamic stochastic block
  models}.
\bjournal{Journal of Statistical Mechanics: Theory and Experiment}
\bvolume{2018}
\bpages{1--18}.
\bdoi{10.1088/1742-5468/aaeb44}
\end{barticle}
\endbibitem

\bibitem{Bassett_Wymbs_Porter_Mucha_Carlson_Grafton_2011}
\begin{barticle}[author]
\bauthor{\bsnm{Bassett},~\bfnm{Danielle~S}\binits{D.~S.}},
  \bauthor{\bsnm{Wymbs},~\bfnm{Nicholas~F}\binits{N.~F.}},
  \bauthor{\bsnm{Porter},~\bfnm{Mason~A}\binits{M.~A.}},
  \bauthor{\bsnm{Mucha},~\bfnm{Peter~J}\binits{P.~J.}},
  \bauthor{\bsnm{Carlson},~\bfnm{Jean~M}\binits{J.~M.}} \AND
  \bauthor{\bsnm{Grafton},~\bfnm{Scott~T}\binits{S.~T.}}
(\byear{2011}).
\btitle{Dynamic reconfiguration of human brain networks during learning}.
\bjournal{Proceedings of the National Academy of Sciences}
\bvolume{108}
\bpages{7641--7646}.
\end{barticle}
\endbibitem

\bibitem{Bhattacharyya_Chatterjee_2018-05-27}
\begin{bmisc}[author]
\bauthor{\bsnm{Bhattacharyya},~\bfnm{Sharmodeep}\binits{S.}} \AND
  \bauthor{\bsnm{Chatterjee},~\bfnm{Shirshendu}\binits{S.}}
(\byear{2018}).
\btitle{Spectral clustering for multiple sparse networks: {I}}.
\bnote{\url{https://arxiv.org/abs/1805.10594}}.
\end{bmisc}
\endbibitem

\bibitem{Bhattacharyya_Chatterjee_2020-04-06}
\begin{bmisc}[author]
\bauthor{\bsnm{Bhattacharyya},~\bfnm{Sharmodeep}\binits{S.}} \AND
  \bauthor{\bsnm{Chatterjee},~\bfnm{Shirshendu}\binits{S.}}
(\byear{2020}).
\btitle{General community detection with optimal recovery conditions for
  multi-relational sparse networks with dependent layers}.
\bnote{\url{https://arxiv.org/abs/2004.03480}}.
\end{bmisc}
\endbibitem

\bibitem{Bickel_Chen_2009}
\begin{barticle}[author]
\bauthor{\bsnm{Bickel},~\bfnm{Peter~J}\binits{P.~J.}} \AND
  \bauthor{\bsnm{Chen},~\bfnm{Aiyou}\binits{A.}}
(\byear{2009}).
\btitle{A nonparametric view of network models and Newman--Girvan and other
  modularities}.
\bjournal{Proceedings of the National Academy of Sciences}
\bvolume{106}
\bpages{21068--21073}.
\end{barticle}
\endbibitem

\bibitem{Billingsley_1961}
\begin{barticle}[author]
\bauthor{\bsnm{Billingsley},~\bfnm{Patrick}\binits{P.}}
(\byear{1961}).
\btitle{Statistical methods in {M}arkov chains}.
\bjournal{Annals of Mathematical Statistics}
\bvolume{32}
\bpages{12--40}.
\bdoi{10.1214/aoms/1177705136}
\end{barticle}
\endbibitem

\bibitem{Dhara_Souvik_Mossel_Sandon_2022}
\begin{binproceedings}[author]
\bauthor{\bsnm{Dhara},~\bfnm{Souvik}\binits{S.}},
  \bauthor{\bsnm{Gaudio},~\bfnm{Julia}\binits{J.}},
  \bauthor{\bsnm{Mossel},~\bfnm{Elchanan}\binits{E.}} \AND
  \bauthor{\bsnm{Sandon},~\bfnm{Colin}\binits{C.}}
(\byear{2022}).
\btitle{Spectral recovery of binary censored block models}.
In \bbooktitle{Proceedings of the 33rd ACM-SIAM Symposium on Discrete
  Algorithms}
\bpages{3389--3416}.
\end{binproceedings}
\endbibitem

\bibitem{Fortunato_2010}
\begin{barticle}[author]
\bauthor{\bsnm{Fortunato},~\bfnm{Santo}\binits{S.}}
(\byear{2010}).
\btitle{Community detection in graphs}.
\bjournal{Physics Reports}
\bvolume{486}
\bpages{75--174}.
\bdoi{http://dx.doi.org/10.1016/j.physrep.2009.11.002}
\end{barticle}
\endbibitem

\bibitem{Fournet_Barrat_2014}
\begin{barticle}[author]
\bauthor{\bsnm{Fournet},~\bfnm{Julie}\binits{J.}} \AND
  \bauthor{\bsnm{Barrat},~\bfnm{Alain}\binits{A.}}
(\byear{2014}).
\btitle{Contact patterns among high school students}.
\bjournal{PLOS ONE}
\bvolume{9}
\bpages{1-17}.
\end{barticle}
\endbibitem

\bibitem{Gao_Ma_Zhang_Zhou_2017}
\begin{barticle}[author]
\bauthor{\bsnm{Gao},~\bfnm{Chao}\binits{C.}},
  \bauthor{\bsnm{Ma},~\bfnm{Zongming}\binits{Z.}},
  \bauthor{\bsnm{Zhang},~\bfnm{Anderson~Y.}\binits{A.~Y.}} \AND
  \bauthor{\bsnm{Zhou},~\bfnm{Harrison~H.}\binits{H.~H.}}
(\byear{2017}).
\btitle{Achieving optimal misclassification proportion in stochastic block
  models}.
\bjournal{Journal of Machine Learning Research}
\bvolume{18}
\bpages{1980--2024}.
\end{barticle}
\endbibitem

\bibitem{Ghasemian_Zhang_Clauset_Moore_Peel_2016}
\begin{barticle}[author]
\bauthor{\bsnm{Ghasemian},~\bfnm{Amir}\binits{A.}},
  \bauthor{\bsnm{Zhang},~\bfnm{Pan}\binits{P.}},
  \bauthor{\bsnm{Clauset},~\bfnm{Aaron}\binits{A.}},
  \bauthor{\bsnm{Moore},~\bfnm{Cristopher}\binits{C.}} \AND
  \bauthor{\bsnm{Peel},~\bfnm{Leto}\binits{L.}}
(\byear{2016}).
\btitle{Detectability thresholds and optimal algorithms for community structure
  in dynamic networks}.
\bjournal{Physical Review X}
\bvolume{6}
\bpages{031005}.
\end{barticle}
\endbibitem

\bibitem{Ghosal_VanDerVaart_2017}
\begin{bbook}[author]
\bauthor{\bsnm{Ghosal},~\bfnm{Subhashis}\binits{S.}} \AND
  \bauthor{\bparticle{Van~der} \bsnm{Vaart},~\bfnm{Aad}\binits{A.}}
(\byear{2017}).
\btitle{Fundamentals of nonparametric Bayesian inference}
\bvolume{44}.
\bpublisher{Cambridge University Press}.
\end{bbook}
\endbibitem

\bibitem{Gosgens_Tikhonov_Prokhorenkova_2021}
\begin{binproceedings}[author]
\bauthor{\bsnm{G{\"o}sgens},~\bfnm{Martijn~M}\binits{M.~M.}},
  \bauthor{\bsnm{Tikhonov},~\bfnm{Alexey}\binits{A.}} \AND
  \bauthor{\bsnm{Prokhorenkova},~\bfnm{Liudmila}\binits{L.}}
(\byear{2021}).
\btitle{Systematic analysis of cluster similarity indices: {How} to validate
  validation measures}.
In \bbooktitle{Proceedings of the 38th International Conference on Machine
  Learning}
(\beditor{\bfnm{Marina}\binits{M.}~\bsnm{Meila}} \AND
  \beditor{\bfnm{Tong}\binits{T.}~\bsnm{Zhang}}, eds.)
\bvolume{139}
\bpages{3799--3808}.
\end{binproceedings}
\endbibitem

\bibitem{Hajek_Wu_Xu_2016}
\begin{barticle}[author]
\bauthor{\bsnm{Hajek},~\bfnm{Bruce}\binits{B.}},
  \bauthor{\bsnm{Wu},~\bfnm{Yihong}\binits{Y.}} \AND
  \bauthor{\bsnm{Xu},~\bfnm{Jiaming}\binits{J.}}
(\byear{2016}).
\btitle{Achieving exact cluster recovery threshold via semidefinite
  programming: {E}xtensions}.
\bjournal{IEEE Transactions on Information Theory}
\bvolume{62}
\bpages{5918-5937}.
\bdoi{10.1109/TIT.2016.2594812}
\end{barticle}
\endbibitem

\bibitem{Han_Xu_Airoldi_2015}
\begin{binproceedings}[author]
\bauthor{\bsnm{Han},~\bfnm{Qiuyi}\binits{Q.}},
  \bauthor{\bsnm{Xu},~\bfnm{Kevin}\binits{K.}} \AND
  \bauthor{\bsnm{Airoldi},~\bfnm{Edoardo}\binits{E.}}
(\byear{2015}).
\btitle{Consistent estimation of dynamic and multi-layer block models}.
In \bbooktitle{Proceedings of the 32nd International Conference on Machine
  Learning}
(\beditor{\bfnm{Francis}\binits{F.}~\bsnm{Bach}} \AND
  \beditor{\bfnm{David}\binits{D.}~\bsnm{Blei}}, eds.)
\bvolume{37}
\bpages{1511--1520}.
\end{binproceedings}
\endbibitem

\bibitem{Hartle_Papadopoulos_Krioukov_2021}
\begin{barticle}[author]
\bauthor{\bsnm{Hartle},~\bfnm{Harrison}\binits{H.}},
  \bauthor{\bsnm{Papadopoulos},~\bfnm{Fragkiskos}\binits{F.}} \AND
  \bauthor{\bsnm{Krioukov},~\bfnm{Dmitri}\binits{D.}}
(\byear{2021}).
\btitle{Dynamic hidden-variable network models}.
\bjournal{Phys. Rev. E}
\bvolume{103}
\bpages{052307}.
\bdoi{10.1103/PhysRevE.103.052307}
\end{barticle}
\endbibitem

\bibitem{Heimlicher_Lelarge_Massoulie_2012}
\begin{binproceedings}[author]
\bauthor{\bsnm{Heimlicher},~\bfnm{Simon}\binits{S.}},
  \bauthor{\bsnm{Lelarge},~\bfnm{Marc}\binits{M.}} \AND
  \bauthor{\bsnm{Massouli\'{e}},~\bfnm{Laurent}\binits{L.}}
(\byear{2012}).
\btitle{Community detection in the labelled stochastic block model}.
In \bbooktitle{NIPS Workshop on Algorithmic and Statistical Approaches for
  Large Social Networks}.
\end{binproceedings}
\endbibitem

\bibitem{Holland_Laskey_Leinhardt_1983}
\begin{barticle}[author]
\bauthor{\bsnm{Holland},~\bfnm{Paul~W.}\binits{P.~W.}},
  \bauthor{\bsnm{Laskey},~\bfnm{Kathryn~Blackmond}\binits{K.~B.}} \AND
  \bauthor{\bsnm{Leinhardt},~\bfnm{Samuel}\binits{S.}}
(\byear{1983}).
\btitle{Stochastic blockmodels: {F}irst steps}.
\bjournal{Social Networks}
\bvolume{5}
\bpages{109--137}.
\bdoi{10.1016/0378-8733(83)90021-7}
\end{barticle}
\endbibitem

\bibitem{Holme_Saramaki_2012}
\begin{barticle}[author]
\bauthor{\bsnm{Holme},~\bfnm{Petter}\binits{P.}} \AND
  \bauthor{\bsnm{Saram{\"a}ki},~\bfnm{Jari}\binits{J.}}
(\byear{2012}).
\btitle{Temporal networks}.
\bjournal{Physics Reports}
\bvolume{519}
\bpages{97--125}.
\end{barticle}
\endbibitem

\bibitem{Janson_Luczak_Rucinski_2000}
\begin{bbook}[author]
\bauthor{\bsnm{Janson},~\bfnm{Svante}\binits{S.}},
  \bauthor{\bsnm{{\L}uczak},~\bfnm{Tomasz}\binits{T.}} \AND
  \bauthor{\bsnm{Ruci\'{n}ski},~\bfnm{Andrzej}\binits{A.}}
(\byear{2000}).
\btitle{Random Graphs}.
\bpublisher{Wiley}.
\bdoi{10.1002/9781118032718}
\end{bbook}
\endbibitem

\bibitem{Jog_Loh_2015}
\begin{binproceedings}[author]
\bauthor{\bsnm{Jog},~\bfnm{Varun}\binits{V.}} \AND
  \bauthor{\bsnm{Loh},~\bfnm{Po-Ling}\binits{P.-L.}}
(\byear{2015}).
\btitle{Recovering communities in weighted stochastic block models}.
In \bbooktitle{2015 53rd Annual Allerton Conference on Communication, Control,
  and Computing (Allerton)}
\bpages{1308--1315}.
\bdoi{doi: 10.1109/ALLERTON.2015.7447159}
\end{binproceedings}
\endbibitem

\bibitem{Kivela_etal_2014}
\begin{barticle}[author]
\bauthor{\bsnm{Kivel{\"a}},~\bfnm{Mikko}\binits{M.}},
  \bauthor{\bsnm{Arenas},~\bfnm{Alex}\binits{A.}},
  \bauthor{\bsnm{Barthelemy},~\bfnm{Marc}\binits{M.}},
  \bauthor{\bsnm{Gleeson},~\bfnm{James~P.}\binits{J.~P.}},
  \bauthor{\bsnm{Moreno},~\bfnm{Yamir}\binits{Y.}} \AND
  \bauthor{\bsnm{Porter},~\bfnm{Mason~A.}\binits{M.~A.}}
(\byear{2014}).
\btitle{{Multilayer networks}}.
\bjournal{Journal of Complex Networks}
\bvolume{2}
\bpages{203-271}.
\bdoi{10.1093/comnet/cnu016}
\end{barticle}
\endbibitem

\bibitem{Lei_Lin_2022}
\begin{barticle}[author]
\bauthor{\bsnm{Lei},~\bfnm{Jing}\binits{J.}} \AND
  \bauthor{\bsnm{Lin},~\bfnm{Kevin~Z}\binits{K.~Z.}}
(\byear{2022}).
\btitle{Bias-adjusted spectral clustering in multi-layer stochastic block
  models}.
\bjournal{Journal of the American Statistical Association, to appear}.
\end{barticle}
\endbibitem

\bibitem{Lei_Rinaldo_2015}
\begin{barticle}[author]
\bauthor{\bsnm{Lei},~\bfnm{Jing}\binits{J.}} \AND
  \bauthor{\bsnm{Rinaldo},~\bfnm{Alessandro}\binits{A.}}
(\byear{2015}).
\btitle{Consistency of spectral clustering in stochastic block models}.
\bjournal{Annals of Statistics}
\bvolume{43}
\bpages{215--237}.
\bdoi{10.1214/14-AOS1274}
\end{barticle}
\endbibitem

\bibitem{Lei_etal_2017}
\begin{barticle}[author]
\bauthor{\bsnm{Lei},~\bfnm{Yang}\binits{Y.}},
  \bauthor{\bsnm{Bezdek},~\bfnm{James~C.}\binits{J.~C.}},
  \bauthor{\bsnm{Romano},~\bfnm{Simone}\binits{S.}},
  \bauthor{\bsnm{Vinh},~\bfnm{Nguyen~Xuan}\binits{N.~X.}},
  \bauthor{\bsnm{Chan},~\bfnm{Jeffrey}\binits{J.}} \AND
  \bauthor{\bsnm{Bailey},~\bfnm{James}\binits{J.}}
(\byear{2017}).
\btitle{Ground truth bias in external cluster validity indices}.
\bjournal{Pattern Recognition}
\bvolume{65}
\bpages{58-70}.
\bdoi{https://doi.org/10.1016/j.patcog.2016.12.003}
\end{barticle}
\endbibitem

\bibitem{Lelarge_Massoulie_Xu_2015}
\begin{barticle}[author]
\bauthor{\bsnm{Lelarge},~\bfnm{Marc}\binits{M.}},
  \bauthor{\bsnm{Massouli{\'e}},~\bfnm{Laurent}\binits{L.}} \AND
  \bauthor{\bsnm{Xu},~\bfnm{Jiaming}\binits{J.}}
(\byear{2015}).
\btitle{Reconstruction in the labelled stochastic block model}.
\bjournal{IEEE Trans. Netw. Sci. Eng.}
\bvolume{2}
\bpages{152--163}.
\bdoi{10.1109/TNSE.2015.2490580}
\end{barticle}
\endbibitem

\bibitem{Levin_Peres_Wilmer_2008}
\begin{bbook}[author]
\bauthor{\bsnm{Levin},~\bfnm{David~A.}\binits{D.~A.}},
  \bauthor{\bsnm{Peres},~\bfnm{Yuval}\binits{Y.}} \AND
  \bauthor{\bsnm{Wilmer},~\bfnm{Elizabeth~L.}\binits{E.~L.}}
(\byear{2008}).
\btitle{Markov Chains and Mixing Times}.
\bpublisher{American Mathematical Society},
  \baddress{\url{http://pages.uoregon.edu/dlevin/MARKOV/}}.
\end{bbook}
\endbibitem

\bibitem{Lewis_Gonzalez_Kaufman_2012}
\begin{barticle}[author]
\bauthor{\bsnm{Lewis},~\bfnm{Kevin}\binits{K.}},
  \bauthor{\bsnm{Gonzalez},~\bfnm{Marco}\binits{M.}} \AND
  \bauthor{\bsnm{Kaufman},~\bfnm{Jason}\binits{J.}}
(\byear{2012}).
\btitle{Social selection and peer influence in an online social network}.
\bjournal{Proceedings of the National Academy of Sciences}
\bvolume{109}
\bpages{68--72}.
\end{barticle}
\endbibitem

\bibitem{Longepierre_Matias_2019}
\begin{barticle}[author]
\bauthor{\bsnm{Longepierre},~\bfnm{L\'{e}a}\binits{L.}} \AND
  \bauthor{\bsnm{Matias},~\bfnm{Catherine}\binits{C.}}
(\byear{2019}).
\btitle{Consistency of the maximum likelihood and variational estimators in a
  dynamic stochastic block model}.
\bjournal{Electronic Journal of Statistics}
\bvolume{13}
\bpages{4157--4223}.
\bdoi{10.1214/19-EJS1624}
\end{barticle}
\endbibitem

\bibitem{Massoulie_2014}
\begin{binproceedings}[author]
\bauthor{\bsnm{Massouli{\'e}},~\bfnm{Laurent}\binits{L.}}
(\byear{2014}).
\btitle{Community detection thresholds and the weak Ramanujan property}.
In \bbooktitle{Proc. 46th annual ACM Symposium on Theory of Computing}
\bpages{694--703}.
\end{binproceedings}
\endbibitem

\bibitem{Mastrandrea_Fournet_Barrat_2015}
\begin{barticle}[author]
\bauthor{\bsnm{Mastrandrea},~\bfnm{Rossana}\binits{R.}},
  \bauthor{\bsnm{Fournet},~\bfnm{Julie}\binits{J.}} \AND
  \bauthor{\bsnm{Barrat},~\bfnm{Alain}\binits{A.}}
(\byear{2015}).
\btitle{Contact patterns in a high school: A comparison between data collected
  using wearable Sensors, contact diaries and friendship surveys}.
\bjournal{PLOS ONE}
\bvolume{10}
\bpages{1--26}.
\end{barticle}
\endbibitem

\bibitem{Matias_Miele_2017}
\begin{barticle}[author]
\bauthor{\bsnm{Matias},~\bfnm{Catherine}\binits{C.}} \AND
  \bauthor{\bsnm{Miele},~\bfnm{Vincent}\binits{V.}}
(\byear{2017}).
\btitle{Statistical clustering of temporal networks through a dynamic
  stochastic block model}.
\bjournal{J. R. Stat. Soc. Ser. B. Stat. Methodol.}
\bvolume{79}
\bpages{1119--1141}.
\bdoi{10.1111/rssb.12200}
\end{barticle}
\endbibitem

\bibitem{Mazzarisi_Barucca_Lillo_Tantari_2020}
\begin{barticle}[author]
\bauthor{\bsnm{Mazzarisi},~\bfnm{P.}\binits{P.}},
  \bauthor{\bsnm{Barucca},~\bfnm{P.}\binits{P.}},
  \bauthor{\bsnm{Lillo},~\bfnm{F.}\binits{F.}} \AND
  \bauthor{\bsnm{Tantari},~\bfnm{D.}\binits{D.}}
(\byear{2020}).
\btitle{A dynamic network model with persistent links and node-specific latent
  variables, with an application to the interbank market}.
\bjournal{European Journal of Operational Research}
\bvolume{281}
\bpages{50--65}.
\bdoi{https://doi.org/10.1016/j.ejor.2019.07.024}
\end{barticle}
\endbibitem

\bibitem{Meila_2007}
\begin{barticle}[author]
\bauthor{\bsnm{Meil{\u a}},~\bfnm{Marina}\binits{M.}}
(\byear{2007}).
\btitle{Comparing clusterings -- An information based distance}.
\bjournal{Journal of Multivariate Analysis}
\bvolume{98}
\bpages{873-895}.
\bdoi{https://doi.org/10.1016/j.jmva.2006.11.013}
\end{barticle}
\endbibitem

\bibitem{Meila_Heckerman_2001}
\begin{barticle}[author]
\bauthor{\bsnm{Meil{\u a}},~\bfnm{Marina}\binits{M.}} \AND
  \bauthor{\bsnm{Heckerman},~\bfnm{David}\binits{D.}}
(\byear{2001}).
\btitle{An experimental comparison of model-based clustering methods}.
\bjournal{Machine Learning}
\bvolume{42}
\bpages{9--29}.
\bdoi{10.1023/A:1007648401407}
\end{barticle}
\endbibitem

\bibitem{Mezard_Montanari_2009}
\begin{bbook}[author]
\bauthor{\bsnm{Mezard},~\bfnm{Marc}\binits{M.}} \AND
  \bauthor{\bsnm{Montanari},~\bfnm{Andrea}\binits{A.}}
(\byear{2009}).
\btitle{Information, Physics, and Computation}.
\bpublisher{Oxford University Press}.
\end{bbook}
\endbibitem

\bibitem{Moore_2017}
\begin{barticle}[author]
\bauthor{\bsnm{Moore},~\bfnm{Cristopher}\binits{C.}}
(\byear{2017}).
\btitle{The computer science and physics of community detection: {L}andscapes,
  phase transitions, and hardness}.
\bjournal{Bulletin of the {EATCS}}
\bvolume{121}.
\end{barticle}
\endbibitem

\bibitem{Mossel_Neeman_Sly_2016}
\begin{binproceedings}[author]
\bauthor{\bsnm{Mossel},~\bfnm{Elchanan}\binits{E.}},
  \bauthor{\bsnm{Neeman},~\bfnm{Joe}\binits{J.}} \AND
  \bauthor{\bsnm{Sly},~\bfnm{Allan}\binits{A.}}
(\byear{2015}).
\btitle{Consistency thresholds for the planted bisection model}.
In \bbooktitle{Proc. 47th annual ACM Symposium on Theory of Computing}
\bpages{69--75}.
\end{binproceedings}
\endbibitem

\bibitem{Mossel_Neeman_Sly_2015}
\begin{barticle}[author]
\bauthor{\bsnm{Mossel},~\bfnm{Elchanan}\binits{E.}},
  \bauthor{\bsnm{Neeman},~\bfnm{Joe}\binits{J.}} \AND
  \bauthor{\bsnm{Sly},~\bfnm{Allan}\binits{A.}}
(\byear{2015}).
\btitle{Reconstruction and estimation in the planted partition model}.
\bjournal{Probability Theory and Related Fields}.
\end{barticle}
\endbibitem

\bibitem{Mossel_Neeman_Sly_2018}
\begin{barticle}[author]
\bauthor{\bsnm{Mossel},~\bfnm{Elchanan}\binits{E.}},
  \bauthor{\bsnm{Neeman},~\bfnm{Joe}\binits{J.}} \AND
  \bauthor{\bsnm{Sly},~\bfnm{Allan}\binits{A.}}
(\byear{2018}).
\btitle{A proof of the block model threshold conjecture}.
\bjournal{Combinatorica}
\bvolume{38}
\bpages{665--708}.
\end{barticle}
\endbibitem

\bibitem{Paul_Chen_2016}
\begin{barticle}[author]
\bauthor{\bsnm{Paul},~\bfnm{Subhadeep}\binits{S.}} \AND
  \bauthor{\bsnm{Chen},~\bfnm{Yuguo}\binits{Y.}}
(\byear{2016}).
\btitle{Consistent community detection in multi-relational data through
  restricted multi-layer stochastic blockmodel}.
\bjournal{Electronic Journal of Statistics}
\bvolume{10}
\bpages{3807--3870}.
\bdoi{10.1214/16-EJS1211}
\end{barticle}
\endbibitem

\bibitem{Paul_Chen_2020}
\begin{barticle}[author]
\bauthor{\bsnm{Paul},~\bfnm{Subhadeep}\binits{S.}} \AND
  \bauthor{\bsnm{Chen},~\bfnm{Yuguo}\binits{Y.}}
(\byear{2020}).
\btitle{Spectral and matrix factorization methods for consistent community
  detection in multi-layer networks}.
\bjournal{Annals of Statistics}
\bvolume{48}
\bpages{230--250}.
\bdoi{10.1214/18-AOS1800}
\end{barticle}
\endbibitem

\bibitem{Peixoto_2019}
\begin{bincollection}[author]
\bauthor{\bsnm{Peixoto},~\bfnm{Tiago~P.}\binits{T.~P.}}
(\byear{2019}).
\btitle{Bayesian stochastic blockmodeling}.
In \bbooktitle{Advances in Network Clustering and Blockmodeling}
(\beditor{\bfnm{P.}\binits{P.}~\bsnm{Doreian}},
  \beditor{\bfnm{V.}\binits{V.}~\bsnm{Batagelj}} \AND
  \beditor{\bfnm{A.}\binits{A.}~\bsnm{Ferligoj}}, eds.)
\bchapter{11},
\bpages{289--332}.
\bpublisher{John Wiley \& Sons Ltd}.
\bdoi{https://doi.org/10.1002/9781119483298.ch11}
\end{bincollection}
\endbibitem

\bibitem{Pensky_2019}
\begin{barticle}[author]
\bauthor{\bsnm{Pensky},~\bfnm{Marianna}\binits{M.}}
(\byear{2019}).
\btitle{Dynamic network models and graphon estimation}.
\bjournal{Annals of Statistics}
\bvolume{47}
\bpages{2378--2403}.
\bdoi{10.1214/18-AOS1751}
\end{barticle}
\endbibitem

\bibitem{Pensky_Zhang_2019}
\begin{barticle}[author]
\bauthor{\bsnm{Pensky},~\bfnm{Marianna}\binits{M.}} \AND
  \bauthor{\bsnm{Zhang},~\bfnm{Teng}\binits{T.}}
(\byear{2019}).
\btitle{Spectral clustering in the dynamic stochastic block model}.
\bjournal{Electronic Journal of Statistics}
\bvolume{13}
\bpages{678--709}.
\bdoi{10.1214/19-EJS1533}
\end{barticle}
\endbibitem

\bibitem{Rastelli_Fop_2020}
\begin{barticle}[author]
\bauthor{\bsnm{Rastelli},~\bfnm{Riccardo}\binits{R.}} \AND
  \bauthor{\bsnm{Fop},~\bfnm{Michael}\binits{M.}}
(\byear{2020}).
\btitle{A stochastic block model for interaction lengths}.
\bjournal{Advances in Data Analysis and Classification}
\bvolume{14}
\bpages{485--512}.
\bdoi{10.1007/s11634-020-00403-w}
\end{barticle}
\endbibitem

\bibitem{Suveges_Olhede_2022}
\begin{bmisc}[author]
\bauthor{\bsnm{S{\"u}veges},~\bfnm{Maria}\binits{M.}} \AND
  \bauthor{\bsnm{Olhede},~\bfnm{Sofia~C.}\binits{S.~C.}}
(\byear{2022}).
\btitle{Networks with correlated edge processes}.
\bnote{\url{https://arxiv.org/abs/2207.02545}}.
\bdoi{10.48550/ARXIV.2207.02545}
\end{bmisc}
\endbibitem

\bibitem{vanErven_Harremoes_2014}
\begin{barticle}[author]
\bauthor{\bsnm{{van Erven}},~\bfnm{T.}\binits{T.}} \AND
  \bauthor{\bsnm{{Harremo\"es}},~\bfnm{P.}\binits{P.}}
(\byear{2014}).
\btitle{R{\'e}nyi divergence and Kullback--Leibler divergence}.
\bjournal{IEEE Transactions on Information Theory}
\bvolume{60}
\bpages{3797-3820}.
\bdoi{10.1109/TIT.2014.2320500}
\end{barticle}
\endbibitem

\bibitem{Wang_Bickel_2017}
\begin{barticle}[author]
\bauthor{\bsnm{Wang},~\bfnm{Y.~X.~Rachel}\binits{Y.~X.~R.}} \AND
  \bauthor{\bsnm{Bickel},~\bfnm{Peter~J.}\binits{P.~J.}}
(\byear{2017}).
\btitle{Likelihood-based model selection for stochastic block models}.
\bjournal{Annals of Statistics}
\bvolume{45}
\bpages{500--528}.
\bdoi{10.1214/16-AOS1457}
\end{barticle}
\endbibitem

\bibitem{Xu_Hero_2014}
\begin{barticle}[author]
\bauthor{\bsnm{Xu},~\bfnm{Kevin~S}\binits{K.~S.}} \AND
  \bauthor{\bsnm{Hero},~\bfnm{Alfred~O}\binits{A.~O.}}
(\byear{2014}).
\btitle{Dynamic stochastic blockmodels for time-evolving social networks}.
\bjournal{IEEE Journal of Selected Topics in Signal Processing}
\bvolume{8}
\bpages{552--562}.
\end{barticle}
\endbibitem

\bibitem{Xu_Jog_Loh_2020}
\begin{barticle}[author]
\bauthor{\bsnm{Xu},~\bfnm{Min}\binits{M.}},
  \bauthor{\bsnm{Jog},~\bfnm{Varun}\binits{V.}} \AND
  \bauthor{\bsnm{Loh},~\bfnm{Po-Ling}\binits{P.-L.}}
(\byear{2020}).
\btitle{Optimal rates for community estimation in the weighted stochastic block
  model}.
\bjournal{Annals of Statistics}
\bvolume{48}
\bpages{183--204}.
\end{barticle}
\endbibitem

\bibitem{Yang_Chi_Zhu_Gong_Jin_2011}
\begin{barticle}[author]
\bauthor{\bsnm{Yang},~\bfnm{Tianbao}\binits{T.}},
  \bauthor{\bsnm{Chi},~\bfnm{Yun}\binits{Y.}},
  \bauthor{\bsnm{Zhu},~\bfnm{Shenghuo}\binits{S.}},
  \bauthor{\bsnm{Gong},~\bfnm{Yihong}\binits{Y.}} \AND
  \bauthor{\bsnm{Jin},~\bfnm{Rong}\binits{R.}}
(\byear{2011}).
\btitle{Detecting communities and their evolutions in dynamic social
  networks---a Bayesian approach}.
\bjournal{Machine Learning}
\bvolume{82}
\bpages{157--189}.
\bdoi{10.1007/s10994-010-5214-7}
\end{barticle}
\endbibitem

\bibitem{Yun_Proutiere_2016}
\begin{binproceedings}[author]
\bauthor{\bsnm{Yun},~\bfnm{Se-Young}\binits{S.-Y.}} \AND
  \bauthor{\bsnm{Prouti{\`e}re},~\bfnm{Alexandre}\binits{A.}}
(\byear{2016}).
\btitle{Optimal cluster recovery in the labeled stochastic block model}.
In \bbooktitle{Proc. 30th International Conference on Neural Information
  Processing Systems}
\bpages{973--981}.
\bpublisher{Curran Associates Inc.}, \baddress{USA}.
\end{binproceedings}
\endbibitem

\bibitem{Zhang_Zhou_2016}
\begin{barticle}[author]
\bauthor{\bsnm{Zhang},~\bfnm{Anderson~Y.}\binits{A.~Y.}} \AND
  \bauthor{\bsnm{Zhou},~\bfnm{Harrison~Huibin}\binits{H.~H.}}
(\byear{2016}).
\btitle{Minimax rates of community detection in stochastic block models}.
\bjournal{Annals of Statistics}
\bvolume{44}
\bpages{2252--2280}.
\end{barticle}
\endbibitem

\bibitem{Zhao_Wang_Li_Wang_Wang_Gao_2014}
\begin{barticle}[author]
\bauthor{\bsnm{Zhao},~\bfnm{Dawei}\binits{D.}},
  \bauthor{\bsnm{Wang},~\bfnm{Lianhai}\binits{L.}},
  \bauthor{\bsnm{Li},~\bfnm{Shudong}\binits{S.}},
  \bauthor{\bsnm{Wang},~\bfnm{Zhen}\binits{Z.}},
  \bauthor{\bsnm{Wang},~\bfnm{Lin}\binits{L.}} \AND
  \bauthor{\bsnm{Gao},~\bfnm{Bo}\binits{B.}}
(\byear{2014}).
\btitle{Immunization of epidemics in multiplex networks}.
\bjournal{PLOS ONE}
\bvolume{9}
\bpages{e112018}.
\end{barticle}
\endbibitem

\end{thebibliography}
